\documentclass[paper=a4,DIV=10,10pt]{scrartcl}

\areaset{18cm}{26cm}

\usepackage[utf8]{inputenc}
\usepackage{amsfonts,amssymb,amsmath}
\usepackage{mathtools}
\usepackage{booktabs}
\usepackage{cite}
\usepackage{float}
\usepackage{tikz}

\usetikzlibrary{intersections}

\usepackage{bm}

\usepackage[bookmarksnumbered=true,hidelinks,psdextra]{hyperref}

\usepackage[normalem]{ulem}

\makeatletter
\g@addto@macro\bfseries{\boldmath}
\makeatother

\newcommand{\QGL}[1]{Q_n^{\operatorname{GL}}( #1 )}

\newcommand{\ud}{\,\mathrm{d}}
\newcommand{\R}{\mathbb{R}}
\newcommand{\N}{\mathbb{N}}

\newcommand{\buhn}{\overline{\vec u}_h^{\,n}}
\newcommand{\bphn}{\overline{p}_h^{\,n}}

\renewcommand{\vec}{\boldsymbol}

\numberwithin{equation}{section}
\numberwithin{figure}{section}
\numberwithin{table}{section}

\newcommand{\FS}[1]{{}{\color{black} #1}{}}  

\renewcommand{\P}{\mathbb{P}}
\newcommand{\dt}{\partial_t}

\newcommand{\unbh}{\overline{\vec u}^n_h}
\newcommand{\pnbh}{\overline p^n_h}
\newcommand{\tnb}{\overline t_n}
\newcommand{\tnbm}{\overline t_{n-1}}
\newcommand{\uth}{\vec u_{\tau,h}}
\newcommand{\utht}{\vec{\widetilde u}_{\tau,h}}
\newcommand{\Vhdiv}{\vec V_h^{\operatorname{div}}}
\newcommand{\vth}{\vec v_{\tau,h}}
\newcommand{\pth}{p_{\tau,h}}
\newcommand{\qth}{q_{\tau,h}}
\newcommand{\ptht}{{\widetilde p}_{\tau,h}}
\newcommand{\pthb}{{\overline p}_{\tau,h}}
\newcommand{\epb}{{\overline e}^{p}}
\newcommand{\eptb}{{\overline e}^{p}_{\tau,h}}
\newcommand{\sprod}[1]{\left\langle #1 \right\rangle}  
\newcommand{\udth}{\vec{\widehat u}_{\tau,h}}
\newcommand{\pdth}{\widehat{p}_{\tau,h}}
\newcommand{\euth}{ \vec e_{\tau,h}^{\vec u} }
\newcommand{\epth}{ e_{\tau,h}^{p} }
\newcommand{\etau}{ \vec \eta^{\vec u} }
\newcommand{\eu}{ \vec e_{\vec u} }

\newcommand{\heuh}{\vec{\widehat e}_{h}^{\, \vec u} }
\newcommand{\heph}{\vec{\widehat e}_{h}^{\, p} }

\newtheorem{defi}{Definition}[section]
\newtheorem{thm}[defi]{Theorem}
\newtheorem{lem}[defi]{Lemma}
\newtheorem{rem}[defi]{Remark}
\newtheorem{cor}[defi]{Corollary}
\newtheorem{ass}[defi]{Assumption}
\newtheorem{prob}[defi]{Problem}

\newenvironment{mproof}{\paragraph{Proof.}}{\hfill$\blacksquare$}


\begin{document}

\title{Optimal pressure approximation for the nonstationary Stokes problem by a
variational method in time with post-processing }

\author{
Mathias Anselmann${}^1$,
Markus Bause${}^1$\thanks{bause@hsu-hh.de (corresponding author)},
Gunar Matthies${}^2$,
Friedhelm Schieweck${}^3$\\
{\small ${}^1$ Helmut Schmidt University, Faculty of 
Mechanical and Civil Engineering, Holstenhofweg 85,}\\[-1ex]
{\small 22043 Hamburg, Germany}\\
{\small ${}^2$ Technische Universität Dresden, Institut für Numerische
Mathematik,}\\[-1ex]
{\small 01062 Dresden, Germany}\\
{\small ${}^3$ University of Magdeburg, Faculty of Mathematics,
Universitätsplatz 2,}\\[-1ex]
{\small 39016 Magdeburg, Germany}\\
}


\maketitle

\begin{abstract}
\textbf{Abstract.}
We provide an error analysis for the solution of the nonstationary
Stokes problem by a variational method in space and time. 
We use finite elements of higher order for the approximation in space
and a Galerkin--Petrov method with first order polynomials
for the approximation in time. We require global continuity of the
discrete velocity trajectory in time, while allowing the discrete pressure
trajectory to be discontinuous at the endpoints of the time intervals. We show existence and uniqueness of the discrete velocity solution, characterize the set of all discrete pressure solutions and prove an optimal second order estimate in time for the pressure error in the  midpoints of the time intervals.
  The key result {and innovation} is the construction of approximations {to} the pressure
  trajectory by means of post-processing together with the proof of
  optimal order error estimates. 
We propose two variants for a post-processed pressure within the set of pressure solutions {based on collocation techniques or interpolation.
Both variants}
guarantee that the pressure error measured in the $L^2$-norm converges with
optimal second order in time and optimal order in space. 
For the discrete velocity solution, we prove  error estimates
of optimal order in time and space. We present some numerical tests to
support our theoretical results.
\end{abstract}



\section{Introduction}
\label{Sec:Introduction}

In this work we study the numerical approximation of the nonstationary Stokes system
\begin{subequations}
\label{Eq:SE}
\begin{alignat}{3}
\label{Eq:SE_1}
\partial_t \vec u - \Delta \vec u + \nabla p & = \vec f && \quad \text{in } \;
\Omega 
\times (0,T]\,,\\[1ex]
\label{Eq:SE_2}
 \nabla \cdot \vec u & = 0 && \quad \text{in } \; \Omega \times (0,T]\,,\\[1ex]
\label{Eq:SE_3}
 \vec u (0) & = \vec u_0 && \quad \text{in } \; \Omega \,,\\[1ex]
\label{Eq:SE_4}
  \vec u & = \vec 0 && \quad \text{on } \; \partial\Omega \times [0,T]\,,
\end{alignat}
\end{subequations}
by a variational method in time and space. 
In \eqref{Eq:SE}, let $\Omega \subset \R^d$, $d=2,3$, be an open,
bounded domain and $T>0$ the final time. Further, we denote by
$\vec u:\overline\Omega\times[0,T]\to\R^d$ and $p:\overline\Omega\times[0,T]\to\R$ the unknown velocity and pressure field, respectively. 
The force $\vec f$ and the initial velocity $\vec u_0$ are prescribed data. 
In \eqref{Eq:SE_1}, the fluid's viscosity $\nu$ is dropped and chosen as $\nu=1$. 
Throughout, we assume that the Stokes system \eqref{Eq:SE} admits a
sufficiently regular solution up to $t=0$. For
the discussion of the existence and regularity of solutions to
Stokes and Navier--Stokes systems, we refer to the literature (cf., e.g.,
\cite{HR82,HR88,J16,T84}). The behavior of Stokes solutions at time $t=0$
with nonlocal compatibility conditions for the data is studied in
\cite{HR82,R83,T82}.

In order to solve the problem~\eqref{Eq:SE} numerically we apply
a variational method with respect to space and time.
To approximate $\vec u$ and $p$ in space, we use the finite element spaces
$\vec V_h$ and $Q_h$, respectively, which have an arbitrary approximation
order $r$ and satisfy the usual inf-sup stability condition.
For thex
approximation with respect to time, we restrict ourselves to the
simplest case of piecewise first order polynomials. 
An extension of our ideas to higher order polynomials in time is
feasible and left for future work. 
For the approximation in time, we
decompose the global time interval $I=(0,T]$ into local
time intervals $I_n = (t_{n-1},t_n]$, $n=1,\dots, N$, where
$0=t_0<t_1< \dots < t_N=T$, and use polynomial spaces
$\P_k(I_n;B)$ of functions $w: I_n \to B$, where $w(t)$ is a
polynomial in $t$ of degree not greater than $k$ with values in a Banach
space $B$.
%
%
Then, we approximate the exact velocity $\vec u(t)$
for $t\in I_n$ by means of a polynomial $\vec u_{\tau,h}(t)$ of first
order in $t$ with values in the discrete velocity space $\vec V_h$,
i.e., we use the ansatz  
  $\vec u_{\tau,h}{}_{| I_n} \in \P_1(I_n; \vec V_h)$.
Since we want to have that the time derivative $\partial_t\vec u_{\tau,h}$
exists globally on $[0,T]$ in the weak $L^2$-sense, we require
that the fully discrete solution trajectory $\vec u_{\tau,h}(t)$ is globally
continuous on $[0,T]$. 
This global continuity will be realized during a successive time-marching process 
over the time intervals by means 
of an explicit condition on each interval $I_n = (t_{n-1},t_n]$ which can be
regarded as a local initial condition for $\vec u_{\tau,h}{}_{|  I_n}$ at
$t=t_{n-1}$.
%
%
In order to get for the pressure the same order of approximation with
respect to time as for the velocity, we approximate the exact pressure $p(t)$ 
for $t\in I_n$ by the ansatz  
  $p_{\tau,h}{}_{| I_n} \in \P_1(I_n; Q_h)$.
However, in contrast to the velocity, we allow that the global discrete pressure
trajectory $\pth(t)$ can be discontinuous at the endpoints of the time
intervals.
To the best of our knowledge, this approach concerning the approximation of the
pressure is new.

The fully discrete approximations $\uth(t)$ and $\pth(t)$ for $t\in I_n$,
according to the ansatz described above,
are determined by means of discrete variational equations.
These variational equations are consistent in the sense that 
they are satisfied by the exact solution $\vec u(t)$ and $p(t)$ as well. 
In the following, we will derive, for time interval $I_n$, the
variational equation for $\vec u(t)$ and $p(t)$  
which is related to~\eqref{Eq:SE_1}.
%
We choose the test function $\vth$ as polynomial in time with values in
$\vec V_h$, where the polynomial degree $k$ will be specified later. Then, for arbitrary fixed $t\in I_n$, we take on both sides
of~\eqref{Eq:SE_1} 
the inner product $\sprod{\cdot,\cdot}$ 
in $L^2(\Omega)^d$ with the test function in space $\vth(t)\in\vec V_h$ and obtain
\begin{equation}
\label{Eq:var_eq_space}
   \sprod{\partial_t \vec u(t), \vth(t)} + \sprod{\nabla\vec u(t),
   \nabla\vth(t)} -
   \sprod{p(t), \nabla\cdot\vth(t)}
  =
   \sprod{\vec f(t), \vth(t)}  \,, 
\end{equation}
where we have applied integration by parts to the terms related to 
$-\Delta \boldsymbol{u}$ and  $\nabla p$ using the homogeneous boundary condition
of $\vth(t)\in \Vec V_h$.
{We note, that for functions in $\Vec V_h$ we require homogeneous
boundary conditions.}
Since~\eqref{Eq:var_eq_space} holds for all $t\in I_n$, we can apply
numerical integration with respect to time over the interval $I_n$ and 
get for all $\vec v_{\tau,h} \in \P_k(I_n; \vec V_h)$ that
\begin{equation}
\label{Eq:var_eq}
  \QGL{
   \sprod{\partial_t \vec u, \vth} + \sprod{\nabla\vec u, \nabla\vth} -
   \sprod{p, \nabla\cdot\vth}
  }
  =
  \QGL{ \sprod{\vec f, \vth} } \,, 
\end{equation}
where $\QGL{g}$ denotes for a time-dependent function $g: I_n \to \R$ the value of the 2-point Gau{ss}--Lobatto quadrature formula approximating the
integral $\int_{I_n} g(t)\ud t$ (see~\eqref{Eq:GLF} below).

Now, we get the variational equation for the discrete velocity and
pressure by replacing in~\eqref{Eq:var_eq} the exact solution $\vec u$ and
$p$ by the discrete solution $\uth{}_{| I_n}$ and $\pth{}_{| I_n}$.
In an analogous way, we can derive from the continuity equation~\eqref{Eq:SE_2}
an associated variational equation using a test function
$\qth \in \P_k(I_n; Q_h)$.
For the polynomial degree $k$, it will turn out that the choice $k=0$ is
suitable.
Thus, setting  $k=0$ in the variational equations derived above,
we obtain the following discrete problem on interval
$I_n= (t_{n-1},t_n]$ to determine the discrete solution $\uth{}_{| I_n}$
and $\pth{}_{| I_n}$:\medskip

\textit{For given $\uth(t_{n-1})$, find 
  $\uth{}_{| I_n} \in \P_1(I_n; \vec V_h)$ and $\pth{}_{| I_n} \in \P_1(I_n; Q_h)$
  such that   }
\begin{subequations}
\label{Eq:dVE}
\begin{alignat}{3}
\label{Eq:dVE1}
  \uth^+(t_{n-1}) & = \uth(t_{n-1})  \,, & 
  \\[1ex]
\label{Eq:dVE2}
\QGL{
   \sprod{\dt\uth , \vth} + \sprod{\nabla\uth , \nabla\vth} - \sprod{\pth, \nabla\cdot\vth}
}  
  & = \QGL{  \sprod{\vec f, \vth}  } 
  & \quad  \forall \vth \in \P_0(I_n; \vec V_h)  \,,
   \\[1ex]
\label{Eq:dVE3}
  \QGL{ \sprod{\nabla\cdot\uth , \qth }  }  & = 0
  & \quad  \forall \qth \in \P_0(I_n; Q_h) \,.
\end{alignat}
\end{subequations}
In~\eqref{Eq:dVE1},
$\uth^+(t_{n-1})$ denotes the right-sided limit of the values
$\uth(t)$ when $t\in I_n$ tends to $t_{n-1}$. 
The value $\uth(t_{n-1})$ is known: for $n>1$ from the previous interval
$I_{n-1}$ since $\uth(t_{n-1})=\uth{}_{| I_{n-1}}(t_{n-1})$ and for $n=1$ from
the discrete initial value, i.e., $\uth(0)=\vec u_{0,h}$, where $\vec
u_{0,h} \in \vec V_h$ is an approximation of the initial value $\vec u_0$.
Thus, condition~\eqref{Eq:dVE1} for all $n=1, \dots, N$ guarantees the
global continuity of discrete velocity trajectory $\uth(t)$ for
$t\in[0,T]$.


In this paper, we study existence and uniqueness of
solutions $(\uth, \pth)$ of the discrete problem~\eqref{Eq:dVE}
related to the interval $I_n$. Moreover, we exactly describe the set
of all solutions (cf.\  Lemma~\ref{Lem:space_time_probl}). We prove that the
velocity solution $\uth$ is uniquely determined. However, we show for
the pressure solution $\pth$ that the $Q_h$-valued degree of freedom
$\pnbh\coloneqq\pth(\tnb)$ with $\tnb \coloneqq (t_{n-1}+t_n)/2$ is
uniquely determined by~\eqref{Eq:dVE}, whereas the degree of freedom
$p^{n-1}_h \coloneqq p^+_{\tau,h}(t_{n-1})$ can be chosen arbitrarily. 
In our error analysis we show that the global continuity of the discrete
velocity trajectory $\uth(t)$ guaranteed by~\eqref{Eq:dVE1}
and the fulfillment of the variational equations~\eqref{Eq:dVE2} and~\eqref{Eq:dVE3}
for all time intervals $I_n$ are sufficient to prove  optimal
estimates for the velocity error 
$\vec e_{\vec u} \coloneqq \vec u - \vec u_{\tau,h}$
of  order $\tau^2 + h^r$ (see   Theorem~\ref{Thm:Erruprime}),
where $\tau$ and $h$ denote the mesh sizes in time and space,
respectively. Precisely, we have
\begin{equation}
\label{Eq:err_est_u_all}
    \sum_{m=1}^n \tau_m \| \partial_t \vec e_{\vec u}(\overline t_m ) \|^2 
    + \sum_{m=1}^n \tau_m \| \nabla \vec e_{\vec u} (\overline t_m) \|^2 
    + \| \nabla \vec e_{\vec u} (t_n) \|^2  
    \le C \big(
      \| \nabla \vec e_{\vec u} (0) \| 
      + \tau^2 + h^{r}
      \big)^2 
\end{equation}
for all $n=1,\ldots , N$,
where $\|\cdot\|$ denotes the $L^2$-norm in space and $\tau_m$ the length of the time 
interval $I_m$.
Note that the initial velocity $\vec u_{0,h}$ can be chosen such that
the initial error norm $\| \nabla \vec e_{\vec u} (0) \|$ is of the
order $h^r$.
It is also important to note that in terms of accuracy in time, the
term $\tau^2$ in the estimate~\eqref{Eq:err_est_u_all} shows a
superconvergence of $\partial_t \vec u_{\tau,h}$ in the midpoints
$\overline t_m$ of the time intervals $I_m$. This fact is an essential key
point to derive an optimal error estimate for the pressure.
We will show in Theorem~\ref{Thm:EstP}, for each
discrete pressure trajectory $\pth(t)$ in the set of pressure solutions,
that the error estimate
\begin{equation}
\label{Eq:err_est_ptnb}
 \sum_{n=1}^N \tau_n \| p(\tnb) - \pth(\tnb)  \|^2 \le C (\tau^2 + h^r)^2 
\end{equation}
is satisfied at the midpoints $\tnb$ of the intervals $I_n$.
A discrete pressure trajectory $\pth$ restricted to an interval  $I_n$ is
fixed by the uniquely determined value $\pth(\tnb)$ and the
freely selectable degree of freedom $p^+_{\tau,h}(t_{n-1})$.
An important goal of this work is to compute, for each interval $I_n$, the
degree of freedom $p^{n-1}_h \coloneqq p^+_{\tau,h}(t_{n-1})$
in such a way that a pressure error estimate similarly 
to~\eqref{Eq:err_est_ptnb} is obtained, where the midpoints $\tnb$ are replaced by the endpoints $t_n$ of the time intervals. 
Then, this also leads to an optimal order error estimate for
the pressure in the $L^2$-norm with respect to space and time.

In our paper, we propose two variants for computing the values
$p^{n-1}_h \coloneqq p^+_{\tau,h}(t_{n-1})$. 
The first variant is based on the idea that the value $p^{n-1}_h$ 
is chosen as the selected pressure value 
$\widetilde p^{n-1}_h \coloneqq \ptht(t_{n-1})$ of a uniquely defined
post-processed pressure trajectory $\ptht(t)$ which is in the set
of the pressure solutions and globally continuous on $[0,T]$.
Together with $\ptht$ we define a post-processed velocity trajectory
$\utht(t)$ which is globally $C^1$-continuous in time and locally in
time characterized by 
  $\utht{}_{| I_n} \in \P_2(I_n; \vec V_h)$.
The connection of the post-processed solution $(\utht,\ptht)$ to the original 
solution $(\uth,\pth)$ on interval $I_n$ is that 
$\utht(t)=\uth(t)$ for $t\in\{ t_{n-1},t_n \}$ and
$\ptht(\tnb) = \pth(\tnb)$.
The remaining degrees of freedom to fix $(\utht,\ptht)$  on $I_n$  are
$\vec a^{n-1}_h \coloneqq \dt \utht(t_{n-1})$ and  $\widetilde p^{n-1}_h \coloneqq \ptht(t_{n-1})$.
They are determined 
such that the pair $(\vec a^{n-1}_h, \widetilde p^{n-1}_h)$
satisfies a collocation condition at  time $t=t_{n-1}$ together 
with the known velocity value $\uth(t_{n-1})$.  
This collocation condition corresponds to  a Stokes-type problem for the
pair $(\vec a^{n-1}_h, \widetilde p^{n-1}_h)$. However, it turns out that such
a collocation problem has to be solved only once for the first interval
$I_1$ at the time $t_0=0$. We will prove that, for all subsequent intervals
$I_n$, $n>1$, the solution $(\vec a^{n-1}_h, \widetilde p^{n-1}_h)$ is equal to
the value $(\dt\utht(t_{n-1}), \ptht(t_{n-1}))$ coming from the previous
interval $I_{n-1}$ such that the post-processing is computationally cheap.
For the resulting globally continuous post-processed pressure trajectory
$\ptht$, we show an optimal order error estimate in the $L^2$-norm with
respect to space and time (see Corollary~\ref{Cor:ErrPPP}).

In the second variant, we use the idea of interpolation based on the
accurate discrete pressure values $\pth(\overline t_m)$ in the midpoints
$\overline t_m$ of the intervals $I_m$ (see
estimate~\eqref{Eq:err_est_ptnb}) to define another post-processed pressure
trajectory $\pthb$, which is also in the set of the pressure solutions and
locally characterized by $\pthb{}_{| I_n} \in \P_1(I_n; Q_h)$. The connection
between this post-processed pressure $\pthb$ and the original
pressure solution $\pth$ on the interval $I_n$ is that $\pthb(\tnb) =
\pth(\tnb)$. The remaining degree of freedom of $\pthb$ on $I_n$ is
the value $p^{n-1}_h \coloneqq \pthb^+(t_{n-1})$, which is determined
by means of linear interpolation (or extrapolation for $n=1$) of two
accurate pressure values $\pth(\overline t_m)$ from suitable intervals
$I_m$ in the neighborhood of $t_{n-1}$. This approach leads to a global
pressure trajectory that is generally discontinuous at the endpoints
of the time intervals. For the post-processed pressure trajectory $\pthb$,
we will prove an optimal order error estimate in the $L^2$-norm with
respect to space and time (cf.\ Theorem~\ref{Thm:EPPI}). 

{To the best of our knowledge, these results for {our} approximation of
the pressure trajectory are new in the literature.} 
In the following, we discuss some links of our work to the existing
literature in this field. The first link is related to the classical
works of Heywood and Rannacher~\cite{HR82,HR86,HR88,HR90} on finite element
approximations of the Navier--Stokes equations. In~\cite{HR90}, an error
analysis is given for a discretization by finite elements in space and the
Crank--Nicolson finite difference scheme in time, which generates, for the
time nodes $t_n$, approximations $\vec u^n_h$ and $p^n_h$ of the exact
solution values $\vec u(t_n)$ and $p(t_n)$, respectively. If we would apply
this discretization to the Stokes equations then the link to our solution
trajectory $(\uth, \pth)$ is that $\vec u^n_h = \uth(t_n)$  and  $p^n_h =
\pth(\tnb)$ hold for all $n$. The Crank--Nicolson value $p^n_h$ can
also be characterized as the value $\pth(t_n)$, where $\pth$ is the
one special pressure trajectory in the set of all pressure solutions that
is piecewise constant on each time interval $I_n$ such that $\pth(t_n) =
\pth(\tnb) = p^n_h$. Together with our pressure error
estimate~\eqref{Eq:err_est_ptnb}, we see that the error of the Crank--Nicolson
scheme $\| p(t_n) - p^n_h \| = \| p(t_n) - \pth(\tnb) \|$ is generally
only of first order with respect to the time step size. 
For this error, a bound of the order $\tau + h^r$ is proved in~\cite{HR90}. 
\FS{
Further, in~\cite{R08} it is pointed out for the nonstationary Navier--Stokes equations that the pressure value $p^n_h$ resulting from the application of the Crank--Nicolson scheme (regarded as an implicit $\theta$-scheme  with $\theta=1/2$) approximates the exact pressure $p$ at time $t=t_{n-1}+\theta\tau_n = \tnb$ and not at $t=t_n$. 
}

A second link to our work is is the work of Sonner and Richter~\cite{SR19},
where the Crank--Nicolson scheme is considered as semidiscretization in
time of the Stokes and Navier--Stokes equations based on a variational
setting. Here, the test space and also the ansatz space for the
semidiscrete pressure $p_{\tau}$ are chosen to be piecewise constant on the
time intervals, whereas the semidiscrete velocity $\vec u_{\tau}$ is
required to be globally continuous and piecewise linear in time. For the
pressure error $\|p(\tnb) - p_{\tau}(\tnb) \|$ in the midpoints $\tnb$ of
the time intervals the order $\tau^2$ is shown.

\FS{
Further links to this work can be found in \cite{HST12, BKRS18}. 
In~\cite{HST12}, the nonstationary Stokes equations are discretized in
time by means of the continuous Galerkin--Petrov method with the
polynomial order $k=1$ and $k=2$. The discretization in space is based
on a finite element method. The difference to our paper is that 
in~\cite{HST12} the numerical integration with respect to time is based
on the $k$-point Gau{ss} formula such that, for $k=1$, this time
discretization differs from the classical Crank--Nicolson scheme.
In~\cite{HST12},  a post-processing of the discrete pressure at the time nodes $t_n$ is proposed which differs from the proposals of this work. A theoretical error analysis is not contained in \cite{HST12}.
In~\cite{BKRS18}, 
the hyperbolic wave problem is discretized in time by the continuous
Galerkin--Petrov method with the Gau{ss}--Lobatto quadrature formula for the
time integration and the finite element method for the space
discretization. 
A post-processing is introduced and analyzed. It lifts the fully discrete approximation in time from a continuous to a continuously
differentiable one.  The Gau{ss}–-Lobatto quadrature in time 
is essential for ensuring the improved convergence behavior. 
}


This paper is organized as follows. In Section~\ref{Sec:NotPrem}, notation,
preliminaries and auxiliary results are introduced. The numerical scheme
along with the error estimate \eqref{Eq:err_est_ptnb} are presented in
Section~\ref{Sec:AppS_ErrGP}. In Section~\ref{Sec:PPCol}, a first variant of
computing a post-processed solution $(\utht, \ptht)$ by means of
collocation is introduced and corresponding error estimates are proved. In
Section~\ref{Sec:PPInt}, a second variant of computing a post-processed
pressure solution $\pthb$  by means of interpolation is presented and a
corresponding pressure error estimate is proved. Finally, in
Section~\ref{Sec:NumExp} we present some numerical experiments which confirm
and illustrate the theoretical results.

\section{Notation and preliminaries}
\label{Sec:NotPrem}

\subsection{Continuous function spaces}

We use standard notation. $H^m(\Omega)$ is the Sobolev space of $L^2(\Omega)$ 
functions with derivatives up to order $m$ in $L^2(\Omega)$, while 
$\langle \cdot,\cdot \rangle$ denotes the inner product in $L^2(\Omega)$
and its vector-valued and matrix-valued versions.
Let $L^2_0(\Omega)\coloneqq \{ q\in 
L^2(\Omega) \mid \int_\Omega q \ud \vec x =0\}$ and  $H^1_0(\Omega)\coloneqq \{u\in 
H^1(\Omega) \mid u=0 \mbox{ on } \partial \Omega\}$. For short, we put 
\[
Q\coloneqq L^2_0(\Omega)\qquad \text{and} \qquad \vec V\coloneqq  \left(H^1_0(\Omega)\right)^d\,.
\]
We note that bold-face letters are used to indicate vector-valued spaces and functions. Further, we define by
\[
\vec V^{\operatorname{div}}\coloneqq  \{\vec u\in \vec V \mid \langle \nabla \cdot \vec 
u, q\rangle = 0 \:\: \forall q\in Q \}
\]
the space of divergence-free functions.
We denote by $\vec V'$ the dual space of $\vec V$ and use the notation
\begin{align*}
\| \cdot \| \coloneqq  \| \cdot\|_{L^2(\Omega)}\,,\qquad 
\| \cdot \|_m \coloneqq  \| \cdot\|_{H^m(\Omega)}, \quad m \in \N,
\end{align*} 
for the norms of the Sobolev spaces, where we do not differ between
scalar-valued, vector-valued, and matrix-valued cases. 
Throughout, the meaning will be obvious from the context.
For a Banach space $B$ and an interval $J\subset [0,T]$, let $L^2(J;B)$,
$C^m(J;B)$, and $C^0(J;B)\coloneqq C(J;B)$ be the Bochner spaces of
$B$-valued functions, equipped with their natural norms.

For the variational setting we define the bilinear forms $a: \vec V \times 
\vec V \rightarrow \R $ and $b:\vec V \times Q \rightarrow \R$ by
\begin{subequations}
\label{Eq:DefBF}
\begin{alignat}{2}
\label{Eq:DefBLFa}
a(\vec u, \vec v) & \coloneqq \langle \nabla \vec u, \nabla \vec v\rangle \,,
&& \quad \vec u, \vec v\in \vec V\,, \\[1ex]
\label{Eq:DefBFb}
b(\vec v,q) & \coloneqq - \langle \nabla \cdot \vec v,q \rangle\,,
&& \quad \vec v\in \vec V\,, \; q \in Q\,. 
\end{alignat}
\end{subequations}

In what follows, for non-negative expressions $a$ and
$b$, the notation $a\lesssim b$ stands for the inequality $a \leq C\, b$ with a generic
constant $C$ that is independent of the sizes of the spatial and temporal
meshes. The value of $C$ can depend on the regularity of the space mesh,
the polynomial degree of the space discretization, and the data
(including $\Omega$).

\subsection{Discrete function spaces and operators}

For the time discretization, we decompose the time interval
$I \coloneqq (0,T]$ into $N$ subintervals 
$I_n \coloneqq (t_{n-1},t_n]$, $n=1,\ldots,N$, where $0=t_0<t_1< \cdots < t_{N-1} < t_N 
= T$ such that $I=\bigcup_{n=1}^N I_n$. We put $\tau \coloneqq  \max_{n=1,\ldots, N} \tau_n$ with 
$\tau_n \coloneqq t_n - t_{n-1}$. Further, the set $\mathcal{M}_\tau \coloneqq  \{I_1,\ldots, I_N\}$
of time intervals is called the time mesh. For a Banach space $B$, an
interval $J\subset [0,T]$, and any $k\in \N_0$,  let 
\begin{equation}
\label{Def:Pk}
  \P_k({J};B) \coloneqq \bigg\{w_\tau : J \to B \::\: w_\tau(t) = \sum_{j=0}^k 
W^j t^j \; \forall t\in J\,, \; W^j \in B\; \forall j \bigg\}
\end{equation}
be the space of $B$-valued polynomials of degree less than or equal to
$k$ on $J$.
For $k\in \N$, we define the space of globally continuous and piecewise
polynomial functions of degree less than or equal to $k$ in time with values in a Banach space $B$ by 
\begin{equation}
\label{Eq:DefXk} 
X^k_\tau (B) \coloneqq  \left\{w_\tau \in C(\overline{I};B) \::\: w_\tau{}_{| I_n} \in
\mathbb P_k(I_n;B)\; \forall I_n\in \mathcal{M}_\tau \right\}\,.
\end{equation}
For $k\in \N_0$, we define the space of piecewise polynomial functions 
of degree less than or equal to $k$ in time with values in a Banach space $B$ by 
\begin{equation}
	\label{Eq:DefYk}
      Y^k_\tau (B) \coloneqq  \left\{{w_\tau : I \to B} \::\:
      w_\tau{}_{| I_n} \in \mathbb 
	P_{k}(I_n;B)\; \forall I_n\in \mathcal M_\tau \right\}\,.
\end{equation}

For any function $w: \overline I\to B$ that is piecewise sufficiently
smooth with respect to the time mesh $\mathcal{M}_{\tau}$, for instance for
$w\in Y^k_\tau (B)$, we define the right-sided and left-sided
limits at a time mesh point $t_n$ by
\begin{equation}
\label{Eq:Defw_In_bdr}
  w^+(t_n) \coloneqq  \lim_{t\to t_n+0} w(t) ,\quad n<N\,,
\qquad\text{and}\qquad
 w^-(t_n) \coloneqq  \lim_{t\to t_n-0} w(t) ,\quad n>0\,.
\end{equation}
Due to the definition of $Y^k_\tau (B)$, a function $w\in Y^k_\tau (B)$
is continuous from the left, i.e.\  $w(t_n)=w^-(t_n)$,
for all points $t_n$, $n=1, \dots, N$.  Moreover, a function $w\in Y^k_\tau (B)$ 
is globally continuous, i.e.\ belongs to the space $X^k_\tau (B)$, if and only if
\begin{equation}
\label{Eq:glob_contin_w}
  w^+(t_n) =  w(t_n) ,\quad n=0, \dots, N-1 \,.
\end{equation}

For the space discretization, let $\{\mathcal{T}_h\}$ be a
shape-regular family of meshes $\mathcal{T}_h$ of the domain $\Omega$
consisting of triangular or quadrilateral mesh cells in two space
dimensions and tetrahedral or hexahedral mesh cells in three
space dimensions with mesh size $h>0$.
By ${\vec V}(K)$ and $Q(K)$, $K\in \mathcal{T}_h$, we
denote the vector-valued and scalar-valued finite element spaces of mapped polynomials
defined by the reference transformation of polynomials on the
reference element with maximum degree uniformly bounded with
respect to the mesh family
$\{\mathcal{T}_h\}$. For convenience, the local spaces ${\vec V}(K)$ and $Q(K)$ are defined implicitly by prescribing their approximation properties in \eqref{Eq:ApFES}. Exemplified, $\vec V(K)$ can be considered to consist of polynomials of maximum degree $r$ and $Q(K)$  of polynomials of maximum degree $r-1$. Then, the finite element spaces to be used
for approximating $\vec V$ and $Q$ are 
\begin{subequations}
\label{Def:VhQh}
\begin{alignat}{2}
\label{Def:Vh}
\vec V_h & \coloneqq  \{\vec v_h \in \vec V \; : \; \vec v_{h}{}_{| K}\in 
\vec V(K) \;\; \text{for all}\; K \in \mathcal{T}_h\}\,, \\[1ex]
\label{Def:Qh}
  Q_h & \coloneqq  \{\vec q_h \in Q \; : \; \vec q_{h}{}_{| K}\in Q(K) \;\; \text{for all}\; K \in \mathcal{T}_h\}\,.
\end{alignat}
\end{subequations}
Further, we {define}
\begin{equation}
	\label{Def:Vdiv}
	\vec V_h^{\text{div}} \coloneqq  \{\vec v_h \in \vec V_h \; : \; b(\vec v_h,q_h) = 0 \; \text{for all } q_h \in Q_h\}
\end{equation}
as space of {the} discretely divergence-free functions.

Now we make further assumptions about the discrete spaces $\vec V_h$ and $Q_h$.  By  imposing conditions on $\vec V_h$ and $Q_h$ about their approximation properties, the polynomial degrees of functions $\vec v_h\in \vec V_h$ and $q_h\in Q_h$, needed to meet these demands, are defined implicitly. 
\begin{ass}
\label{Ass:FES}
The spaces $\vec V_h$ and $Q_h$ satisfy the uniform inf-sup (or LBB) condition
\begin{equation}
\label{Eq:InfSupCod}
0<\beta\le \inf_{q_h \in Q_h\setminus \{0\}}
\sup_{\vec v_h\in \vec V_h\setminus \{\vec 0\}}
\dfrac{b(\vec v_h,q_h)}{\|\vec v_h \|_1 \, \| q_h\|}
\end{equation}
for some constant $\beta$ independent of $h$. For $\vec v\in \vec V\cap
  \vec H^{r+1}(\Omega)$ and $q\in  Q \cap H^r(\Omega)$, there exist
  approximations $i_h \vec v \in \vec V_h$ and $ j_h q \in Q_h$ such that 
\begin{equation}
\label{Eq:ApFES}
\| \vec v - i_h \vec v \| + h \| \nabla (\vec v-i_h \vec v)\| \lesssim h^{r+1} \qquad \text{and} \qquad \| q - j_h q\| \lesssim  h^r\,.
\end{equation}
\end{ass}
Examples of pairs of finite element spaces satisfying
Assumption~\ref{Ass:FES} can be found, e.g., in~\cite[{Ch.~3}]{J16}. In particular, the well-known family of Taylor--Hood pairs of finite element spaces satisfies the assumptions.

Finally, the discrete operators $A_h:\vec V \rightarrow \vec V_h'$, $B_h: \vec 
V \rightarrow Q_h'$, and $B_h': Q \rightarrow \vec V_h'$ are defined by means of
\begin{subequations}
	\label{Def:A_hB_h}
\begin{alignat}{2}
	\label{Def:A_h}
	\langle A_h \vec u , \vec v_h\rangle_{\vec V_h', \vec V_h} & \coloneqq  a(\vec u,\vec v_h)\,, 
	&& \hspace*{0.5cm} \text{for } \vec u\in \vec V\,,\; \vec v_h\in \vec V_h\,,\\[1ex]
	\label{Def:B_h}
	\langle B_h \vec v, q_h  \rangle_{Q_h',Q_h} & \coloneqq  b(\vec v, q_h)\,, && \hspace*{0.5cm}  
	\text{for } 
	\vec v\in \vec V\,,\; q_h \in Q_h\,,\\[1ex]
	\label{Def:B_h^p}
	\langle B_h' q, \vec v_h \rangle_{\vec V_h',\vec V_h} & \coloneqq  b(\vec v_h, q)\,, &&  
	\hspace*{0.5cm} \text{for } 
	\vec v_h\in \vec V_h\,,\; q \in Q\,,
\end{alignat}
\end{subequations}
with 
the bilinear forms $a$ and $b$
introduced in \eqref{Eq:DefBF}. On $\vec V_h$, the operator $B_h$ is thus the discrete (negative) divergence operator. By 
the representation 
theorem of Riesz (cf.\ \cite[Theorem\ B.3]{J16}), the spaces $\vec V_h'$ and $Q_h'$ can be 
identified with $\vec V_h$ and $Q_h$, respectively. Therefore, for the sake
of brevity, we 
skip the indices of the duality pairings. Usually, it holds
 $\nabla \cdot 
\vec v_h \not\in Q_h$ for $\vec v_h \in \vec V_h$. On $\vec V_h$, the kernel of $B_h$ is 
the space $\vec V_h^{\text{div}}$ defined in \eqref{Def:Vdiv}. On $Q_h$, the 
operator $B_h'$ is the discrete gradient operator and, thus, the
adjoint operator of the discrete divergence.

\subsection{Quadrature formulas, interpolation operators and Stokes projection} 
\label{Subsec:QuadInterpol}

In this work, the 2-point Gau{ss}--Lobatto formula (trapezoidal rule)
and 1-point Gauss quadrature formula (midpoint rule) are applied on
each time interval  $I_n=(t_{n-1},t_n]$.  On $I_n$,  the 2-point
Gau{ss}--Lobatto formula is defined by
\begin{equation}
	\label{Eq:GLF}
      Q_n^{\operatorname{GL}}(g) \coloneqq  \frac{\tau_n}{2}
      \big(g^+(t_{n-1})+ g^-(t_n)\big) \approx 
	\int_{I_n} g(t) \ud t\,,
\end{equation}
where the values $g^+(t_{n-1})$ and $g^-(t_{n})$  denote the
corresponding one-sided limits of $g(t)$ from the interior of $I_n$
(cf.\  \eqref{Eq:Defw_In_bdr}). Formula \eqref{Eq:GLF} is exact for all
polynomials in $\mathbb P_{1} (I_n;\R)$.  On $I_n$, with its midpoint 
\begin{equation}
\label{Eq:Midpt}
\overline t_{n}\coloneqq  \frac{t_{n-1}+t_n}{2}\,, 
\end{equation}  
we denote by 
\begin{equation}
	\label{Eq:GF}
	Q_n^{\operatorname G}(g) \coloneqq   \tau_n\, g(\overline t_{n}) \approx \int_{I_n} g(t) \ud t  
\end{equation}
the $1$-point Gau{ss} quadrature formula. Formula \eqref{Eq:GF} is exact for all polynomials in $\mathbb P_{1} (I_n;\R)$. 

Next, we define some interpolation operators with respect to the time
variable. For this, let $B$ be a Banach space, typically $B\in
\{L^2(\Omega),H^1_0(\Omega)\}$. Firstly, we introduce Lagrange
interpolation operators. For a given function $w\in C(\overline{I};B)$,
we define on each time interval $I_n$ with midpoint $\overline t_n$
defined in \eqref{Eq:Midpt} the local Lagrange interpolants
$I_{n,1}^{\operatorname{GL}}w\in  \mathbb P_1({\overline I_n};B)$  and
$I_{n,2}^{\operatorname{GL}}w\in  \mathbb P_2({\overline I_n};B)$ by
\begin{subequations}
	\begin{alignat}{2}
	\label{Eq:LocLagIntOp1}
	I_{n,1}^{\operatorname{GL}} w(s) & = w(s)\,,  && \quad
        \;  s\in \{{t_{n-1}},t_n\}\,,\\[1ex]
	\label{Eq:LocLagIntOp2}
	I_{n,2}^{\operatorname{GL}} w(s) & = w(s)\,,  && \quad
        \;  s\in \{{t_{n-1}},\overline t_n,t_n\}\,.
\end{alignat}
\end{subequations}

The global Lagrange interpolants $I_{1}^{\operatorname{GL}}w \in
X^{1}_\tau(B)$ and $I_{2}^{\operatorname{GL}}w\in  X^{2}_\tau(B)$,
 are  given by 
\begin{equation}
\label{Eq:GloLagIntOp1}
    I_{k}^{\operatorname{GL}} w(t)  \coloneqq  I_{n,k}^{\operatorname{GL}} w(t)  
    \quad t\in I_n\,, \quad n= 1, \ldots, N\,, \quad k=1,2\,,
\end{equation}
along with $I_{1}^{\operatorname{GL}} w(0)\coloneqq w(0)$ and $I_{2}^{\operatorname{GL}} w(0)\coloneqq w(0)$.

For the Lagrange interpolation operators $I_{n,1}^\mathrm{GL}$ and
$I_{n,2}^\mathrm{GL}$, the error bounds
\begin{subequations}
	\begin{alignat}{2}
		\label{Eq:ApLag_0}
		\| v - I_{n,1}^\mathrm{GL}\, v \|_{C(\overline I_n;B)} & \lesssim \tau_n^2 \| v \|_{C^{2}(\overline{I}_n;B)}\,, \\[1ex]
		\label{Eq:ApLag_1}
		\| v - I_{n,2}^\mathrm{GL}\, v \|_{C(\overline I_n;B)} & \lesssim \tau_n^3 \| v \|_{C^{3}(\overline{I}_n;B)}\,, \\[1ex]
		\label{Eq:ApLag_2}
		\| \partial_t v - \partial_t I_{n,2}^\mathrm{GL} \, v \|_{C(\overline I_n;B)} & \lesssim \tau_n^2 \| v \|_{C^{3}(\overline{I}_n;B)} \,.
	\end{alignat}
\end{subequations}
hold (cf.\ \cite[Theorem\ 1 and 3]{H91}).

Further, we construct a special interpolant in time. 
For a given function $v\in C(\overline{I};B)$, we define, for each $n=1, \dots, N$, a local interpolant
$R_n^1 v\in \P_1({\overline I_n};B)$ by
\begin{equation}
	\label{Eq:IntOpR}
   R_n^1 v(t_{n-1}) \coloneqq  I_{n,2}^\mathrm{GL}\, v(t_{n-1})  = v(t_{n-1}) 
	\quad \text{and} \quad 
	\dt R_n^1 v(\overline t_n) \coloneqq  \dt I_{n,2}^\mathrm{GL} \, v(\overline t_n)
\end{equation}
and a global interpolant $R_\tau^1 v\in Y^1_{\tau}(B)$ by
the initial value $R_\tau^1 v(0)\coloneqq v(0)$ and
\begin{equation}
	\label{Eq:IntOpR2}
   R_\tau^1 v(t) \coloneqq  R_n^1 v(t)\,, \quad t\in I_n\,,
  \quad n=1,\dots,N\,.
\end{equation}

In the notation, we skip to mention the Banach space $B$ in that $R_\tau^1$ is applied. The meaning is always clear from the context. In this work, $R_\tau^1$ is employed in Hilbert spaces $H\subset L^2(\Omega)$. In such spaces, the global operator $R_\tau^1$ and its image have the following properties.
\begin{lem}
	\label{Lem:AppPropR}
	Let $H\subset L^2(\Omega)$ be a Hilbert space with inner product denoted by $\langle \cdot, \cdot \rangle_H$. 
      For $v\in C(\overline{I};H)$, the function $R_\tau^1 v \in C(\overline I;H)$ is continuous in time on $\overline I$ with $R_\tau^1 v(t_n) = v(t_n)$ for all $n=0,\ldots, N$.
\end{lem}

\begin{mproof}
To prove the assertion, we firstly show that $R_n^1 v(t_n) =
I_{n,2}^\mathrm{GL}\, v(t_n) = v(t_n)$ for all $n=1,\ldots, N$. For an
arbitrary time-independent test function $w\in H\subset L^2(\Omega)$, we
obtain by the fundamental theorem of calculus that 
\begin{equation*}
\begin{aligned}
\langle  R_n^1 v(t_n) , w  \rangle_H 	& = \langle  R_n^1 v(t_{n-1}) , w
\rangle_H + \int_{I_n} \langle \partial_t R_n^1 v, w \rangle_H   \ud t \\[1ex]
& = \langle I_{n,2}^\mathrm{GL} \, v (t_{n-1}) , w  \rangle_H + Q_n^{\text{G}} (\langle \partial_t R_n^1  v , w \rangle_H) \\[1ex]
& = \langle I_{n,2}^\mathrm{GL} \, v (t_{n-1}) , w  \rangle_H + Q_n^{\text{G}} (\langle \partial_t   I_{n,2}^\mathrm{GL} \, v, w \rangle_H) \,. 
\end{aligned}	
\end{equation*}
Since $\langle \partial_t   I_{n,2}^\mathrm{GL} \, v, w \rangle  \in
\mathbb P_1(I_n;\R)$ and the Gauss
quadrature formula is exact for all polynomials in $\mathbb
P_1(I_n;\R)$, it follows 
\begin{equation*}
 \langle  R_n^1 v(t_n) , w  \rangle = \langle I_{n,2}^\mathrm{GL} \, v (t_{n-1}) , w  \rangle_H +  \int_{I_n} \langle \partial_t   I_{n,2}^\mathrm{GL} \, v, w \rangle_H \ud t  = \langle I_{n,2}^\mathrm{GL} \, v (t_{n}) , w  \rangle_H
\qquad \forall w \in H\,,
\end{equation*}
proving that $R_n^1 v(t_n) = I_{n,2}^\mathrm{GL} \, v (t_{n}) = v(t_n)$ for all $n=1,\ldots, N$. 

Now, we show the continuity of $R_\tau^1  \,v$ on $[0,T]$ . The function $R_\tau^1  \,v$ is continuous at $0$, since  $(R_\tau^1 \, v)^+(0) =  I_{1,2}^\mathrm{GL}\, v(0)  = v(0) = R_\tau^1 \, v(0)$. From the definition  \eqref{Eq:IntOpR2} along with \eqref{Eq:IntOpR} we get, for $n>1$, that $(R_\tau^1 \,v)^+(t_{n-1}) = R_n^1 \,v(t_{n-1}) = v(t_{n-1})=R_{n-1}^1 \,v(t_{n-1}) =R_\tau^1 \,v(t_{n-1})$, which proves the continuity of $R_\tau^1\, v$ at the point $t_{n-1}$ for all $n=2, \dots, N$ and, thus, the assertion of the lemma. 
\end{mproof}

We note that Lemma \ref{Lem:AppPropR} implies that the operators
$R^1_{\tau}$ and $I_1^{\mathrm{GL}}$ coincide.

\begin{lem}
	\label{Lem:AppPropR2}
	For all $n=1,\ldots, N$ and all $v\in C^{2}(\overline{I}_n;H)$, there holds
	\begin{equation}
		\label{Eq:AppPropR}
		\|v - R_n^1 v \|_{C(\overline{I}_n;H)} \lesssim 
		\tau_n^{2}\|v\|_{C^{2}(\overline{I}_n;H)}\,.
	\end{equation}
\end{lem}

\begin{mproof}
The result follows from \cite[Lemma\ 4.4]{ES16} or, since $R^1_{\tau}$ and $I_1^{\mathrm{GL}}$ coincide, from \eqref{Eq:ApLag_0}.
\end{mproof}

For the error estimation we also need a special projection operator mapping
from $\vec V^{\operatorname{div}}\times Q$ to $\vec
V_h^{\operatorname{div}} \times Q_h$, the so-called Stokes projection (cf.~\cite[Eq.~(4.54)]{J16}), that is defined by means of an auxiliary problem of Stokes type.  

\begin{defi}[Stokes projection]
\label{Def:StokesProj}
For $(\vec u, p)\in \vec V^{\operatorname{div}} \times Q$, the Stokes
projection $R_h(\vec u,p) = \big(R_h^{\vec u}(\vec u,p), R_h^{p}(\vec
u,p)\big)$ is defined by
\begin{equation}
R_h^{\vec u}(\vec  u,p)\coloneqq   \vec{\widehat u}_h
\quad \text{and} \quad
R_h^{p}(\vec u,p)\coloneqq  \widehat p_h\,,
\end{equation}
where the pair $(\vec{\widehat u}_h,\widehat p_h)\in \vec V_h \times Q_h$
satisfies
\begin{subequations}
\label{Eq:StokesProj1}
\begin{alignat}{2}
\label{Eq:StokesProj_1}
\langle A_h \vec{\widehat u}_h + B_h' \widehat p_h ,\vec v_h \rangle
& = \langle A_h \vec u + B_h' p  ,\vec v_h \rangle
& \qquad \forall \vec v_h & \in \vec V_h\,,  \\[1ex]
\label{Eq:StokesProj_2}
\langle B_h \vec{\widehat u}_h, q_h \rangle & = 0
& \qquad \forall q_h & \in Q_h\,.
\end{alignat}
\end{subequations}
\end{defi}
Based on~\eqref{Eq:StokesProj_2}, the velocity part $R_h^{\vec u}(\vec u,p)$ of the Stokes projection $R_h(\vec u,p)$ belongs to the space ${\vec V}_h^{\operatorname{div}}$
of discretely divergence-free functions.


\begin{rem} 
\label{Rem:Com_R_R}	
For arbitrary open intervals $J\subset [0,T]$ and functions
  $(\vec v, q) \in C^1(J; {\vec V^{\operatorname{div}} } \times Q)$, the Stokes projection of
Def.~\ref{Def:StokesProj} yields the commutation property
\begin{equation}
\label{Eq:comm_R_dt}
\partial_t R_h (\vec v, q)  = 
R_h(\partial_t \vec v,\partial_t q)
\end{equation}
in the interval $J$. The commutation of projection in space and
differentiation in time can be proved easily by taking the limit of the
difference quotient for the time derivative and exploiting the linearity of
the Stokes projection $R_h$. The commutation property can be extended in
the same way to higher order temporal derivatives.
\end{rem}


\begin{lem}[Error estimates for the Stokes projection]
\label{Lem:StokesApprox}
Let the stationary Stokes problem be $H^2$-regular in the sense of (A1) in
\cite{HR82}. Then, the Stokes projection
$R_h(\vec u,p) = (\vec{\widehat u}_h,\widehat p_h)\in
\vec V_h^{\operatorname{div}} \times Q_h$ of Def.~\ref{Def:StokesProj}
provides for $\vec u\in \vec V^{\operatorname{div}}\cap
\vec H^{r+1}(\Omega )$ and $p\in Q\cap H^r(\Omega)$ the estimates
\begin{subequations}
\label{Eq:StokesApprox0}
\begin{align}
\label{Eq:StokesApprox1}
\|\vec u - \vec{\widehat u}_h \| + h \|\nabla (\vec u - \vec{\widehat u}_h ) \|  & \lesssim 
h^{r+1}\,,\\[1ex]
\label{Eq:StokesApprox2}
\| p - \widehat p_h \| & \lesssim h ^r\,.
\end{align}
\end{subequations}
\end{lem}

\begin{mproof}
The estimation of $\|\nabla (\vec u - \vec{\widehat u}_h ) \|$ and $\| p -
\widehat p_h \|$ directly follows from \cite[Lemma~4.43]{J16} along with the
Assumption~\ref{Ass:FES} about the approximation properties of $\vec V_h$
and $Q_h$. The estimation $	\|\vec u - \vec{\widehat u}_h \| $ is then
obtained by an duality argument (cf.~\cite[Lemma~4.2 and 4.7]{HR82}).
\end{mproof}

For the rest of the paper, we assume the $H^2$-regularity of
the stationary Stokes problem.

\section{Approximation scheme and first error estimates}
\label{Sec:AppS_ErrGP}


In Subsection~\ref{Subsec:DS}, we present the fully discrete approximation
scheme for the Stokes system \eqref{Eq:SE} that uses
linear in time approximations for both velocity and pressure.
The approximation
of the velocity field will be globally continuous in time on $\overline I$,
whereas the discrete pressure trajectory may have jumps at the
endpoints of the time subintervals.
Concerning the
velocity, the proposed approach resembles the continuous
Galerkin--Petrov method that is studied for parabolic problems, e.g., in
\cite{AM89,S10}. 
We will show that our approximation scheme on each time interval
determines uniquely the discrete velocity, while the discrete pressure trajectory is only partially determined.
The remaining degrees of freedom for the discrete pressure will be
determined either by means of a post-processing, which is addressed in 
Section~\ref{Sec:PPCol}, or by interpolation, which is studied in
Section~\ref{Sec:PPInt}.
In Subsection~\ref{Subsec:ErrPGN},  we prove optimal order error estimates
for the discrete velocity in the endpoints of the time intervals 
and 
for the discrete pressure in the Gauss quadrature nodes related to the midpoint
rule.

\subsection{A linear in time  Galerkin method with finite elements in space}
\label{Subsec:DS}

For the discretization of the Stokes system \eqref{Eq:SE}, 
we solve in a time marching process, for increasing $n=1,\ldots,N$,
the sequence of the following local problems 
for the linear in time approximations 
$(\vec u_{\tau,h}(t),p_{\tau,h}(t)) \in \vec V_h \times Q_h$
of the velocity and pressure trajectories $(\vec u(t), p(t))$
defined on the subintervals $I_n$.
For each interval $I_n =(t_{n-1},t_n]$, we assume
that the trajectory $\uth(t)$ has been already computed by the time
marching process for all $t\in [0,t_{n-1}]$, starting with an
approximation $\uth(0)=\vec u_{0,h}\in \vec V_h^{\operatorname{div}}$ of
the initial value $\vec u_0$.

\begin{prob}[Space-time problem related to {$I_n =(t_{n-1},t_n]$}]
\label{Prob:FdS}
\mbox{}\\
For given $\vec u_{h}^{n-1}\coloneqq \uth(t_{n-1})\in \vec
  V_h^{\operatorname{div}}$ with  $\uth(t_0) \coloneqq \vec u_{0,h}$, find
  $(\vec u_{\tau,h},p_{\tau,h})\in \mathbb P_1(I_n; \vec V_h) \times
  \mathbb P_1(I_n;Q_h)$, 
  such that  $\vec u_{\tau,h}^+(t_{n-1})= \vec
  u_{h}^{n-1}$ and 
\begin{subequations}
\label{Eq:FdS1}
\begin{alignat}{2}
\label{Eq:FdS2}
\int_{t_{n-1}}^{t_n}  \Big(\langle \partial_t \vec u_{\tau,h} , \vec v_{\tau,h} \rangle 
  + { \langle A_h \vec u_{\tau,h} , \vec v_{\tau,h}  \rangle }
  + { \langle B_h'  p_{\tau,h} , \vec v_{\tau,h}  \rangle }
  \Big) \ud t  & 
    = Q^{\operatorname{GL}}_n (\langle \vec f, \vec v_{\tau,h} \rangle)\,, \\[1ex]
\label{Eq:FdS3}
\int_{t_{n-1}}^{t_n} 
  { \langle B_h \vec u_{\tau,h}, q_{\tau,h} \rangle }
  \ud t  & = 0
\end{alignat}
\end{subequations}
for all $(\vec v_{\tau,h},q_{\tau,h}) \in \mathbb P_0(I_n;\vec V_h)\times \mathbb P_0(I_n;Q_h)$.
\end{prob}

By expanding  the discrete solution $(\vec u_{\tau,h}(t),p_{\tau,h}(t))$,
$t\in I_n$, in terms of temporal basis functions, we can
recast~\eqref{Eq:FdS1} as a space problem related to $I_n$ for the
corresponding $\vec V_h$-valued and $Q_h$-valued degrees of freedom.
To this end, let the Lagrange basis polynomials $\varphi_{n,k} \in \mathbb
P_1(\overline I_n;\R)$, $k=0,1$,  on
$\overline I_n$ be defined by the conditions 
\begin{subequations}
\label{Eq:LagBas}
\begin{alignat}{4}
  \varphi_{n,0}(t_{n-1}) & = 1\,, && \qquad \varphi_{n,0}(\overline t_n) & = 0\,, \\[1ex]
  \varphi_{n,1}(t_{n-1}) & = 0\,, && \qquad \varphi_{n,1}(\overline t_n) & = 1\,, 	
\end{alignat}
\end{subequations} 
where $\tnb$ denotes the midpoint of $I_n$ (see 
\eqref{Eq:Midpt}).
In terms of these basis functions we can uniquely represent
$\uth\in \mathbb{P}_1(I_n; \vec V_h)$ and $\pth\in \mathbb{P}_1(I_n;
Q_h)$
for all $t\in I_n$ by 
\begin{equation}
\label{Eq:Rep_u_p}
\uth(t) = \vec u_h^{n-1}\varphi_{n,0}(t) +  \buhn\varphi_{n,1}(t)
\qquad\text{and}\qquad 
\pth(t) = p_h^{n-1}\varphi_{n,0}(t) +  \bphn\varphi_{n,1}(t)\,,
\end{equation}
where the coefficients $\vec u_h^{n-1},\buhn\in \vec V_h$ and $p_h^{n-1},\bphn\in Q_h$ 
have the following interpretation as degrees of freedom of the
discrete solution $(\vec u_{\tau,h},p_{\tau,h})$ on $I_n$,
\begin{equation}
\label{Eq:coeff_u_p}
 \vec u^{n-1}_h \coloneqq \vec u^+_{\tau,h}(t_{n-1})\,, \quad
\unbh \coloneqq \uth(\tnb) \,, \qquad\quad
 p^{n-1}_h \coloneqq p^+_{\tau,h}(t_{n-1})\,, \quad
\bphn \coloneqq \pth(\tnb) \,.
\end{equation}

We can now substitute the expansions \eqref{Eq:Rep_u_p} into the
variational equations \eqref{Eq:FdS1}.  Then, we obtain the following
space problem for the degrees of freedom $\unbh$  and $\bphn$
of the discrete solution $(\vec u_{\tau,h},p_{\tau,h})$ on $I_n$. 
For the derivation of Problem~\ref{Prob:FdSAF}, we refer to the proof of
Lemma~\ref{Lem:space_time_probl}. 


\begin{prob}[{Space problem related to $I_n$}]
\label{Prob:FdSAF}
\mbox{}\\
For given $\vec u_{h}^{n-1}\coloneqq \uth(t_{n-1})\in \vec V_h^{\operatorname{div}}$ with  $\uth(t_0) \coloneqq \vec u_{0,h}$, 
find $( \unbh ,\bphn) \in \vec V_h \times Q_h$ 
such that 
\begin{subequations}
\label{Eq:FdSAF_0}
\begin{align}
\label{Eq:FdSAF_1}
\langle \buhn -\vec u_h^{n-1}, \vec v_h \rangle + \frac{\tau_n}{2} \langle A_h \buhn + B_h' 
  \bphn,\vec v_h  \rangle & = \frac{\tau_n}{4}\langle \vec f(t_{n-1}) +  \vec f(t_n) , \vec v_h 
\rangle\,, \\[1ex]
\label{Eq:FdSAF_2}
\langle B_h \buhn,q_h\rangle & = 0
\end{align}
\end{subequations}
for all $(\vec v_h,q_h)\in \vec V_h\times Q_h$.
\end{prob}

Regarding  existence and uniqueness of solutions to
Problem~\ref{Prob:FdSAF} we have the following result. 

\begin{lem}[Well-posedness of Problem~\ref{Prob:FdSAF}]
\label{Lem:space_probl}
Problem \ref{Prob:FdSAF} admits a unique solution  $( \unbh ,\bphn) \in \vec V_h \times Q_h$. Moreover, we have that $\unbh\in\vec V_h^{\operatorname{div}}$.
\end{lem}

\begin{mproof}
Since Problem~\ref{Prob:FdSAF} is equivalent to a finite dimensional
linear system of equations for the coefficients of $\unbh$ and $\pnbh$
in their basis representation, 
the existence of a solution is a direct consequence of its uniqueness.
To prove uniqueness, we suppose that $( \overline{\vec u}^{n,k}_h,\overline{p}^{n,k}_h) \in \vec V_h \times Q_h$,
  $k=1,2$, are two solutions satisfying \eqref{Eq:FdSAF_0} with
$\vec u^{n-1,1}_h = \vec u^{n-1,2}_h = \vec u_{\tau,h}(t_{n-1})$.
  Then, for their difference it holds 
\begin{subequations}
	\label{Eq:FdSAF_3}
	\begin{align}
		\label{Eq:FdSAF_4}
		\langle  \overline{\vec u}^{n,1}_h -  \overline{\vec u}^{n,2}_h , \vec v_h \rangle + \frac{\tau_n}{2} \langle A_h ( \overline{\vec u}^{1,2}_h- \overline{\vec u}^{n,2}_h) + B_h'  (\overline{p}^{n,1}_h-\overline{p}^{n,2}_h) ,\vec v_h  \rangle & = 0\,, \\[1ex]
		\label{Eq:FdSAF_5}
		\langle B_h ( \overline{\vec u}^{n,1}_h -  \overline{\vec u}^{n,2}_h ),q_h\rangle & = 0
	\end{align}
\end{subequations}
for all $(\vec v_h,q_h)\in \vec V_h\times Q_h$ since the difference
 $\vec u^{n-1,1}_h - \vec u^{n-1,2}_h$ vanishes.
  We choose the test function 
  $\vec v_h=\overline{\vec u}^{n,1}_h -  \overline{\vec u}^{n,2}_h$ 
  in \eqref{Eq:FdSAF_4} and use that 
  $\vec v_h \in \vec V_h^{\operatorname{div}}$ due to \eqref{Eq:FdSAF_5}.
Then, the pressure term in  \eqref{Eq:FdSAF_4} vanishes and we obtain 
  $\| \vec v_h \|^2 + (\tau_n/2) \| \nabla\vec v_h \|^2 = 0 $ 
which implies $\vec v_h = \vec 0$ 
and, thus, $\overline{\vec u}^{n,1}_h = \overline{\vec u}^{n,2}_h$.  
Inserting this again into  \eqref{Eq:FdSAF_4} we get 
$
\langle B_h'(\overline{p}^{n,1}_h-\overline{p}^{n,2}_h), \vec v_h \rangle = 0
$
for all $\vec v_h \in \vec V_h $.
  The inf-sup stability condition
  \eqref{Eq:InfSupCod}  then yields the existence of some 
  $ \vec v_h\in \vec V_h \setminus\{\vec 0\}$ such that 
\begin{equation*}
  \beta \, \| \overline{p}^{n,1}_h-\overline{p}^{n,2}_h \| \,\le\,
\dfrac{b(\vec v_h, \overline{p}^{n,1}_h-\overline{p}^{n,2}_h)}%
{\| \vec v_h \|_1 } =
\dfrac{\langle B_h'(\overline{p}^{n,1}_h-\overline{p}^{n,2}_h), \vec v_h
\rangle}{\| \vec v_h \|_1 }
= 0  \,.
\end{equation*}
Hence, we get that $\overline{p}^{n,1}_h=\overline{p}^{n,2}_h$ and the unique solvability follows.
Furthermore, the property $\unbh\in\vec V_h^{\operatorname{div}}$ is a direct consequence of~\eqref{Eq:FdSAF_2}.
Thus, Lemma~\ref{Lem:space_probl} is proved.
\end{mproof}

\begin{lem}[Set of solutions of Problem~\ref{Prob:FdS}]
\label{Lem:space_time_probl}
The set of all solutions 
$ (\vec u_{\tau,h},p_{\tau,h})\in 
          \mathbb P_1(I_n; \vec V_h) \times \mathbb P_1(I_n;Q_h) 
$
of the space-time Problem~\ref{Prob:FdS} on interval $I_n =(t_{n-1},t_n]$
is characterized by the representation~\eqref{Eq:Rep_u_p}, where
the coefficient 
    $\vec u^{n-1}_h = \vec u^+_{\tau,h}(t_{n-1})$ satisfies
$\vec u_{h}^{n-1}=\uth(t_{n-1})$ and the pair of coefficients
$
  (\unbh ,\bphn) =  ( \uth(\tnb) , \pth(\tnb) ) \in \vec V_h \times Q_h
$
is the unique solution of Problem~\ref{Prob:FdSAF}, whereas the coefficient
$p^{n-1}_h =  p^+_{\tau,h}(t_{n-1})\in Q_h$
can be chosen arbitrarily.
\end{lem}

\begin{mproof}
Let 
$ (\vec u_{\tau,h},p_{\tau,h})\in 
          \mathbb P_1(I_n; \vec V_h) \times \mathbb P_1(I_n;Q_h) 
$
be given.
Then, this pair has the unique representation~\eqref{Eq:Rep_u_p}, where the
coefficients $\vec u_h^{n-1},\buhn\in \vec V_h$ and $p_h^{n-1},\bphn\in Q_h$ 
are defined by~\eqref{Eq:coeff_u_p}.
Let further $\vec v_{\tau,h} \in \mathbb P_0(I_n;\vec V_h)$
and $q_{\tau,h} \in \mathbb P_0(I_n;Q_h)$ be given with the degrees of
freedom  $\vec v_h \coloneqq  \vec v_{\tau,h}(\tnb) \in \vec V_h$
and  $q_h \coloneqq  q_{\tau,h}(\tnb) \in Q_h$. 

  We will now show that the system \eqref{Eq:FdSAF_0} of
  Problem~\ref{Prob:FdSAF} arises from the system \eqref{Eq:FdS1} of
  Problem~\ref{Prob:FdS}.
Under the above assumptions, the integrands of the time integrals on the left-hand side of the
variational equations of the space-time Problem~\ref{Prob:FdS} 
are first order polynomials in time. 
Since the one-point Gauss quadrature
formula \eqref{Eq:GF} is exact for these integrals, we obtain
that the variational equations~\eqref{Eq:FdS1}
are equivalent to
\begin{subequations}
\begin{align*}
  \tau_n \langle \partial_t \vec u_{\tau,h}(\tnb), \vec v_h \rangle 
  + \tau_n \langle A_h \buhn + B_h' \bphn,\vec v_h  \rangle 
  & = 
  \frac{\tau_n}{2}\langle \vec f(t_{n-1}) + \vec f(t_n), \vec v_h \rangle\,, \\[1ex]
 \tau_n \langle B_h \buhn,q_h\rangle & = 0
\end{align*}
\end{subequations}
for all $(\vec v_h,q_h)\in \vec V_h\times Q_h$.
From \eqref{Eq:Rep_u_p} we get for the time derivative of $\vec u_{\tau,h}$
that
$
  \partial_t \vec u_{\tau,h}(\tnb) = (\unbh - \vec u^{n-1}_h)/(\tau_n/2)
$.
Therefore,  the variational equations~\eqref{Eq:FdS1}
of the space-time formulation in Problem~\ref{Prob:FdS} are equivalent to the
variational equations~\eqref{Eq:FdSAF_0} of the 
space problem for the coefficients $(\unbh, \bphn)$
where 
    $\vec u^{n-1}_h \coloneqq  \vec u^+_{\tau,h}(t_{n-1})$.

Let now $(\vec u_{\tau,h},p_{\tau,h})$  be a solution of the space-time
Problem~\ref{Prob:FdS} on interval $I_n$. Then, it must satisfy the
condition 
$\vec u^+_{\tau,h}(t_{n-1}) = \vec u_{\tau,h}(t_{n-1})$, i.e., 
the coefficient $\vec u^{n-1}_h$ in the representation~\eqref{Eq:Rep_u_p}  
of $\vec u_{\tau,h}$  must satisfy
$\vec u^{n-1}_h = \vec u_{\tau,h}(t_{n-1})$.
Furthermore, the variational equations~\eqref{Eq:FdS1},
which are equivalent to the variational equations~\eqref{Eq:FdSAF_0}
in space, imply that
the pair of the coefficients $(\unbh ,\bphn) $ in~\eqref{Eq:Rep_u_p}
is a solution of the space problem~\eqref{Eq:FdSAF_0} on $I_n$. 
Since the coefficient $p^{n-1}_h$ of $p_{\tau,h}$ does not occur in
Problem~\ref{Prob:FdSAF}, one can choose $p^{n-1}_h\in Q_h$
arbitrarily without violating the variational equations~\eqref{Eq:FdSAF_0} of the 
space problem which are equivalent to the variational equations~\eqref{Eq:FdS1}
of the space-time formulation in  Problem~\ref{Prob:FdS}.

Therefore, if the coefficients 
$\vec u_h^{n-1},\buhn\in \vec V_h$ and $p_h^{n-1},\bphn\in Q_h$ 
in the representation~\eqref{Eq:Rep_u_p}  
of $\vec u_{\tau,h}$ and  $p_{\tau,h}$ are chosen such that 
$\vec u^{n-1}_h \coloneqq  \vec u_{\tau,h}(t_{n-1})$, the pair $(\unbh ,\bphn)$ is
the unique solution of Problem~\ref{Prob:FdSAF} on $I_n$ (due to
Lemma~\ref{Lem:space_probl}) and $p^{n-1}_h$  is an arbitrary element of
$Q_h$, then the first order polynomials in time  
$(\vec u_{\tau,h}, p_{\tau,h})$, which are built with these coefficients,
solve Problem~\ref{Prob:FdS}.
\end{mproof}

\begin{rem}
\label{Rem:P}
The statement of Lemma~\ref{Lem:space_time_probl} provides the following
  conclusions.
\begin{itemize}
\item[(i)]
  Lemma~\ref{Lem:space_time_probl} shows that, under all solutions of the 
  space-time Problem~\ref{Prob:FdS}, the velocity
part $\uth$ on $I_n$ is uniquely determined since the two coefficients 
$\vec u_{h}^{n-1}\coloneqq \uth(t_{n-1})$ and $\unbh$ in the 
representation~\eqref{Eq:Rep_u_p} are uniquely determined.
Due to \eqref{Eq:coeff_u_p}
we have $\vec u^{n-1}_h \coloneqq \uth(t_{n-1}) = \vec u^+_{\tau,h}(t_{n-1})$ 
which implies the
continuity of the global velocity trajectory $\uth(t)$ at $t=t_{n-1}$.
Thus, the time marching process solving the local space-time 
problems~\ref{Prob:FdS} generates a unique globally continuous velocity
trajectory  $\vec u_{\tau,h}\in X_\tau^1(\vec V_h)$. 
\item[(ii)]
Contrary to the velocity part of the solutions of Problem~\ref{Prob:FdS}
on $I_n$, the pressure part $\pth$ is not uniquely 
determined since its coefficient 
$p^{n-1}_h \coloneqq  p^+_{\tau,h}(t_{n-1})$ 
in the representation~\eqref{Eq:Rep_u_p} can be chosen arbitrarily.
This also implies that there is in general no continuity property of the 
global pressure trajectory $p_{\tau,h}(t)$ at $t=t_{n-1}$ which leads to
globally discontinuous pressure solutions $p_{\tau,h}\in Y_\tau^1(Q_h)$. 
However, Lemma~\ref{Lem:space_time_probl} also shows that the
coefficient $ \bphn \coloneqq   \pth(\tnb) $ is uniquely determined as the
pressure part of the unique solution of the space 
problem~\ref{Prob:FdSAF} on $I_n$.
Thus, all the different pressure solutions $\pth$ of Problem~\ref{Prob:FdS} 
coincide in the value  $\pth(\tnb)$.
This should be taken into account in Section~\ref{Subsec:ErrPGN},
where error estimates will be proved for the 
pressure point values $\pth(\tnb)$.
\end{itemize}
\end{rem}

\subsection{Estimation of the error in the Gauss quadrature nodes}
\label{Subsec:ErrPGN}

In the following, let $\vec u_{\tau,h}\in X_\tau^1(\vec V_h)$ denote
the unique  global velocity part and $p_{\tau,h}\in Y_\tau^1(Q_h)$
any of the different global pressure solutions created by solving 
sequentially the local systems of Problem~\ref{Prob:FdS}  within the time marching process.

Now we shall derive estimates for the velocity and pressure errors
\begin{equation}
\label{Def:Err_0}
  \vec e_{\vec u} \coloneqq  \vec u - \vec u_{\tau,h}
  \qquad\text{and}\qquad
  e_p \coloneqq   p -  p_{\tau,h}
\end{equation}
with the solution $(\vec u,p)$ of \eqref{Eq:SE} and its discrete approximation $(\uth, \pth)\in X_\tau^1(\vec V_h)\times 
Y_\tau^1(Q_h)$ being defined locally on $I_n$ by Problem~\ref{Prob:FdS}
or \ref{Prob:FdSAF}, respectively (cf.\ Remark~\ref{Rem:P}).
We split the errors \eqref{Def:Err_0} into  
\begin{subequations}
\label{Eq:SpErr0}
\begin{alignat}{3}
\label{Eq:SpErr1}
\eu & = \etau + \euth\,, & \qquad \etau & \coloneqq  \vec u - \udth, &
\quad  \euth & \coloneqq  \udth - \uth \\[1ex]
\label{Eq:SpErr2}
e_p & = \eta^p + \epth\,, & \qquad \eta^p & \coloneqq  p - \pdth\, &
\quad \epth & \coloneqq  \pdth - \pth \,,
\end{alignat}
\end{subequations}
with the space-time interpolants $(\udth,\pdth)\in X_\tau^1(\vec
V_h^{{\operatorname{div}}})\times X_\tau^1(Q_h)$ defined by 
\begin{equation}
\label{Eq:DefSP}
  \udth \coloneqq    R_h^{\vec u}( R_\tau^1 \vec u, R_\tau^1 p )
  \qquad\text{and}\qquad
  \pdth \coloneqq    R_h^{p}( R_\tau^1 \vec u, R_\tau^1 p ) \,,
\end{equation}
where $R_h^{\vec u}$ and $R_h^{p}$ are the velocity and pressure
components of the Stokes projection $R_h$ given in
Def.~\ref{Def:StokesProj},
while $R_\tau^1$ is the interpolant in time defined in~\eqref{Eq:IntOpR}.
We refer 
to $(\vec \eta^{\vec u},\eta^p)$ as the interpolation error and to $ (\vec 
e_{\tau,h}^{\vec u},e_{\tau,h}^p)$ as the discrete error. 

%
\begin{rem}
\label{Rem:DiscDivFree}
The discrete velocity error $\euth$ is discretely divergence-free pointwise
in time. This follows from~\eqref{Eq:SpErr1}, the properties of the Stokes
projection, the representation~\eqref{Eq:Rep_u_p}, and $\vec u_h^{n-1},
  \buhn\in \vec V_h^{\operatorname{div}}$  (cf.\ Lemma~\ref{Lem:space_probl} for
  the last condition). Since the discrete velocity error $\euth$ is piecewise
linear and discretely divergence-free pointwise in time, 
  {we obtain that the time derivative $\partial_t \euth$ is also discretely 
  divergence-free in the interior of each time interval.}
\end{rem}

\begin{lem}[Estimate of the error $(\vec \eta^{\vec u}, \eta^p)$]
\label{Lem:ApproxErr}
For the error $(\vec \eta^{\vec u},\eta^p)=(\vec u - \vec{\widehat u}_{\tau,h},p- 
\widehat p_{\tau,h})$,
there holds 
\begin{equation}
\label{Eq:ApproxErr00} 
|Q_n^{\operatorname G}(\langle \partial_t \vec \eta^{\vec u} + 
A_h \vec \eta^{\vec u} + B_h' \eta^p , \vec v_{\tau,h}\rangle)| \lesssim \tau_n 
(\tau_n^2 + h^{r+1}) \| \vec 
v_{\tau,h}(\overline t_n)\|
\end{equation}
for all $\vec v_{\tau,h} \in \mathbb P_0(I_n;\vec V_h)$.
\end{lem}

\begin{mproof}
Let $T_n^{\vec \eta}
\coloneqq  Q_n^{\operatorname G}(\langle \partial_t \vec \eta^{\vec u} +
A_h \vec \eta^{\vec u} + B_h' \eta^p , \vec v_{\tau,h}\rangle)$. It follows that 
$T_n^{\vec \eta} = T_n^{\vec \eta,1} + T_n^{\vec \eta,2}$ with
\begin{align*}
T_n^{\vec \eta,1} \coloneqq  Q_n^{\operatorname G}(\langle \partial_t \vec u -
\partial_t  \vec{\widehat  u}_{\tau,h}, \vec v_{\tau,h}\rangle),\qquad
T_n^{\vec \eta,2} \coloneqq  Q_n^{\operatorname G}(\langle A_h (\vec u -
\vec{\widehat u}_{\tau,h}) + B_h' (p- \widehat p_{\tau,h}), \vec
v_{\tau,h}\rangle )
\,. 
\end{align*}
For $T_n^{\vec \eta,1}$, we conclude by \eqref{Eq:IntOpR},
\eqref{Eq:comm_R_dt},  \eqref{Eq:ApLag_2}, and \eqref{Eq:StokesApprox0}, applied to $(\partial_t R_\tau^1 \vec u,\partial_t R_\tau^1 p)$, that 
\begin{align*}
 T_n^{\vec \eta,1} & = 
 Q_n^{\operatorname G}(\langle \partial_t \vec u - \partial_t R_ 
 \tau^1 \vec u, \vec v_{\tau,h}\rangle) + Q_n^{\operatorname G}(\langle \partial_t 
  R_\tau^1 \vec u - \partial_t R_h^{\vec u}(R_\tau^1 \vec u,R_\tau^1 p),
  \vec v_{\tau,h}\rangle)\\[1ex]
& = Q_n^{\operatorname G}(\langle \partial_t \vec u - \partial_t
I_{n,2}^\mathrm{GL}\, \vec u, \vec v_{\tau,h}\rangle) + Q_n^{\operatorname G}(\langle \partial_t R_\tau^1 \vec u - 
R_h^{\vec u}(\partial_t R_\tau^1 \vec u,\partial_t R_\tau^1 p), \vec
v_{\tau,h}\rangle)\\[1ex]
& \lesssim \tau_n  (\| \partial_t \vec u - \partial_t  I_{n,2}^\mathrm{GL}
\vec u \|_{C(\overline I_n;\vec L^2)} + \| \partial_t R_\tau^1 \vec u - 
R_h^{\vec u}(\partial_t R_\tau^1 \vec u,\partial_t R_\tau^1 p)\| _{C(\overline I_n;\vec L^2)}) \| \vec 
v_{\tau,h}(\overline t_n) \|\\[1ex]
& \lesssim \tau_n (\tau_n^2 + h^{r+1}) \| \vec 
v_{\tau,h}(\overline t_n) \|\,.
\end{align*}
For $T_n^{\vec \eta,2}$, we find by \eqref{Eq:StokesProj_1}, applied to
$(R_\tau^1 \vec u,R_\tau^1 p)$, \eqref{Def:A_h}, \eqref{Def:B_h^p}, and property \eqref{Eq:AppPropR} of $R_\tau^1 $ that 
\begin{align*}
 T_n^{\vec \eta,2} & = Q_n^{\operatorname G}(\langle A_h (\vec u - R_\tau^1 \vec u) + B_h' 
(p - R_\tau^1 p), \vec v_{\tau,h}\rangle)\\[1ex]
& \quad + Q_n^{\operatorname G}(\langle A_h (
R_\tau^1 \vec u - R_h^{\vec u}(R_\tau^1 \vec u,R_\tau^1 p)) + B_h'
(R_\tau^1 p - R_h^{p}(R_\tau^1 \vec u,R_\tau^1 p)),{\vec v_{\tau,h}}
\rangle) \\[1ex]
& = Q_n^{\operatorname G}(\langle A_h (\vec u - R_\tau^1 \vec u) + B_h' 
(p - R_\tau^1 p), \vec v_{\tau,h}\rangle)\\[1ex]
& = Q_n^{\operatorname G}(\langle - \Delta (\vec u - R_\tau^1 \vec u) + \nabla 
(p - R_\tau^1 p), \vec v_{\tau,h}\rangle)\\[1ex]
& \lesssim \tau_n \tau_n^2  \| \vec v_{\tau,h}(\overline t_n) \|\,.
\end{align*}
Finally, combining the previous estimates yields the assertion \eqref{Eq:ApproxErr00}. 
\end{mproof}

\begin{lem}[Estimate of the discrete error $\vec e_{\tau,h}^{\vec u}$]
\label{Lem:DiscErr}
For the discrete velocity error $\vec e_{\tau,h}^{\vec u} = \udth -\vec 
u_{\tau,h}$, there holds 
\begin{subequations}
\label{Eq:DiscErr000}	
\begin{alignat}{2}
\label{Eq:DiscErr00} 
\sum_{m=1}^n \tau_m \| \partial_t \vec e_{\tau,h}^{\vec u}(\overline t_m ) \|^2 + \| 
\nabla \vec e_{\tau,h}^{\vec u} (t_n) \|^2 & \lesssim \| \nabla \vec
e_{\tau,h}^{\vec u} (0) 
\|^2 + \tau^4 + h^{2(r+1)}\,,\\[1ex]
\label{Eq:DiscErr01}
\| \vec e_{\tau,h}^{\vec u} (t_n) \|^2 + \sum_{m=1}^n \tau_m \| \nabla 
\vec e_{\tau,h}^{\vec u} (\overline t_m) \|^2 & \lesssim  \| \vec
e_{\tau,h}^{\vec u} (0) \|^2 + \tau^4 +  h^{2(r+1)} 
\end{alignat}
\end{subequations}
for all $n=1,\ldots , N$.
\end{lem}

\begin{mproof}
For a sufficiently regular solution $(\vec u,p)$ of \eqref{Eq:SE}, equation
\eqref{Eq:SE_1} is satisfied for all $t\in \overline I$. Together 
with \eqref{Eq:FdS2}, this yields the error equation 
\begin{align*}
 Q_n^{\text G}(\langle \partial_t \vec e_{\vec u} + A_h \vec e_{\vec u} + B_h' e_p , \vec 
v_{\tau,h} \rangle) ={Q_n^{\text G}(\langle \vec f, \vec v_{\tau,h}\rangle)} - 
Q_n^{\text{GL}} (\langle \vec f, \vec v_{\tau,h}\rangle)
\end{align*}
for all $\vec v_{\tau,h}\in \mathbb P_0(I_n;\vec V_h)$. The
splitting~\eqref{Eq:SpErr0} of the errors then implies 
\begin{equation}
\label{Eq:EQDS}
Q_n^{\text G}(\langle \partial_t \vec e_{\tau,h}^{\vec u}
+ A_h \vec e_{\tau,h}^{\vec u} 
+ B_h' e_{\tau,h}^p , \vec v_{\tau,h} \rangle)
= T_n^{\vec f} + T_n^{\vec \eta}
\end{equation}
with
\begin{equation*}
T_n^{\vec f} \coloneqq  \tau_n \langle \vec \delta_{\vec f}^n 
,  \vec 
v_{\tau,h}(\overline t_n) \rangle\,,\qquad
T_n^{\vec \eta} \coloneqq  - Q_n^{\text G}(\langle \partial_t \vec \eta^{\vec u} + 
A_h \vec \eta^{\vec u} + B_h' \eta^p , \vec v_{\tau,h}\rangle)
\end{equation*}%
for all $\vec v_{\tau,h}\in \mathbb P_0(I_n;\vec V_h)$, where 
$\vec \delta_{\vec f}^n\coloneqq  \vec 
f(\overline t_n)-\frac{1}{2}\big(\vec f(t_{n-1})+\vec f(t_n)\big)$.

By the approximation property \eqref{Eq:ApLag_0} of
$I_{n,1}^{\operatorname{GL}}$ we have  
\begin{align}
\nonumber
| T_n^{\vec f} |  & \leq  \tau_n \left\| \vec f(\overline
t_n)-\frac{1}{2}\big(\vec 
f(t_{n+1})+ \vec f(t_n)\big)\right\| \, \|  \vec 
v_{\tau,h}(\overline t_n) \|\\[1ex]
\label{Eq:DiscErr045} 
  & = \tau_n \| \vec f(\overline t_n)- I_{n,1}^{\text{GL}}  \vec f(\overline t_n )\| \, \| 
 \vec 
v_{\tau,h}(\overline t_n) \|  \lesssim \tau_n \tau_n^2\| 
\vec 
v_{\tau,h}(\overline t_n) \|\,.
\end{align}
By Lemma~\ref{Lem:ApproxErr} it holds
\begin{equation*}
 |T_n^{\vec \eta} |  \lesssim \tau_n (\tau_n^2 +h^{r+1})\| \vec 
 v_{\tau,h}(\overline t_n) \|\,.
\end{equation*}
Thus, we have  
\begin{equation}
\label{Eq:DiscErr05} 
Q_n^{\text G}(\langle \partial_t \vec e_{\tau,h}^{\vec u} + A_h \vec e_{\tau,h}^{\vec u} 
+ B_h' e_{\tau,h}^p , \vec v_{\tau,h} \rangle) \lesssim \tau_n (\tau_n^2 +h^{r+1})\| 
\vec v_{\tau,h}(\overline{t}_n)\|
\end{equation}
 for all $\vec v_{\tau,h}\in \mathbb{P}_0(I_n; \vec V_h)$.

Firstly, we choose $\vec v_{\tau,h} = \partial_t \euth$ as test function
in~\eqref{Eq:DiscErr05}. Since $\partial_t \euth$ belongs to
$\mathbb{P}_0(I_n; V_h^{\operatorname{div}})$  (cf.\ Remark~\ref{Rem:DiscDivFree}) we have
\begin{equation}
\label{Eq:B_h^p_eq_0}
\langle B_h' q_h, \partial_t \euth(t) \rangle = 0 
\end{equation}%
for all $q_h\in Q_h$ and $t\in I_n$.  From \eqref{Eq:DiscErr05} and \eqref{Eq:B_h^p_eq_0} we then get 
\begin{equation}
\label{Eq:DiscErr06} 
Q_n^{\text G}(\langle \partial_t \vec e_{\tau,h}^{\vec u}, \partial_t \vec 
e_{\tau,h}^{\vec u}\rangle + \langle A_h \vec e_{\tau,h}^{\vec u}, \partial_t \vec 
e_{\tau,h}^{\vec u}\rangle)  \lesssim \tau_n (\tau_n^2 +h^{r+1})\| \partial_t \vec 
e_{\tau,h}^{\vec u} (\overline t_n) \|\,.
\end{equation}
By the definition of the Gauss quadrature formula along with 
\begin{equation*}
Q_n^{\text G}(\underbrace{\langle A_h \vec e_{\tau,h}^{\vec u}, \partial_t \vec 
e_{\tau,h}^{\vec u}\rangle}_{\in \mathbb P_1(I_n;\R)}) = \int_{I_n} \langle A_h \vec e_{\tau,h}^{\vec u}, \partial_t 
\vec 
e_{\tau,h}^{\vec u}\rangle \ud t = \int_{I_n}\dfrac{1}{2}\dfrac{d}{dt} \langle A_h 
\vec e_{\tau,h}^{\vec u}, \vec e_{\tau,h}^{\vec u}\rangle \ud t \,,
\end{equation*}
it follows from 
\eqref{Eq:DiscErr06} 
that 
\begin{equation*} 
\tau_n \| \partial_t \vec e_{\tau,h}^{\vec u}(\overline t_n) \|^2 + \frac{1}{2}\| \nabla 
\vec 
e_{\tau,h}^{\vec u} (t_n) \|^2 - \frac{1}{2}\| \nabla \vec e_{\tau,h}^{\vec u} (t_{n-1}) 
\|^2 \lesssim \tau_n (\tau_n^2 +h^{r+1})\| \partial_t \vec 
e_{\tau,h}^{\vec u}(\overline t_n) \|\,.
\end{equation*}
 By absorption we get 
\begin{equation}
 \label{Eq:DiscErr07}
 \frac{1}{2}\tau_n \| \partial_t \vec e_{\tau,h}^{\vec u}(\overline t_n) \|^2 + 
\frac{1}{2}\| \nabla 
\vec e_{\tau,h}^{\vec u} (t_n) \|^2 - \frac{1}{2}\| \nabla \vec e_{\tau,h}^{\vec u} 
(t_{n-1}) \|^2 \lesssim \tau_n (\tau_n^4 +h^{2(r+1)}) \,.
\end{equation}
Changing now the index $n$ to $m$ and summing up the resulting inequality from $m=1$ to $n$ show 
\begin{equation}
\label{Eq:DiscErr08}
\sum_{m=1}^n \tau_m \| \partial_t \vec e_{\tau,h}^{\vec u}(\overline t_m) \|^2 + \| 
\nabla  \vec e_{\tau,h}^{\vec u} (t_n) \|^2 \lesssim \| \nabla \vec e_{\tau,h}^{\vec u} (0) 
\|^2 + t_n (\tau^4 +h^{2(r+1)}) \,.
\end{equation}
This proves~\eqref{Eq:DiscErr00}.

Secondly, we choose
the test function $\vec v_{\tau,h}$ in~\eqref{Eq:DiscErr05} as those
function from $\mathbb{P}_0(I_n; \vec V_h)$ that has the constant value
$\vec e_{\tau,h}^{\vec u}(\overline{t}_n)$ in $I_n$.
By the arguments
of~\eqref{Eq:B_h^p_eq_0}
with $\vec e_{\tau,h}^{\vec u}(\overline{t}_n)$ replacing
$\partial_t \vec e_{\tau,h}^{\vec u}$, we then obtain
\begin{equation}
\label{Eq:DiscErr095}
Q_n^{\text G}(\langle \partial_t \vec e_{\tau,h}^{\vec u}, \vec 
e_{\tau,h}^{\vec u}\rangle) + Q_n^{\text G}(\langle A_h \vec 
e_{\tau,h}^{\vec u}, \vec e_{\tau,h}^{\vec u}\rangle)  \lesssim \tau_n 
(\tau_n^2 +h^{r+1})\| \vec e_{\tau,h}^{\vec u} (\overline t_n)\|\,.
\end{equation}
For the first term on the left-hand side of \eqref{Eq:DiscErr095}, we get by
the exactness of the 1-point Gauss quadrature formula for all polynomials
in $\mathbb P_1(I_n;\R)$ that 
\begin{equation}
 \label{Eq:DiscErr10}
 Q_n^{\text G}(\langle \partial_t \vec e_{\tau,h}^{\vec u}, \vec 
e_{\tau,h}^{\vec u}\rangle) = \int_{I_n} \langle \partial_t \vec 
e_{\tau,h}^{\vec u}, \vec e_{\tau,h}^{\vec u}\rangle \ud t = \frac{1}{2}\| \vec 
e_{\tau,h}^{\vec u} (t_{n})\|^2 - \frac{1}{2}\| \vec 
e_{\tau,h}^{\vec u} (t_{n-1})\|^2\,.
\end{equation}
Combining \eqref{Eq:DiscErr095} with \eqref{Eq:DiscErr10}, applying the inequality of Cauchy--Young to the term on the right-hand side of  \eqref{Eq:DiscErr095}, using the inequality of Poincar\'e to absorb $\| \vec e_{\tau,h}^{\vec u} (\overline t_n)\|$ by the left-hand side, changing the index $n$ to $m$ in the resulting inequality and, finally, summing up the thus obtained relation from $m=1$ to $n$ show 
\begin{equation}
 \label{Eq:DiscErr11}
\| \vec e_{\tau,h}^{\vec u} (t_n) \|^2 + \sum_{m=1}^n \tau_m \| \nabla \vec 
e_{\tau,h}^{\vec u} (\overline t_m) \|^2 \lesssim  \| \vec e_{\tau,h}^{\vec u} (0) \|^2 + 
 t_n (\tau^4 + h^{2(r+1)})\,. 
\end{equation}
This proves \eqref{Eq:DiscErr01}.%
\end{mproof}

\begin{lem}[Estimate of the discrete pressure error $e_{\tau,h}^p$]
\label{Lem:EstPth}
For the discrete pressure error $e_{\tau,h}^p =  \pdth  -p_{\tau,h}$, there 
holds 
\begin{equation}
\label{Eq:EstP00}
\sum_{m=1}^n \tau_m \| e_{\tau,h}^p(\overline t_m)\|^2  \lesssim 
\| \nabla \vec e_{\tau,h}^{\vec u}(0) \|^2
+ \tau^4 + h^{2(r+1)}
\end{equation}
for $n=1,\ldots N$. 
\end{lem}

\begin{mproof}
By the choice of an inf-sup stable pair of finite element spaces $\vec V_h \times Q_h$, fulfilling the Assumption~\ref{Ass:FES}, 
for $\overline e_{n,h}^{p}\coloneqq  e_{\tau,h}^p(\overline t_n)\in Q_h$ there exists some 
function $\vec v_h \in \vec V_h\setminus \{\vec 0\}$ such that 
\begin{equation}
\label{Eq:EstP0}
\beta \|\overline e_{n,h}^p\| \le 
\frac{ \langle B_h' \overline e_{n,h}^p, \vec v_h \rangle}%
{\| \vec v_h\|_1}  \,,
\end{equation}
where the constant $\beta >0$ is independent of $h$. Combining
\eqref{Eq:EstP0} with equation \eqref{Eq:EQDS} yields 
\begin{equation}
\label{Eq:EstP05}
\begin{aligned}
\|\overline e_{n,h}^p\|  \lesssim \frac{1}{\tau_n}
\frac{1}{\|\vec v_h\|_1} \Big( & \tau_n \langle 
\vec \delta_{\vec f}^n , \vec v_h\rangle - \langle \partial_t \vec e_{\tau,h}^{\vec u}(\overline 
t_n), \vec v_h \rangle - \langle A_h \vec e_{\tau,h}^{\vec u}(\overline t_n),\vec v_h \rangle \\[1ex]
& - Q_n^{\text G}(\langle \partial_t \vec \eta^{\vec u} + 
A_h \vec \eta^{\vec u} + B_h' \eta^p , \vec v_{\tau,h}\rangle) \Big)\,,
\end{aligned}
\end{equation}
with $\vec \delta_{\vec f}^n$ being defined in the proof of
Lemma~\ref{Lem:DiscErr}, where $\vec v_{\tau,h}\in \mathbb{P}_0(I_n;
\vec V_h)$ is the function taking the constant value $\vec v_h$ in $I_n$.
By 
the inequality of Cauchy--Schwarz along with the inequalities \eqref{Eq:ApproxErr00} and \eqref{Eq:DiscErr045}  we deduce from \eqref{Eq:EstP05} that 
\begin{equation}
\label{Eq:EstP055}
\|\overline e_{n,h}^p\|  \lesssim
(\tau_n^2 + h^{r+1}) + \|\partial_t 
\vec e_{\tau,h}^{\vec u}(\overline t_n) \| + \|\nabla \vec e_{\tau,h}^{\vec u}(\overline t_n)\| \,.
\end{equation}
Multiplying the squared form of \eqref{Eq:EstP055} with $\tau_n$, using the
inequality of Cauchy--Young, then changing the index from $n$ to $m$ and
summing up the resulting inequality from $m=1$ to $n$ yield  
\begin{equation}
\label{Eq:EstP2}
\sum_{m=1}^n \tau_m \|\overline e_{n,h}^p\|^2  \lesssim t_n 
\big(\tau_n^4 + h^{2(r+1)}\big) + \sum_{m=1}^n 
\tau_m \|\partial_t \vec e_{\tau,h}^{\vec u}(\overline t_m) \|^2 + \sum_{m=1}^n \tau_m \|\nabla 
\vec e_{\tau,h}^{\vec u}(\overline t_m)\|^2 \,.
\end{equation}
Combining \eqref{Eq:EstP2} with the estimates \eqref{Eq:DiscErr000}
  of Lemma~\ref{Lem:DiscErr} then proves the assertion \eqref{Eq:EstP00}.
\end{mproof}

\begin{thm}[Estimate of the pressure error $e_p$]
\label{Thm:EstP}
For the pressure error $e_p =  p  -p_{\tau,h}$, there holds 
\begin{equation}
\label{Eq:EstP}
\sum_{m=1}^n \tau_m \| e_p(\overline t_m)\|^2  \lesssim 
\| \nabla \vec e_{\tau,h}^{\vec u}(0)\|^2 + \tau^4 + h^{2r} 
\end{equation}
for $n=1,\ldots N$. 
\end{thm}

\begin{mproof}
From the splitting \eqref{Eq:SpErr2} of the error $e_p$ we directly conclude the assertion of the theorem by the triangle inequality along with the estimate \eqref{Eq:EstP00} and the estimate 
\begin{align}
\label{Eq:EstP1}
\| \eta_p(t)\| \leq \| p(t) - R_\tau^1 p(t)\| +  \| R_\tau^1 p(t) - R_h^p (R_\tau^1 \vec u(t),R_\tau^1p(t))\| \lesssim \tau_n^2 + h^{r}
\end{align}
for $t\in I_n$, where \eqref{Eq:AppPropR} and \eqref{Eq:StokesApprox2} are used for the second estimate in \eqref{Eq:EstP1}.
\end{mproof}


Similarly to the previous theorem, we deduce from Lemma~\ref{Lem:DiscErr} the following error estimate for the approximation of the velocity field. 

\begin{thm}[Estimate of the error $\vec e_{\vec u}$]
	\label{Thm:Erruprime}
	For the velocity error $\vec e_{\vec u} = \vec u - \vec
      u_{\tau,h}$, there holds  
	\begin{subequations}
		\label{Eq:Err000}	
		\begin{alignat}{2}
			\label{Eq:Err00} 
			\sum_{m=1}^n \tau_m \| \partial_t \vec e_{\vec u}(\overline t_m ) \|^2 + \| 
			\nabla \vec e_{\vec u} (t_n) \|^2 & \lesssim \| \nabla 
                  \vec e_{\vec u} (0) 
			\|^2 + \tau^4 + h^{2r}\,,\\[1ex]
			\label{Eq:Err01}
			\| \vec e_{\vec u} (t_n) \|^2 + \sum_{m=1}^n \tau_m \| \nabla 
                  \vec e_{\vec u} (\overline t_m) \|^2 & \lesssim
                  \| \vec e_{\vec u} (0) \|^2 + 
                  \tau^4 +  h^{2r} 
		\end{alignat}
	\end{subequations}
	for all $n=1,\ldots , N$.
\end{thm}

In view of the error bounds \eqref{Eq:EstP} and  \eqref{Eq:Err000},
it remains to find a suitable discrete initial velocity $\vec u_{0,h}$.
Possible choices, preserving the convergence rates of the estimates,
will be discussed in Section~\ref{Sec:DisIniVal} and Section~\ref{Sec:PPInt}.
We note that \eqref{Eq:Err00} shows a superconvergence by one order in time
of $\partial_t \vec u_{\tau,h}$ in the Gauss quadrature nodes
$\overline{t}_n$, $n=1,\ldots,N$.

\section{Post-processing by collocation}
\label{Sec:PPCol}

In this section, we propose a post-processing of the velocity and
pressure values computed by Problem \ref{Prob:FdS} or \ref{Prob:FdSAF},
respectively, with the aim of defining a piecewise polynomial pressure trajectory
	$\widetilde p_{\tau,h}(t)$
	that is continuous in time on the entire interval
	$[0,T]$ and converges of optimal order with respect to the
	mesh sizes in space and time. For the post-processing, we propose and
	analyze a collocation technique.
	On the interval $I_n=(t_{n-1},t_n]$, we use a collocation condition at
	time $t=t_{n-1}$ in order to compute optimal order approximations 
	$\widetilde p_h^{n-1}$ of $p(t_{n-1})$ and 
	$\vec{a}_h^{n-1}$ of $\partial_t \vec u(t_{n-1})$.
	The post-processing does not change the optimal values 
	$\overline p_h^n = p_{\tau,h}(\overline t_n)$ of the original pressure trajectory
	$p_{\tau,h}(t)$ in the midpoints $\overline t_n$ of the time intervals
	$I_n$
	and preserves also the optimal values of the original velocity trajectory 
	$\vec u_{\tau,h}(t)$ at the time mesh points $t=t_n$.

\subsection{Post-processing of velocity and pressure}

\begin{defi}[{Local post-processing on $I_n$ by collocation}]
	\label{Def:PP}
	
		Let $\vec u_{\tau,h} \in \mathbb P_1(I_n;\vec V_h^{\operatorname{div}})$ and
		$p_{\tau,h} \in \mathbb P_1(I_n;Q_h)$, represented by the expansions
		in \eqref{Eq:Rep_u_p} with  an arbitrarily chosen 
		$p_h^{n-1} = \pth^+(t_{n-1}) \in Q_h$,
		be a solution of Problem \ref{Prob:FdS} on $I_n$. 
		Then, we define the local post-processed solution
		$\vec{\widetilde u}^n_{\tau,h}\in \mathbb P_2(\overline I_n;\vec V_h^{\operatorname{div}})$ 
		and 
		$\widetilde p^n_{\tau,h} \in \mathbb P_1(\overline I_n;Q_h)$ 
		by
	\begin{subequations}
		\label{Eq:PP_up}
		\begin{align}
			\label{Eq:PP_u}
			\vec{\widetilde u}^n_{\tau,h}(t) & \coloneqq  \vec u_{\tau,h}(t) + \vec c_n^{\vec u} \, \vartheta_n(t)\,,\\[1ex]
			\label{Eq:PP_p}
			\widetilde p^n_{\tau,h}(t) & \coloneqq  p_{\tau,h}(t) + d^p_n \, \vartheta_n'(t)
		\end{align}
	\end{subequations}
	for $t\in I_n$. The values of 
		$\vec{\widetilde u}^n_{\tau,h}$, $\widetilde p^n_{\tau,h}$ 
		and  $\partial_t \vec{\widetilde u}^n_{\tau,h}$ at time $t=t_{n-1}$ are
		defined as the right-sided limits for 
        $t\to t_{n-1}+0$ of the corresponding functions with $t\in I_n$. 
	In \eqref{Eq:PP_up}, the polynomial $\vartheta_n\in P_2(\overline I_n;\R)$ 
	is given by 
	\begin{equation}
		\label{Eq:Defvt}
		\vartheta_n(t) \coloneqq  \alpha_n (t-t_{n-1})(t-t_n)\,,	
	\end{equation}
	with the constant $\alpha_n\in \R$ such that
        $\vartheta_n'(t_{n-1})=1$.
	The coefficients 
	$\vec c_n^{\vec u}\in \vec V_h^{\operatorname{div}}$ 
		and $d_n^p \in Q_h$, appearing in~\eqref{Eq:PP_up}, are chosen such
	that at $t= t_{n-1}$ the collocation conditions
	\begin{subequations}
		\label{Eq:ColCond}
		\begin{alignat}{2}	
			\label{Eq:ColCondIp_1}
			\langle \partial_t \vec{\widetilde u}^n_{\tau,h}(t_{n-1}) + A_h
			\vec{\widetilde u}^n_{\tau,h}(t_{n-1}) + B_h' \widetilde p^n_{\tau,h}(t_{n-1}), \vec v_h \rangle & = \langle \vec f(t_{n-1}), \vec v_h \rangle  
			\,, \\[1ex]
			\label{Eq:ColCondIp_2}
			\langle B_h \vec{\widetilde u}^n_{\tau,h}(t_{n-1}), q_h \rangle & = 0 \,,
		\end{alignat}
	\end{subequations}
	hold for all $(\vec v_h,q_h) \in \vec V_h \times Q_h$.
\end{defi}
Since $\vec{\widetilde{u}}^n_{\tau,h}(t_{n-1}) = \vec{u}_{\tau,h}(t_{n-1})$
by \eqref{Eq:PP_u} and \eqref{Eq:Defvt}, condition \eqref{Eq:ColCondIp_2} is
automatically satisfied. 
Using \eqref{Eq:PP_up}, \eqref{Eq:Defvt}, and
	$\vartheta_n'(t_{n-1})=1$, we see that 
	\eqref{Eq:ColCondIp_1}
	is equivalent to
\begin{subequations}
	\label{Eq:PPsys}
	\begin{equation}
		\label{Eq:PPsys1}
		\langle  \vec{c}^{\vec{u}}_{n} 
		+ B_h' d^p_n, \vec v_h \rangle
		= \langle \vec f(t_{n-1})
		- \partial_t \vec{u}_{\tau,h}(t_{n-1})
		- A_h \vec{u}_{\tau,h}(t_{n-1})
		- B_h' p^+_{\tau,h}(t_{n-1})
		, \vec{v}_h \rangle
	\end{equation}
	for all $\vec{v}_h\in \vec V_h$. In order to ensure 
	$\vec{\widetilde u}^n_{\tau,h}\in \mathbb P_2(\overline I_n;\vec V_h^{\operatorname{div}})$, the
	condition
	\begin{equation}
		\label{Eq:PPsys2}
		\langle B_h \vec{c}_n^{\vec{u}} , q_h \rangle = 0
		\qquad\forall q_h\in Q_h
	\end{equation}
\end{subequations}
has to be satisfied. A closer look to \eqref{Eq:PPsys} shows that a
Stokes-like saddle-point problem with the mass matrix instead of the
stiffness matrix has to be solved to obtain the coefficients 
	$\vec{c}^{\vec{u}}_n\in \vec V_h^{\operatorname{div}}$ 
	and $d^p_n \in Q_h$.

\begin{rem}[Uniqueness of the post-processed solution 
  $( \vec{\widetilde u}^n_{\tau,h}, \widetilde p^n_{\tau,h} )$]
	\label{Rem:PPcoll}
   Although the pressure solution $\pth{}_{| I_n}$ of Problem~\ref{Prob:FdS}
   is not unique (see Lemma\ref{Lem:space_time_probl}), we will show that
   the post-processed solution $(\utht^n$, $\ptht^n)$ due to 
   Def.~\ref{Def:PP} is unique in the entire set of solutions.
        The definition of $\vartheta_n$ in \eqref{Eq:Defvt} gives
        \begin{equation}
        \label{Eq:PropTheta}
	\vartheta_n(t_{n-1})=0, \qquad \vartheta_n(t_n) = 0, \qquad \vartheta'_n(\tnb)=0\,.
        \end{equation}
		Hence, it holds 
		\begin{equation}
			\label{Eq:post_up_remains}
			\vec{\widetilde u}^n_{\tau,h}(t_m) = \vec{u}_{\tau,h}(t_m) \,,
			\quad m\in\{n-1, n\}, \qquad
			\widetilde{p}^n_{\tau,h}(\tnb) = p_{\tau,h}(\tnb) \,,
			\qquad 
			\partial_t \vec{\widetilde u}^n_{\tau,h}(\tnb)
				= \partial_t \vec{u}_{\tau,h}(\tnb) \,.
		\end{equation}
	
	By the condition \eqref{Eq:ColCond}, the quantities 
	\begin{subequations}
		\label{Eq:Defahph}
		\begin{alignat}{2}
		\label{Eq:Defahph_a}
		\vec{a}_h^{n-1} & \coloneqq 
		\partial_t \vec{\widetilde u}^n_{\tau,h}(t_{n-1})
		= \partial_t \vec{u}^+_{\tau,h}(t_{n-1}) + \vec{c}^{\vec{u}}_n
		\in  \vec V_h^{\operatorname{div}} \,,\\[1ex]
		\label{Eq:Defahph_b}
		\widetilde{p}_h^{n-1} & \coloneqq 
		\widetilde p^n_{\tau,h}(t_{n-1})  = p^+_{\tau,h}(t_{n-1}) + d^p_n 
		\in Q_h
		\end{alignat} 
	\end{subequations}
	are defined uniquely and independent of the choice of the
      degree of freedom $p^+_{\tau,h}(t_{n-1})$ in the set of
      pressure solutions $\pth$.
		In fact, if we subtract in \eqref{Eq:ColCond}  the known data
		$A_h \vec{\widetilde u}^n_{\tau,h}(t_{n-1}) = A_h \vec{u}_{\tau,h}(t_{n-1})$
                tested with $\vec v_h$
		on both sides of the identity, we get that the pair
		$(\vec{a}_h^{n-1},\widetilde{p}_h^{n-1})
                \in \vec V_h^{\operatorname{div}} \times Q_h$
		solves the problem
		\begin{equation}
			\label{Eq:probl_a_p}
			\langle \vec{a}_h^{n-1}
			+ B_h'  \widetilde{p}_h^{n-1}     , \vec v_h \rangle
			= \langle \vec f(t_{n-1})
			- A_h \vec{u}_{\tau,h}(t_{n-1}) , \vec{v}_h \rangle
			\qquad\forall\, \vec{v}_h\in \vec V_h \,,
		\end{equation}
		where the right-hand side is independent of the choice of
            $p^+_{\tau,h}(t_{n-1})$. 
             Existence and uniqueness of {the solution $( \vec{a}_h^{n-1},\widetilde{p}_h^{n-1}  )$ 
             of~\eqref{Eq:probl_a_p}} can be proved similarly to
            Lemma~\ref{Lem:space_probl} by showing the uniqueness.
      Since the quantities $\widetilde{p}_h^{n-1} = \ptht^n(t_{n-1})$        
      and $\overline p^n_h = \pth(\tnb) = \ptht(\tnb)$ 
      are uniquely determined (cf.\ \eqref{Eq:post_up_remains}), we conclude that the post-processed pressure
       $\widetilde p^n_{\tau,h} \in \mathbb P_1(\overline I_n;Q_h)$, 
      which is given by	\begin{equation}
	  \label{Eq:ptht_represent}
       \ptht^n(t) =   \widetilde{p}_h^{n-1} \varphi_{n,0}(t) 
            +  \bphn\varphi_{n,1}(t)\,, \quad t\in\overline I_n \,,
       \end{equation}
      (see~\eqref{Eq:Rep_u_p})
      is uniquely determined as well.
      Due to~\eqref{Eq:Defahph_a}, the coefficient $\vec{c}^{\vec{u}}_n$ is
      given by
   $\vec{c}^{\vec{u}}_n = \vec{a}_h^{n-1} - \partial_t \vec{u}^+_{\tau,h}(t_{n-1}) $. 
      Its uniqueness follows from the uniqueness of $\vec{a}_h^{n-1}$ 
      and of $\uth{}_{| I_n}$.
      Thus, the definition~\eqref{Eq:PP_u} implies that $\utht$ is
      determined uniquely  and independent of the choice of 
      $p^+_{\tau,h}(t_{n-1})$. 
\end{rem}

A key observation now is that the collocation condition
\eqref{Eq:ColCondIp_1} at the left boundary point $t_{n-1}$ of the
interval $I_n$ implies a corresponding collocation
condition also at the right boundary point $t_n$. This will be proved
in the next lemma.

\begin{lem}[Collocation condition for $t=t_n$]
	\label{Lem:ColCond}
	On $\overline I_n$, let $(\vec{\widetilde u}^n_{\tau,h},\widetilde
	p^n_{\tau,h})\in \mathbb P_2(\overline I_n;\vec V_h^{\operatorname{div}})
	\times \mathbb P_1(\overline I_n;Q_h)$ be defined by \eqref{Eq:PP_up} along
	with \eqref{Eq:ColCond}. Then, for the time $t=t_n$ the collocation conditions
	\begin{subequations}
		\label{Eq:ColCondtn}	
		\begin{alignat}{2}
			\label{Eq:ColCondtn_1}	
			\langle \partial_t \vec{\widetilde u}^n_{\tau,h}(t_{n}) + A_h
			\vec{\widetilde u}^n_{\tau,h}(t_{n}) + B_h' \widetilde p^n_{\tau,h}(t_{n}), \vec v_h \rangle & = \langle \vec f(t_{n}), \vec v_h \rangle  \,,\\[1ex]
			\label{Eq:ColCondtn_2}	
			\langle B_h \vec{\widetilde u}^n_{\tau,h}(t_{n}), q_h \rangle & = 0 
		\end{alignat}
	\end{subequations}
   hold for all $(\vec v_h,q_h) \in \vec V_h\times Q_h$. 
\end{lem}

\begin{mproof}
	We start with proving \eqref{Eq:ColCondtn_1}. 
        The properties~\eqref{Eq:PropTheta} of $\vartheta_n$ and the
        exactness of the Gauss--Lobatto quadrature formula \eqref{Eq:GLF}
        for all polynomials in $\mathbb P_1(I_n;\R)$ provide 
	\begin{align*}
		Q_n^{\operatorname{GL}}\big(\langle \partial_t \vec{\widetilde
			u}^n_{\tau,h}, \vec v_{\tau,h}\rangle\big) & = \int_{I_n}\langle \partial_t
		\vec{\widetilde u}^n_{\tau,h},{\vec v_{\tau,h}}\rangle \ud t =
		\int_{I_n}\langle \partial_t \vec{u_{\tau,h}},{\vec
                v_{\tau,h}}\rangle \ud t
		+ \underbrace{\int_{I_n}\langle \vec c_n^{\vec u}\,
                \vartheta_n',{\vec v_{\tau,h}}\rangle \ud t}_{=Q_n^{\operatorname{G}}(\langle \vec c_n^{\vec
				u}\, \vartheta_n',{\vec v_{\tau,h}} \rangle)= 0}\,,\\[1ex]
		Q_n^{\operatorname{GL}}\big(\langle A_h \vec{\widetilde
			u}^n_{\tau,h},{\vec v_{\tau,h}} \rangle\big) & =
		Q_n^{\operatorname{GL}}\big(\langle A_h \vec{
                u_{\tau,h}},{\vec v_{\tau,h}}
		\rangle\big) + \underbrace{Q_n^{\operatorname{GL}}\big(\langle A_h\vec
			c_n^{\vec{u}} \, \vartheta_n,{\vec v_{\tau,h}} \rangle\big)}_{=0} =
		\int_{I_n}\langle A_h \vec{ u_{\tau,h}},{\vec
                v_{\tau,h}} \rangle \ud t\,,
		\\[1ex] Q_n^{\operatorname{GL}}\big(\langle B_h' \widetilde
		p^n_{\tau,h}, \vec v_{\tau,h}\rangle \big) & = \int_{I_n}\langle B_h'
		p_{\tau,h}, \vec v_{\tau,h}\rangle \ud t +
                \underbrace{\int_{I_n}\langle B_h' d_n^p \, \vartheta_n',
                \vec v_{\tau,h}\rangle \ud
                t}_{{=Q_n^{\operatorname{G}}(\langle B_h' d_n^p \,
                \vartheta_n', \vec v_{\tau,h}\rangle)}=0}
	\end{align*}
	for all $\vec v_{\tau,h} \in \mathbb{P}_0(I_n; \vec V_h)$. By these relations and identity
	\eqref{Eq:FdS2} we then get 
	\begin{equation*}
		Q_n^{\operatorname{GL}}\big(\langle \partial_t \vec{\widetilde
			u}^n_{\tau,h}+A_h \vec{\widetilde u}^n_{\tau,h} + B_h' \widetilde
		p^n_{\tau,h}, \vec v_{\tau,h}\rangle\big) =
		Q_n^{\operatorname{GL}}\big(\langle \vec f ,{\vec
                v_{\tau,h}}\rangle\big) 
	\end{equation*}
	and, therefore, that
	\begin{equation}
		\label{Eq:ColCondtn_4}	
		\sum_{m=n-1}^n \frac{\tau_n}{2} \, \langle \partial_t \vec{\widetilde
			u}^n_{\tau,h}(t_m)+A_h \vec{\widetilde u}^n_{\tau,h}(t_m) + B_h'
		\widetilde p^n_{\tau,h}(t_m) ,\vec v_{h}\rangle = \sum_{m=n-1}^n \frac{\tau_n}{2} \, \langle \vec f(t_m), \vec v_h\rangle
	\end{equation}
	for all $\vec v_h\in \vec V_h$, since $\vec v_{\tau,h}$ takes on
        $I_n$ the constant value $\vec v_h \coloneqq
        \vec v_{\tau,h}(\overline{t}_n)$. Now, from \eqref{Eq:ColCondtn_4} along with \eqref{Eq:ColCondIp_1} we conclude the assertion \eqref{Eq:ColCondtn_1}.	
	
	Since $\vec{\widetilde u}^n_{\tau,h}\in \mathbb P_2(\overline I_n;
	\vec V_h^{\operatorname{div}})$ by construction, condition
	\eqref{Eq:ColCondtn_2} is satisfied.
\end{mproof}

As a consequence of Lemma~\ref{Lem:ColCond}, we get in 
	Lemma~\ref{Lem:ConExt} the result
	that the local post-processed solutions 
	$\partial_t \vec{\widetilde u}^n_{\tau,h}$ and $\widetilde p^n_{\tau,h}$ 
	on the intervals $\overline I_n$ fit continuously together and form  globally
	continuous in time trajectories 
	$\vec{\widetilde u_{\tau,h}}\in X_\tau^2(\vec V_h^{\operatorname{div}})$,
	$\partial_t \vec{\widetilde u_{\tau,h}}\in X_\tau^1(\vec V_h^{\operatorname{div}})$ and 
	$\widetilde p_{\tau,h}\in X_\tau^1(Q_h)$.

\begin{lem}[Global continuity of the $\partial_t \vec{\widetilde u}^n_{\tau,h}$
	and $\widetilde p^n_{\tau,h}$]
	\label{Lem:ConExt}
	On $\overline I_m$, for $m \in \{n,n+1\}$, $n=1, \dots, N-1$,
	let $\vec{\widetilde
		u}^m_{\tau,h}\in \mathbb P_2(\overline I_m;\vec V_h^{\operatorname{div}})$
	and $\widetilde p^m_{\tau,h} \in \mathbb P_1(\overline I_m;Q_h)$ be given by
	Def.~\ref{Def:PP}.
	Then, the continuity conditions
	\begin{subequations}
		\label{Eq:CE}
		\begin{alignat}{2}
			\label{Eq:CE_1}
			\partial_t \vec{\widetilde u}^{n}_{\tau,h}(t_n) & =	\partial_t \vec{\widetilde u}^{n+1}_{\tau,h}(t_n)\,,\\[1ex]
			\label{Eq:CE_2}
			\widetilde p^n_{\tau,h}(t_n) & = \widetilde p^{n+1}_{\tau,h}(t_n)
		\end{alignat}
	\end{subequations}
	hold true for  $n=1, \dots, N-1$.
\end{lem}

\begin{mproof}
	By Lemma~\ref{Lem:ColCond}, the pair $( \vec{\widetilde
		u}^n_{\tau,h},{\widetilde p}^n_{\tau,h})\in P_2(\overline I_n;\vec
	V_h^{\operatorname{div}}) \times \mathbb P_1(\overline I_n;Q_h)$ satisfies the
	collocation conditions \eqref{Eq:ColCondtn} at the endpoint $t_n$
	of $\overline I_n$.  By Def.~\ref{Def:PP}, the pair
	$( \vec{\widetilde u}^{n+1}_{\tau,h},
	{\widetilde p}^{n+1}_{\tau,h})\in P_2(\overline I_{n+1};
	\vec V_h^{\operatorname{div}}) \times \mathbb
	P_1(\overline I_{n+1};Q_h)$ fulfills
	the collocation conditions \eqref{Eq:ColCond} at the initial time point
	$t_n$ of $\overline I_{n+1}$. Thus,  for the differences
	\begin{equation*}
		\begin{aligned}
			\vec a_h & \coloneqq  
			\partial_t \vec{\widetilde u}^{n}_{\tau,h}(t_n) - \partial_t \vec{\widetilde u}^{n+1}_{\tau,h}(t_n)
			\,,\\[1ex]
			s_h   & \coloneqq  
			\widetilde p^n_{\tau,h}(t_n) - \widetilde p^{n+1}_{\tau,h}(t_n) \,,
		\end{aligned}
	\end{equation*}
	we conclude from the difference of the identities \eqref{Eq:ColCond}
	and  \eqref{Eq:ColCondtn} along with \eqref{Eq:PP_up} and
	\eqref{Eq:Defvt} that 
	\begin{subequations}
		\label{Eq:ColCondDif}	
		\begin{alignat}{2}
			\label{Eq:ColCondDif_1}	
			\langle \vec a_h + B_h' s_h, \vec v_h \rangle & = 0 \,,\\[1ex]
			\label{Eq:ColCondDif_2}	
			\langle B_h \vec a_h, q_h \rangle & = 0 
		\end{alignat}
	\end{subequations}
      for all $(\vec v_h,q_h) \in \vec V_h\times Q_h$. Then, the unique
      solvability of \eqref{Eq:ColCondDif}
      (cf. Remark~\ref{Rem:PPcoll}) proves the assertion.
\end{mproof}

\begin{rem}[Efficient post-processing on the subsequent interval $\overline I_{n+1}$]
Let us assume that the post-pro\-cessed solution $(\utht^n , \ptht^n)$ on the
interval $\overline I_n$ has been already computed. Then, the
post-processed solution on $\overline I_{n+1}$ can be determined in a
computationally cheaper way without solving a system corresponding
to~\eqref{Eq:PPsys} for the coefficients $\vec c_{n+1}^{\vec u}$ and
$d_{n+1}^{p}$. This follows from the following observation. On $\overline
I_{n+1}$, let $\vec{\widetilde u}^{n+1}_{\tau,h}\in \mathbb P_2(\overline
I_{n+1};\vec V_h^{\operatorname{div}})$ and $\widetilde p^{n+1}_{\tau,h}
\in \mathbb P_1(\overline I_{n+1};Q_h)$ be given by
	\begin{subequations}
		\label{Eq:PP_upe}
		\begin{align}
			\label{Eq:PP_ue}
			\vec{\widetilde u}^{n+1}_{\tau,h}(t) & = \vec
			u_{\tau,h}(t) + \vec c_{n+1}^{\vec u} \,
			\vartheta_{n+1}(t)   \,,\\[1ex]
			\label{Eq:PP_pe}
			\widetilde p^{n+1}_{\tau,h}(t) & = p_{\tau,h}(t) +
			d^p_{n+1} \, \vartheta_{n+1}'(t) \,,
		\end{align}
	\end{subequations}
	with coefficients 
      $(\vec c_{n+1}^{\vec u},d^p_{n+1} ) \in \vec V_h^{\operatorname{div}}\times Q_h$. 
	By $\vartheta_{n+1}'(t_n)=1$ and the continuity conditions \eqref{Eq:CE},  
      we obtain from~\eqref{Eq:PP_upe} that  
	\begin{subequations}
		\label{Eq:CE_3}
		\begin{alignat}{2}
			\label{Eq:CE_4}
			\vec c_{n+1}^{\vec u} & =
			\partial_t \vec{\widetilde{u}}^n_{\tau,h}(t_n) -
                        \partial_t \vec{u}_{\tau,h}^+(t_n) \,,\\[1ex]
			d_{n+1}^{p} & = \widetilde{p}^n_{\tau,h}(t_n) -
                        p_{\tau,h}^{+}(t_n)   \,.
		\end{alignat}
	\end{subequations}
As already noted before in Remark~\ref{Rem:PPcoll}, the coefficient
$d_{n+1}^{p}$ changes with the choice of degree of freedom $\pth^+(t_n)$ in
the set of the pressure solutions $\pth$ on $I_{n+1}$. However, the
resulting post-processed pressure solution $\ptht^{n+1}$ is independent of
the choice of $\pth^+(t_n)$ and can be computed also by means of
\begin{equation}
\ptht^{n+1}(t) = \widetilde{p}_h^n  \varphi_{n+1,0}(t) + \overline p_h^{n+1}
\varphi_{n+1,1}(t)\,, \quad t\in I_{n+1}\,,	
\end{equation}
where 
$\widetilde{p}_h^n \coloneqq \ptht^{n}(t_n) = \ptht^{n+1}(t_n)$
and $\overline p_h^{n+1} \coloneqq \pth(\overline t_{n+1}) =
\ptht^{n+1}(\overline t_{n+1}) $.
\end{rem}

We will comment on the implementation of the post-processing in the first subinterval. 

\begin{rem}[Post-processing on interval $\overline I_1$]
	\label{Rem:PPcoll_one}
  We distinguish two cases.
  The first case is that a discrete initial value
  $\vec u_{0,h} \in \Vhdiv$ for $\vec u_0$ is available, such that $\vec u_{0,h} \in \Vhdiv$ admits optimal order (spatial) approximation properties. Then, we compute  
	 $\vec{\widetilde u}^{1}_{\tau,h}\in \mathbb P_2(\overline
	                       I_{1};\vec V_h^{\operatorname{div}})$ 
      from the definition~\eqref{Eq:PP_u}
      with the coefficient  (cf.\ \eqref{Eq:Defahph_a})
	\begin{equation}
        \label{Eq:cu_one}
	  \vec c_{1}^{\vec u} \coloneqq \vec{a}_h^{0} 
		- \partial_t \vec{u}^+_{\tau,h}(0) 
	\end{equation}
      and $\widetilde p^{1}_{\tau,h} \in \P_1(\overline I_{1};Q_h)$ from
	\begin{equation}
		\label{Eq:ptht_one}
        \ptht^{1}(t) = \widetilde{p}_h^0  \varphi_{1,0}(t) + \overline p_h^{1}
        \varphi_{1,1}(t) \,, \quad t\in \overline I_{1}\,,	
	\end{equation}
      with 
      $\,\overline p_h^{1} \coloneqq \pth(\overline t_{1})$
      and 
      $(\vec{a}_h^{0} , \widetilde{p}_h^{0} ) \in \vec V_h^{\operatorname{div}} \times Q_h$
	as solution of problem~\eqref{Eq:probl_a_p} for $n=1$ and
      $\uth(t_0) \coloneqq \vec u_{0,h}$.

      In the second case without any discrete initial value $\vec u_{0,h} \in \Vhdiv$, we first determine an approximation 
      $ \vec{a}_{0,h}  \in \Vhdiv$ of $\vec a_0 = \partial_t\vec u (0)$
      by solving the Stokes problem~\eqref{Eq:a0hsh3} (cf.\ Section~\ref{Sec:DisIniVal}
      below). Using this approximation $\vec{a}_{0,h}$, we compute the
      discrete initial velocity and pressure
            $(\vec u_{0,h},p_{0,h})\in \Vhdiv\times Q_h$
      as solution of the Stokes problem~\eqref{Eq:u0hp0h1}.
      Then, we put 
      $(\vec{a}_h^{0} , \widetilde{p}_h^{0} ) \coloneqq ( \vec a_{0,h},p_{0,h} ) $. 
      This pair solves \eqref{Eq:probl_a_p} for $n=1$ and
      with $\uth(t_0) \coloneqq \vec u_{0,h}$, similarly to the first case.
      Therefore, we again get the 
      post-processed velocity $\utht^{1}$ by~\eqref{Eq:PP_u},
      with the coefficient  $\vec c_{1}^{\vec u}$ from~\eqref{Eq:cu_one}
      and the post-processed pressure $\ptht^{1}$ by~\eqref{Eq:ptht_one}.
\end{rem}

Due to the continuity properties of the local post-processed trajectories 
	$\vec{\widetilde u}^n_{\tau,h}(t)$, $\partial_t \vec{\widetilde u}^n_{\tau,h}(t)$, 
	and $\widetilde p^n_{\tau,h}(t)$  given in Lemma~\ref{Lem:ConExt}, we will
	omit in the following the upper index $n$ indicating the interval
	$I_n$ where these local trajectories have been computed.
	For a given time $t\in I$, this index $n$ is determined by the
	condition $t\in I_n$. The setting is completed by choosing
        $n=1$ for $t=0$.
	Thus, we will write in the following for $t\in \overline I$ simply 
	$\vec{\widetilde u}_{\tau,h}(t)$, $\partial_t \vec{\widetilde u}_{\tau,h}(t)$,
	and $\widetilde p_{\tau,h}(t)$ 
	to denote the global post-processed solution trajectories.

\subsection{Preparation of the error analysis for the post-processed variables}	

Here we present some auxiliary results that will be used below in
Section~\ref{Subsec:ErrPostPresCol} for the error estimation. Firstly,
we derive an operator equation that is satisfied by the pair
$(\vec{\widetilde u}_{\tau,h},\widetilde p_{\tau,h})$ for $t\in
\overline I_n$. In particular, at $t=\overline{t}_n$ a collocation
condition for the time derivative $(\partial_t \vec{\widetilde u}_{\tau,h},
\partial_t \widetilde p_{\tau,h})$ of the post-processed discrete
solution $(\vec{\widetilde u}_{\tau,h}, \widetilde p_{\tau,h})$
is thus obtained. Since $\vec{\widetilde u}_{\tau,h}(t) =
\vec{u}_{\tau,h}(t)$ for $t\in \{t_{n-1},t_n\}$ by means of
\eqref{Eq:post_up_remains},
we recast the collocation
conditions \eqref{Eq:ColCond} and \eqref{Eq:ColCondtn} as 
\begin{equation}
	\label{Eq:Aux01}
	\partial_t \vec{\widetilde u}_{\tau,h}(t) + A_h \vec{u}_{\tau,h}(t) + B_h' \widetilde p_{\tau,h}(t) = P_h I_{n,1}^{\operatorname{GL}}\vec f(t)\,, \quad t\in \{t_{n-1},t_n\}\,,
\end{equation}
where $I_{n,1}^{\operatorname{GL}}$  denotes
	the local Gauss--Lobatto interpolation operator defined in \eqref{Eq:LocLagIntOp1}
	and $P_h: \vec L^2(\Omega) \to \vec V_h$ is the $L^2$-projection
        operator into $\vec V_h$ defined by
	\begin{equation}
		\label{Eq:PhDef}
		\langle  P_h \vec{w}, \vec v_h \rangle  = \langle \vec w , \vec v_{h}\rangle
		\quad\forall \vec v_h\in \vec V_h \,.
	\end{equation}
Since both sides of \eqref{Eq:Aux01} are elements of $\mathbb P_1(I_n; \vec
V_h)$, we obtain 
\begin{equation}
	\label{Eq:Aux02}
	\partial_t \vec{\widetilde u}_{\tau,h}(t) + A_h \vec{u}_{\tau,h}(t) + B_h' \widetilde p_{\tau,h}(t) = P_h I_{n,1}^{\operatorname{GL}}\vec f(t)\,, \quad t\in \overline I_n\,.
\end{equation} 

Next, we derive a collocation condition at $t=\overline t_n$ that is
satisfied by the time derivatives $(\partial_t \vec{\widetilde
	u}_{\tau,h},\partial_t \widetilde p_{\tau,h})$.  For this, we need to
introduce a further interpolant of mixed Lagrange and Hermite type. For
$\vec f \in C^1(\overline I;\vec L^2(\Omega))$, we define the local
interpolant  $L_{n} \vec f \in \mathbb P_2(\overline I_n;\vec
L^2(\Omega))$ by 
\begin{equation}
	\label{Eq:DefLtau_01}
	L_n \vec f (t)= \vec f(t)\,, \quad t\in \{t_{n-1},t_n\}\,, \quad \partial_t L_n  \vec f (t_{n-1}) = \partial_t \vec f (t_{n-1})\,,
\end{equation}	
and the global interpolant $L_\tau \vec f \in X_\tau^2(\vec L^2(\Omega))$ 
by the initial value $L_\tau \vec f(0)\coloneqq \vec f(0)$ and 
\begin{equation*}
	\label{Eq:DefLtau_03}
	L_\tau \vec f (t) \coloneqq  L_n \vec f(t)\,, \quad t \in I_n\,, \quad n=1,\ldots ,N\,. 
\end{equation*}
From \eqref{Eq:DefLtau_01} and the global continuity of $L_\tau \vec f$  on $\overline I$ we get  
\begin{equation*}
	\label{Eq:DefLtau_04}
	L_\tau \vec f (t) - I_1^{\operatorname{GL}}\vec f (t) = 0 \,, \quad t \in \{t_{n-1},t_n\}\,.  
\end{equation*}
Therefore, the local interpolant $L_n\vec f \in \mathbb P_2(I_n;\vec L^2(\Omega))$ admits the representation 
\begin{equation}
	\label{Eq:DefLtau_02}
	L_n \vec f(t) = I_{n,1}^{\operatorname{GL}}\vec f
        (t) + \vec b_n^{\vec f}\, \vartheta_n(t)\,, \quad t\in I_n\,,
\end{equation}

where the coefficient $\vec b_n^{\vec f}\in \vec L^2(\Omega)$ is chosen
such that the second condition in \eqref{Eq:DefLtau_01} is
satisfied. For the  interpolant $L_\tau \vec f$ of mixed Lagrange
and Hermite type, the error estimate
\begin{equation}
	\label{Eq:DefLtau_05}
	\| \partial_t (L_n \vec f - \vec f) \|_{C(\overline I_n;\vec L^2(\Omega))} \lesssim \tau_n^2 \| \partial_t^{3} \vec f\|_{C(\overline I_n;\vec L^2(\Omega))} 
\end{equation}
holds if $\vec f$ is sufficiently regular. Its proof follows the
lines of \cite[Lemma~5.5]{ES16}.

\begin{lem}
	For $(\partial_t \vec{\widetilde u}_{\tau,h},\partial_t \widetilde
      p_{\tau,h})\in \mathbb P_1(I_n;\vec V_h^{\operatorname{div}})
      \times \mathbb P_0(I_n;Q_h)$,  there holds 
	\begin{equation}
		\label{Eq:Aux05}
		\langle \partial_t (\partial_t \vec{\widetilde u}_{\tau,h})(\overline t_n) + A_h (\partial_t \vec{\widetilde u}_{\tau,h})(\overline t_n) + B_h' (\partial_t \widetilde p_{\tau,h})(\overline t_n) , \vec v_h\rangle = \langle \partial_t  L_n \vec f(\overline t_n), \vec v_{h} \rangle 
	\end{equation}
	for all $\vec v_h \in \vec V_h$. 
\end{lem}

\begin{mproof}
	Taking the time derivative of \eqref{Eq:Aux02}, we get   
	\begin{equation}
		\label{Eq:Aux025}
		\partial_t \partial_t \vec{\widetilde u}_{\tau,h}(t) + A_h \partial_t \vec{u}_{\tau,h}(t) + B_h' \partial_t \widetilde p_{\tau,h}(t) = P_h \partial_t I_{n,1}^{\operatorname{GL}}\vec f(t)\,, \quad t\in (t_{n-1},t_n)\,.
	\end{equation} 
	By \eqref{Eq:PP_u} there holds 
	\begin{equation*}
		\partial_t \vec{\widetilde u}_{\tau,h}(t) = \partial_t
                \vec{ u}_{\tau,h}(t)+ \vec c_n^{\vec{ u}} \,\vartheta_n'(t)
	\end{equation*}
	with some coefficient $ \vec c_n^{\vec{ u}}\in  \vec
        V_h^{\operatorname{div}}$. The properties~\eqref{Eq:PropTheta}
        provide 
	\begin{equation}
		\label{Eq:Aux03}
		\int_{I_n} \langle A_h \partial_t \vec{\widetilde u}_{\tau,h}, \vec v_{\tau,h}\rangle \ud t = \int_{I_n} \langle A_h \partial_t \vec{u}_{\tau,h}, \vec v_{\tau,h}\rangle \ud t +  \underbrace{\int_{I_n} \langle \vec c_n^{\vec{u}}\, \vartheta_n', \vec v_{\tau,h}\rangle \ud t}_{= Q_n^{\operatorname{G}}( \langle \vec c_n^{\vec{u}} \,\vartheta_n', \vec v_{\tau,h}\rangle )=0}
	\end{equation}
	for all $\vec v_{\tau,h}\in \mathbb P_0 (I_n;\vec V_h)$. Using \eqref{Eq:Aux025} along with \eqref{Eq:Aux03} yields 
	\begin{equation}
		\label{Eq:Aux04}
		\int_{I_n} \langle \partial_t (\partial_t \vec{\widetilde u}_{\tau,h}) + A_h (\partial_t \vec{\widetilde u}_{\tau,h}) + B_h' (\partial_t \widetilde p_{\tau,h}), \vec v_{\tau,h} \rangle \ud t = \int_{I_n} \langle \partial_t I_{n,1}^{\operatorname{GL}}\vec f, \vec v_{\tau,h} \rangle \ud t
	\end{equation}
	for all $\vec v_{\tau,h}\in \mathbb P_0 (I_n;\vec V_h)$.
	
	For $L_n \vec f$, we have
	\begin{equation}
		\label{Eq:DefLtau_10}
		\int_{I_n} \langle \partial_t L_n \vec f, \vec v_{\tau,h} \rangle \ud t = \int_{I_n} \langle \partial_t I_{n,1}^{\operatorname{GL}} \vec f, \vec v_{\tau,h} \rangle \ud t+
		\underbrace{\int_{I_n}  \langle \vec b_n^{\vec f}\, \vartheta_n', \vec v_{\tau,h} \rangle \ud t}_{=Q_n^{\operatorname{G}}(\langle \vec b_n^{\vec f} \,\vartheta_n', \vec v_{\tau,h} \rangle )=0} 	
	\end{equation}
	for all $\vec v_{\tau,h} \in \mathbb P_0(I_n;\vec V_h)$. From \eqref{Eq:Aux04} and \eqref{Eq:DefLtau_10}  we then deduce  
	\begin{equation}
		\label{Eq:Aux045}
		\int_{I_n} \langle \partial_t (\partial_t \vec{\widetilde u}_{\tau,h}) + A_h (\partial_t \vec{\widetilde u}_{\tau,h}) + B_h' (\partial_t \widetilde p_{\tau,h}), \vec v_{\tau,h} \rangle \ud t = \int_{I_n} \langle \partial_t L_n \vec f, \vec v_{\tau,h} \rangle \ud t 
	\end{equation}
	for all $\vec v_{\tau,h} \in \mathbb P_0(I_n;\vec V_h)$. All time integrals in \eqref{Eq:Aux045} are computed exactly by the Gauss formula \eqref{Eq:GF}. Putting $\vec v_h \coloneqq \vec v_{\tau,h}(\overline t_n)$, we then conclude from \eqref{Eq:Aux045} the assertion \eqref{Eq:Aux05}.
\end{mproof}

\subsection{Error estimates for the post-processed variables}	
\label{Subsec:ErrPostPresCol}
Here, we derive error estimates for the time derivative $(\partial_t
\vec{\widetilde u}_{\tau,h}, \partial_t \widetilde p_{\tau,h})$ of the
post-processed discrete solution $(\vec{\widetilde u}_{\tau,h}, \widetilde
p_{\tau,h})$.
For brevity, we put
\begin{equation}
	\label{Eq:DefEup}
	\vec E_{\vec u} \coloneqq  \partial_t \vec u - \partial_t \vec{\widetilde u}_{\tau,h}\qquad \text{and} \qquad 
	E_p  \coloneqq   \partial_t p- \partial_t \widetilde p_{\tau,h}\,.
\end{equation}
We split the errors $\vec E^{\vec u}$ and $E_n^p $ of \eqref{Eq:DefEup} by means of 
\begin{subequations}
	\label{Eq:DefErSplitp}
	\begin{alignat}{3}
		\label{Eq:DefErSplitpu}
		\vec E_{\vec u}& =  \vec{\widehat \eta^{\vec u}} + \vec E^{\vec u}_{\tau,h}\,, & \qquad
		\vec{\widehat \eta^{\vec u}} & \coloneqq  \partial_t \vec u - \vec{\widehat U}_{\tau,h}\,, & \quad
		\vec E^{\vec u}_{\tau,h} & \coloneqq  \vec{\widehat U}_{\tau,h}-\partial_t \vec{\widetilde u_{\tau,h}}\,,\\[1ex]
		\label{Eq:DefErSplitpp}
		E_p &  = \widehat{\eta}^p + E^p_{\tau,h}\,, & \qquad
		\widehat{\eta}^p & \coloneqq  \partial_t p- \widehat{P}_{\tau,h}\,, & \quad
		E^p_{\tau,h} & \coloneqq  \widehat{P}_{\tau,h} - \partial_t \widetilde p_{\tau,h}\,,
	\end{alignat}
\end{subequations}
%
%
with $\vec{\widehat U}_{\tau,h}\in
X_\tau^1(\vec V_h^{\operatorname{div}})$ and $\widehat{P}_{\tau,h}\in
X_\tau^1(Q_h)$ defined by
\begin{equation}
	\label{Eq:DefEup_01}
	 \vec{\widehat U}_{\tau,h}\coloneqq  R_h^{\vec u}( R_\tau^1 \partial_t\vec u,   R_\tau^1 \partial_t p  )
		\qquad \text{and} \qquad 
		\widehat{P}_{\tau,h}\coloneqq  R_h^{p}( R_\tau^1 \partial_t \vec u,  R_\tau^1 \partial_t p  ) \,,  
\end{equation}
where $R_h^{\vec u}$ and $R_h^{p}$ are the velocity and pressure components of the Stokes projection $R_h$ defined in  Def.~\ref{Def:StokesProj} and $R_\tau^1$ is the interpolant in time given in~\eqref{Eq:IntOpR}.

Firstly, we have the following result for $\vec E^{\vec u}_{\tau,h}$. 

\begin{lem}
	\label{Lem:DiscErrp} 
	For the discrete velocity error $\vec E^{\vec u}_{\tau,h}$, defined in \eqref{Eq:DefErSplitpu}, there holds 
	\begin{subequations}
		\label{Eq:ErrUp}
		\begin{align}
			\label{Eq:ErrUp_01}
			\sum_{m=1}^n \tau_m \| \partial_t \vec E^{\vec u}_{\tau,h} (\overline
			t_m)\|^2 + \| \nabla \vec E^{\vec u}_{\tau,h} (t_n)\|^2 & 
			\lesssim  
			\| \nabla \vec E^{\vec u}_{\tau,h}(0)\|^2 
			+ \tau^4 + h^{{2(r+1)}}\,,\\[1ex] 
			\label{Eq:ErrUp_02}
			\| \vec E^{\vec u}_{\tau,h} (t_n)\|^2  + \sum_{m=1}^n \tau_m \| \nabla
			\vec E^{\vec u}_{\tau,h} (\overline t_m)\|^2 &  
			\lesssim  
			\| \vec E^{\vec u}_{\tau,h}(0)\|^2 
                        + \tau^4 + h^{2(r+1)}
		\end{align}
	\end{subequations}
	for all $n=1,\ldots, N$.
\end{lem}

\begin{mproof}
	For a sufficiently regular solution $(\vec u,p)$ of \eqref{Eq:SE}, the time derivative of \eqref{Eq:SE_1} is satisfied for all $t\in \overline I$. Together with \eqref{Eq:Aux05}, this yields the error equation 
	\begin{equation}
		\label{Eq:EPP_01}
		Q_n^{\operatorname G}\big(\langle \partial_t \vec E_{\vec
                u} + A_h \vec E_{\vec u} + B_h' E_p, \vec
                v_{\tau,h}\rangle\big) = Q_n^{\operatorname G}\big(\langle
                \partial_t \vec f,{\vec v_{\tau,h}}\rangle  \big)-
                Q_n^{\operatorname G}\big(\langle \partial_t L_n \vec
                f,{\vec v_{\tau,h}}\rangle \big) 
	\end{equation}
	for all $\vec v_{\tau,h} \in \mathbb{P}_0(I_n; \vec V_h)$. The splitting \eqref{Eq:DefErSplitp} of the errors then implies  
	\begin{subequations}
		\label{Eq:EPP_02}
		\begin{equation}
			Q_n^{\operatorname G}\big(\langle \partial_t \vec
                        E_{\tau,h}^{\vec u}  + A_h \vec E_{\tau,h}^{\vec u}
                        + B_h' E_{\tau,h}^p, \vec v_{\tau,h}\rangle \big)
			= T_n^{\partial_t \vec f} + T_n^{ \vec{\widehat \eta}}\,,
		\end{equation}
		where
		\begin{equation}
			T_n^{\partial_t \vec f} \coloneqq  \tau_n  \langle \vec{\widehat \delta}_n^{\vec f} ,\vec v_h \rangle\,,\qquad
			T_n^{ \vec{\widehat \eta}} \coloneqq  - Q_n^{\operatorname G}\big( \langle \partial_t \vec{\widehat \eta}^{\vec u} + A_h \vec{\widehat \eta}^{\vec u}+ B_h'
                        \widehat \eta^p , \vec v_{\tau,h} \rangle  \big)
		\end{equation}
	\end{subequations}%
	%
        with $\vec v_h \coloneqq \vec v_{\tau,h}(\overline{t}_n)$ and $\vec{\widehat \delta}_n^{\vec f} \coloneqq  \partial_t \vec f(\overline t_n) - \partial_t L_n \vec f(\overline t_n)$.  By \eqref{Eq:DefLtau_05} we obtain 
	\begin{equation}
		\label{Eq:EPP_04}
		|   T_n^{\partial_t \vec f}  |   \lesssim  \tau_n \| \partial_t (\vec f - L_n \vec f)\|_{C(\overline I_n;\vec L^2(\Omega))}   \| \vec v_h\| \lesssim \tau_n \tau_n^2\| \vec v_h\|\,.
	\end{equation}
	Similarly to Lemma~\ref{Lem:ApproxErr}, it can be shown
	\begin{equation}
		\label{Eq:EPP_05}
		|   T_n^{ \vec{\widehat \eta}}  |  \lesssim \tau_n (\tau_n^2 + h^{r+1}) \| \vec v_h\| \,.
	\end{equation}
	Following the lines of Lemma~\ref{Lem:DiscErr}, we deduce from \eqref{Eq:EPP_02} to \eqref{Eq:EPP_05} the assertion \eqref{Eq:ErrUp}.

\end{mproof}


It remains to provide suitable bounds for norms of $\vec E_{\tau,h}^{\vec u}(0)$. 
This will be addressed in Section~\ref{Sec:DisIniVal}.

\begin{lem}
	\label{Lem:DiscErrpp} 
	For the discrete pressure error $E^p_{\tau,h}$, defined in \eqref{Eq:DefErSplitpp}, there holds 
	\begin{align}
		\label{Eq:Errpp_20}
		\sum_{n=1}^{N} \tau_n  \| E^p_{\tau,h}(\overline t_n) \|^2 &\lesssim
                 \| \nabla \vec E_{\tau,h}^{\vec u}(0)\|^2
		+  \tau^4 + h^{{2(r+1)}} \,.
	\end{align}
\end{lem}

\begin{mproof}
	We follow the lines of the proof for Lemma~\ref{Lem:EstPth}. By the choice of an inf-sup stable pair of finite element spaces $\vec V_h \times Q_h$, fulfilling the Assumption~\ref{Ass:FES},  for $\overline E_{n,h}^{p}\coloneqq  E^p_{\tau,h}(\overline t_n) \in Q_h$ there exists some 
	function $\vec v_h \in \vec V_h\setminus \{\vec 0\}$ such that 
	\begin{equation}
		\label{Eq:Errpp_21}	
                \beta \| \overline E_{n,h}^{p}  \| \le
		\dfrac{\langle B_h'  \overline E_{n,h}^{p} , \vec v_h \rangle }%
		{\| \vec v_h \|_1}  \,,
	\end{equation}
	where the constant $\beta>0$ is independent of $h$. Combining \eqref{Eq:Errpp_21} with the error equation \eqref{Eq:EPP_02} yields  
	\begin{equation}
		\label{Eq:Errpp_22}	
		\| \overline E_{n,h}^{p}  \|  \lesssim \frac{1}{\| \vec v_h  \|_1 } \big(|\langle \partial_t \vec E_{\tau,h}^{\vec u} (\overline t_n) + A_h \vec E_{\tau,h}^{\vec u}(\overline t_n),
		\vec v_h\rangle| +\tau_n^{-1}|T_n^{\partial_t \vec f}| + \tau_n^{-1}|T_n^{ \vec{\widehat \eta}}|\big)\,.
	\end{equation}
	By the inequality of Cauchy--Schwarz along with \eqref{Eq:EPP_04} and \eqref{Eq:EPP_05}, we deduce from \eqref{Eq:Errpp_22} that 
	\begin{equation}
		\label{Eq:Errpp_23}
		\| \overline E_{n,h}^{p}  \|  \lesssim \tau_n^2 + h^{r+1}+ \| \partial_t
			\vec E_{\tau,h}^{\vec u} (\overline t_n) \|  + \| \nabla  \vec E_{\tau,h}^{\vec u}(\overline t_n)\|\,.
	\end{equation}
	Squaring \eqref{Eq:Errpp_23}, multiplying the resulting inequality with
	$\tau_n$, summing up from
	$n=1$ to $N$ and, finally, applying Lemma \ref{Lem:DiscErrp} yield the assertion \eqref{Eq:Errpp_20}.
\end{mproof}	

We are now able to prove the desired error estimate for the
post-processed pressure $\widetilde p_{\tau,h}$, with $\widetilde
p_{\tau,h}{}_{| I_n}\in \mathbb P_1(I_n;Q_h)$, that is defined by \eqref{Eq:PP_p} along with \eqref{Eq:ColCond}.

\begin{thm}[Error estimate for the post-processed pressure]
	\label{Thm:ErrPPP}
	For the error of the post-pro\-cessed pressure $\widetilde
      p_{\tau,h}$, defined by \eqref{Eq:PP_p} and \eqref{Eq:ColCond},
      and $\tau=\max_{n=1,\ldots,N} \tau_n$, there holds
	\begin{equation}
		\label{Eq:ErrPPP_00}
		\left(\int_0^T 	\| p(t) - \widetilde p_{\tau,h}(t) \|^2 \ud t\right)^{1/2} \lesssim
                \|\nabla \euth(0)\| + \tau \| \nabla \vec E_{\tau,h}^{\vec u} (0)\|
                +  \tau^2 + h^r\,.
	\end{equation}
\end{thm}

\begin{mproof}
	By the definition of $E^p_{\tau,h}$ in \eqref{Eq:DefErSplitpp} there holds 
	\begin{align}
		\label{Eq:ErrPPP_005}
		E^p_{\tau,h}(\overline t_n) & = \widehat P_{\tau,h}	(\overline t_n) -
		\partial_t \widetilde p_{\tau,h}(\overline t_n)
		= R_h^p(R_\tau^1\partial_t \vec u(\overline t_n), R_\tau^1 \partial_t
		p(\overline t_n) )- \partial_t \widetilde
		p_{\tau,h}(\overline t_n)\,.
	\end{align}
	Setting $\overline E_{n,h}^{p}\coloneqq  E^p_{\tau,h}(\overline t_n)
	$, it follows from \eqref{Eq:ErrPPP_005} along with \eqref{Eq:AppPropR} and \eqref{Eq:StokesApprox0} that 
        \begin{align}
          \| \partial_t p(\overline t_n) -  \partial_t \widetilde p_{\tau,h}(\overline t_n) \|
		& \le  \| \partial_t p(\overline t_n) -  R_\tau^1 \partial_t p (\overline t_n) \|  
		+ \| R_\tau^1 \partial_t p (\overline t_n)  - R_h^p(R_\tau^1\partial_t \vec u(\overline t_n)  , R_\tau^1 \partial_t p (\overline t_n)  ) \| 
		+ \| \overline E_{n,h}^{p}  \| \nonumber \\[1ex] 	
		& \lesssim \tau_n^2 + h^{r} + \| \overline E_{n,h}^{p} \| \,.
		\label{Eq:ErrPPP_01}
        \end{align}
	
	Let now $p_\tau \in \mathbb P_1(I_n;Q)$ denote the linear Taylor polynomial
	\begin{equation}
		\label{Eq:ErrPPP_02}
		p_\tau (t) \coloneqq p(\overline t_n) + \partial_t p(\overline t_n) (t-\overline t_n)\,, \quad t\in I_n\,.
	\end{equation}
	 From \eqref{Eq:ErrPPP_02} we directly get  
		\begin{equation}
			\label{Eq:ErrPPP_05}
			\| p(t) - p_{\tau}(t)\| \lesssim \tau_n^2\,, \quad t\in I_n\,.
		\end{equation}
	
	Since $\widetilde p_{\tau,h}{}_{| I_n}\in \mathbb P_1(I_n;Q_h)$ and $ \widetilde p_{\tau,h}(\overline t_n) = p_{\tau,h}(\overline t_n) $ by \eqref{Eq:PP_p}, we have  
	\begin{equation}
		\label{Eq:ErrPPP_03}
		p_\tau (t) - \widetilde p_{\tau,h}(t) = \underbrace{p (\overline t_n) - \widetilde p_{\tau,h}(\overline t_n)}_{= e_p(\overline t_n)} 
		+ \big(\partial_t p(\overline t_n) - \partial_t \widetilde p_{\tau,h}(\overline t_n)\big) (t-\overline t_n)\,.
	\end{equation}
	 We conclude from \eqref{Eq:ErrPPP_03} along with \eqref{Eq:ErrPPP_01} that 
	\begin{equation}
		\label{Eq:ErrPPP_04}
		\| p_\tau (t_n) - \widetilde p_{\tau,h}(t_n)  \| \lesssim
                \| e_p(\overline t_n) \| + \tau_n \big(\tau_n^2 + h^{r} +
                \| \overline E_{n,h}^{p}  \|\big)\,.  
	\end{equation}
	From \eqref{Eq:ErrPPP_04} we find  
	\begin{equation}
		\label{Eq:ErrPPP_06}
		\tau_n \| p_\tau (t_n) - \widetilde p_{\tau,h}(t_n)  \|^2
		\lesssim   \tau_n \tau_n^2 \big( \tau_n^4 + h^{2r} \big) 
		+ \tau_n  \| e_p(\overline t_n) \|^2 + \tau_n \tau_n^2 \| \overline E_{n,h}^{p}  \|^2\,.
	\end{equation}
	Now, we sum up inequality \eqref{Eq:ErrPPP_06} from $n=1$ to $N$. 
	We apply Theorem~\ref{Thm:EstP} and Lemma~\ref{Lem:DiscErrpp} to obtain 
	\begin{align}
		\label{Eq:ErrPPP_07}
		\sum_{n=1}^{N} \tau_n \| p_\tau (t_n) - \widetilde p_{\tau,h}(t_n) \|^2 & \lesssim \tau^4 + h^{2r}
                + \|\nabla \euth(0)\|^2 + \tau^2 \| \nabla \vec E_{\tau,h}^{\vec u} (0)\|^2
                \,. 
	\end{align}
	Exploiting the fact that $p_\tau, \widetilde p_{\tau,h}{}_{| I_n}\in \mathbb
	P_1(I_n;Q)$, we can easily show the estimate
	\[
	\int_{I_n} 	\| p_\tau(t) - \widetilde p_{\tau,h}(t) \|^2 \ud t 
	\; \lesssim \; 
	\tau_n \| p_\tau (t_n) - \widetilde p_{\tau,h}(t_n) \|^2 
	+
	\tau_n \| p_\tau (\tnb) - \widetilde p_{\tau,h}(\tnb) \|^2 \,.
	\]
	We sum up from $n=1$ to $N$,  
	use   
	$\, p_\tau (\tnb) - \widetilde p_{\tau,h}(\tnb) = e_p(\tnb) \,$,
	which follows from \eqref{Eq:ErrPPP_03}, 
	apply  \eqref{Eq:ErrPPP_07} and Theorem~\ref{Thm:EstP}, and get 
	\begin{equation}
		\label{Eq:ErrPPP_08}
		\left(\int_0^T 	\| p_\tau(t) - \widetilde p_{\tau,h}(t) \|^2 \ud t\right)^{1/2} \lesssim
                 \|\nabla \euth(0)\| + \tau \|  \nabla \vec E_{\tau,h}^{\vec u} (0)\| + \tau^2 + h^r\,.
	\end{equation}
	Using  estimate  \eqref{Eq:ErrPPP_05} for $p_\tau$,
	we obtain
	\[
	\int_{I_n} 	\| p(t) - p_\tau(t) \|^2 \ud t 
	\; \lesssim \; \tau_n \tau_n^4\,,
	\qquad n=1, \dots, N
	\,.
	\]
	By summation from $n=1$ to $N$ and applying the triangle inequality
	together with \eqref{Eq:ErrPPP_08},
	we finally get the assertion of the theorem. 
\end{mproof}

\begin{rem}
	\label{Rem:est_dt_utilde}
Similarly to Theorem~\ref{Thm:ErrPPP}, we deduce from Lemma~\ref{Lem:DiscErrp}
an error estimate for the time derivative of the post-processed approximation of the velocity field. In particular, we have  
\begin{equation}
	\label{Eq:FErrPPu}
	\left(\int_0^T 	\| \partial_t \vec u (t) - \partial_t \vec{\widetilde u}_{\tau,h} (t) \|^2 \ud t\right)^{1/2} \lesssim
 \|\nabla \euth(0)\| + \tau \|  \nabla \vec E_{\tau,h}^{\vec u} (0)\| + \tau^2 + h^r\,.
\end{equation}
\end{rem}

The discrete initial error norms $\|\nabla \euth(0)\|$ and $\|
\nabla \vec E_{\tau,h}^{\vec u} (0)\|$ in \eqref{Eq:ErrPPP_00} and
\eqref{Eq:FErrPPu} will be estimated below in Section~\ref{Sec:DisIniVal}. 
For this, the choice of the discrete initial value $\vec u_{0,h}\in \vec
V_h^{\operatorname{div}}$ in Problem~\ref{Prob:FdSAF} is addressed
first. An optimal order approximation of the pressure trajectory is
finally proved by \eqref{Eq:ErrPPP_Order}.

\subsection{Approximation of $\partial_t \texorpdfstring{\vec u(0)}{u(0)}$
and $p(0)$}	
\label{Sec:DisIniVal}
We start with some observation about the solution $(u,p)$ of~\eqref{Eq:SE}.
Provided  $(u,p)$ of~\eqref{Eq:SE} is sufficiently smooth, we have
\begin{equation}
\label{Eq:Collt0}
\vec{a}_0+ \nabla p_0 = \vec{f}_0 + \Delta \vec{u}_0\,,
\end{equation}
with the initial acceleration $\vec{a}_0 \coloneqq \partial_t
\vec{u}(0)$, the initial pressure $p_0 \coloneqq p(0)$, and $\vec{f}_0
\coloneqq \vec{f}(0)$ as right-hand side at the initial time. We note that $\vec{a}_0 \in \vec V^{\operatorname{div}}$ is satisfied.  
For the initial pressure $p_0$, there holds 
\begin{equation}
\label{Eq:atob}
\langle A_h \nabla p_0, \vec{v}_h \rangle = \langle B_h' (-\Delta p_0),
\vec{v}_h \rangle
\qquad \forall \vec{v}_h \in \vec{V}_h\,,
\end{equation}
which follows from integration by parts along with the homogeneous
Dirichlet boundary condition of $\vec{v}_h \in \vec{V}_h$. 

Now, we will construct an approximation $\vec a_{0,h}\in\Vhdiv$ of 
$\vec a_0$. To this end, let
$(\vec a_{0,h},s_h)\coloneqq R_h(\vec a_0,-\Delta p_0)$
with the Stokes projection $R_h$ introduced in
Def.~\ref{Def:StokesProj}, be the unique solution of the problem  
\begin{subequations}
\label{Eq:a0hsh0}
\begin{alignat}{3}
\label{Eq:a0hsh1}
\langle A_h \vec{a}_{0,h} + B_h' s_{h}, \vec{v}_h\rangle
& = 
\langle A_h \vec{a}_0 + B_h' (-\Delta p_0), \vec{v}_h \rangle
& \qquad \forall \vec{v}_h &\in \vec{V}_h,\\[1ex]
\label{Eq:a0hsh2}
\langle B_h \vec{a}_{0,h}, q_h \rangle & = 0 & \qquad \forall q_h & \in Q_h \,.
\end{alignat}
\end{subequations}
Using~\eqref{Eq:Collt0} and~\eqref{Eq:atob} together with the
definitions of the discrete operator $A_h$ and the bilinear form $a$ given
in~\eqref{Def:A_h} and~\eqref{Eq:DefBLFa}, respectively, we get that
$(\vec{a}_{0,h},s_h) \in \vec{V}_h\times Q_h$ is the unique solution of the
Stokes problem
\begin{subequations}
	\label{Eq:a0hsh3}
	\begin{alignat}{3}
		\label{Eq:a0hsh4}
		\langle A_h \vec{a}_{0,h} + B_h' s_{h}, \vec{v}_h\rangle
		& = 
\langle \nabla (\vec{f}_0 + \Delta \vec{u}_0), \nabla \vec{v}_h \rangle
		& \qquad \forall \vec{v}_h &\in \vec{V}_h,\\[1ex]
		\label{Eq:a0hsh5}
		\langle B_h \vec{a}_{0,h}, q_h \rangle & = 0 & \qquad \forall q_h & \in Q_h \,.
	\end{alignat}
\end{subequations}
We note that the right-hand side in \eqref{Eq:a0hsh4} is computable for
given data $\vec f_0,  \Delta \vec{u}_0 \in \vec H^1(\Omega)$. 
From the approximation
property \eqref{Eq:StokesApprox1} of the Stokes projection $(\vec
a_{0,h},s_h)= R_h(\vec a_0,-\Delta p_0)$, we then get the error estimate   
\begin{equation}
\label{Eq:Erra0h}
\| \vec a_0 - \vec a_{0,h}\| \lesssim h^{r+1}\,.
\end{equation}

Using the approximation $\vec a_{0,h}$ of 
$\vec a_0$ computed before, we are now able
to define an approximation $(\vec u_{0,h}, p_{0,h})$ of
$(\vec u_0, p_0)$.
Then, this function $\vec u_{0,h}\in  \vec V_h^{\operatorname{div}}$ 
will be used as the discrete initial velocity in the time marching process. 
\begin{defi}[Discrete initial value $\vec u_{0,h}\in  \vec V_h^{\operatorname{div}}$]
\label{Def:DiscIniVal}
Let $\vec a_{0,h}\in  \vec V_h^{\operatorname{div}}$ be given by
  \eqref{Eq:a0hsh3}. We define $(\vec u_{0,h},p_{0,h})\in \vec
  V_h^{\operatorname{div}}\times Q_h$ as the solution of the Stokes
  system
\begin{subequations}
	\label{Eq:u0hp0h1}
\begin{alignat}{3}
	\label{Eq:u0hp0h2}
\langle A_h \vec{u}_{0,h} + B_h' p_{0,h}, \vec{v}_h\rangle
& = \langle \vec f_0 - \vec a_{0,h}, \vec v_h\rangle
& \qquad \forall \vec{v}_h &\in \vec{V}_h,\\[1ex]
\label{Eq:u0hp0h3}
\langle B_h \vec{u}_{0,h}, q_h \rangle & = 0 & \qquad \forall q_h & \in Q_h\,.
\end{alignat}
\end{subequations}	
\end{defi}

We note that \eqref{Eq:u0hp0h2} corresponds to the collocation condition
\eqref{Eq:ColCondIp_1} at time $t=0$ from the definition of the 
post-processed solution $(\utht^1, \ptht^1)$ on interval $\overline I_1$.
Now, we can
estimate the errors $\|\nabla \vec{e}_{\tau,h}^{\vec{u}}(0)\|$  and
$\|\nabla \vec{E}_{\tau,h}^{\vec{u}}(0)\|$ of \eqref{Eq:ErrPPP_00}.

\begin{lem}
\label{Lem:InitError}
Let the discrete initial velocity $\vec u_{0,h}\in  \vec
  V_h^{\operatorname{div}}$  be defined according to
  Def.~\ref{Def:DiscIniVal}. For the discrete initial errors
  $\vec{e}_{\tau,h}^{\vec{u}}(0)$ and $\vec{E}_{\tau,h}^{\vec{u}}(0)$,
  defined in~\eqref{Eq:SpErr1} and~\eqref{Eq:DefErSplitpu}, we have 
\begin{equation}
\label{Eq:InitError}
\|\nabla \vec{e}_{\tau,h}^{\vec{u}}(0)\| \lesssim h^{r+1}\quad \text{and} \quad 
\quad 
\|\nabla \vec{E}_{\tau,h}^{\vec{u}}(0)\| \lesssim h^{r}\,.
\end{equation}
\end{lem}

\begin{mproof}
From \eqref{Eq:SE} at time $t=0$ we get  	
\begin{subequations}
	\label{Eq:SEt01}
	\begin{alignat}{3}
		\label{Eq:SEt02}
		\langle A_h \vec{u}_{0} + B_h' p_{0}, \vec{v}_h\rangle
            & = \langle \vec f_0 - \vec a_0, \vec v_h\rangle
		& \qquad \forall \vec{v}_h &\in \vec{V}_h,\\[1ex]
		\label{Eq:SEt03}
		\langle B_h \vec{u}_{0,}, q_h \rangle & = 0 & \qquad \forall q_h & \in Q_h \,.
	\end{alignat}
\end{subequations}	
Using \eqref{Eq:SEt01}, we get that the Stokes projection $ (\widehat{\vec{u}}_{h},\widehat p_h) = R_h(\vec u_0,p_0)\in \vec V_h^{\operatorname{div}}\times Q_h$ of $(\vec u_0,p_0)$, according to Def.~\ref{Def:StokesProj}, satisfies the system
\begin{subequations}
	\label{Eq:StProjt01}
	\begin{alignat}{3}
		\label{Eq:StProjt01a}
            \langle A_h \widehat{\vec{u}}_{h} + B_h' \widehat p_h, \vec{v}_h\rangle
            & = \langle \vec f_0 - \vec a_{0}, \vec v_h\rangle
		& \qquad \forall \vec{v}_h &\in \vec{V}_h,\\[1ex]
		\label{Eq:StProjt01b}
		\langle B_h \widehat{\vec{u}}_{h}, q_h \rangle & = 0 & \qquad \forall q_h & \in Q_h\,.
	\end{alignat}
\end{subequations}	
Subtracting \eqref{Eq:u0hp0h1} from \eqref{Eq:StProjt01} and using the
  notation the notation $\heuh \coloneqq \widehat{\vec{u}}_{h} - \vec
  u_{0,h}$ and $\heph\coloneqq \widehat{p}_{h} - p_{0,h}$
  we obtain 
\begin{subequations}
	\label{Eq:StProjt02}
	\begin{alignat}{3}
		\label{Eq:StProjt02a}
		\langle A_h \heuh + B_h' \heph, \vec{v}_h\rangle
		& = \langle \vec a_{0,h} - \vec a_{0}, \vec v_h\rangle
		& \qquad \forall \vec{v}_h &\in \vec{V}_h,\\[1ex]
		\label{Eq:StProjt02b}
		\langle B_h  \heuh, q_h \rangle & = 0 & \qquad \forall q_h & \in Q_h\,.
	\end{alignat}
\end{subequations}	
Choosing the test function $\vec v_h = \heuh \in \vec V_h^{\operatorname{div}}$ 
in \eqref{Eq:StProjt02a}, it follows
\begin{equation*}
\|  \heuh \|_1^2 =  \langle A_h \heuh , \heuh \rangle = \langle \vec a_{0,h} - \vec a_{0}, \heuh \rangle \,.
\end{equation*} 
Together with \eqref{Eq:Erra0h} and the inequality of Cauchy--Schwarz, this implies 
\begin{equation*}
\|\nabla  \heuh \| \lesssim h^{r+1}\,.
\end{equation*} 
  Using the notation introduced in \eqref{Eq:SpErr1}, we get 
\[
  \euth(0) \coloneqq \udth(0) - \uth(0)
            = R_h^{\vec u}( \vec u(0), p(0) ) - \vec{u}_{0,h}
            = \heuh \,,
\]
which proves the first assertion in \eqref{Eq:InitError}.
To prove the second bound in \eqref{Eq:InitError}, {we use}
  the definition~\eqref{Eq:DefErSplitpu} of $\vec{E}_{\tau,h}^{\vec{u}}$,
  {yielding that} 
\begin{align}
\label{eq:InitError:a}
\vec{E}_{\tau,h}^{\vec{u}}(0)
& = \widehat{\vec{U}}_{\tau,h}(0)
- \partial_t \widetilde{\vec{u}}_{\tau,h}(0)
= R_h^{\vec{u}}\big( (R_\tau^1 \partial_t \vec{u})(0),
  (R_\tau^1 \partial_t p)(0)\big) -\vec{a}_{0,h} 
\nonumber\\[1ex]
& = R_h^{\vec{u}}\big(\vec{a}_0, (\partial_t p)(0)\big) - \vec{a}_{0,h}\,.
\nonumber
\end{align}
Hence, error estimate~\eqref{Eq:StokesApprox1} for the velocity part of
the Stokes projection of Def.~\ref{Def:StokesProj} directly proves the
second assertion in \eqref{Eq:InitError}.
\end{mproof}

Now, Lemma~\ref{Lem:InitError} and Theorem~\ref{Thm:ErrPPP}
result in the following corollary. 

\begin{cor}[Error estimate for the post-processed pressure]
	\label{Cor:ErrPPP}
	Let $\vec u_{0,h}\in \vec V_h^{\operatorname{div}}$ be given by Def.~\ref{Def:DiscIniVal}. 
	Then, for the post-pro\-cessed pressure $\widetilde
	p_{\tau,h}$, defined by \eqref{Eq:PP_p} and \eqref{Eq:ColCond},
        the error estimate
	\begin{equation}
		\label{Eq:ErrPPP_Order}
		\left(\int_0^T \| p(t) - \widetilde p_{\tau,h}(t) \|^2 \ud t\right)^{1/2} 
		\lesssim  \tau^2 + h^r
	\end{equation}
	is satisfied. 
\end{cor}


\section{Post-processing of the pressure by interpolation}
\label{Sec:PPInt}

In this section we propose a post-processing of the pressure values
$\pnbh = \pth(\tnb)$ from the solution of Problem~\ref{Prob:FdSAF}  by a
computationally cheap combined interpolation and extrapolation
technique, again with the aim of defining a piecewise polynomial pressure trajectory
$\overline p_{\tau,h}: \overline I \to Q_h$ with optimal order
approximation properties with respect to the mesh sizes in time and
space. The trajectory $\overline p_{\tau,h}$ computed
here admits discontinuities in time at the discrete time nodes $t_n$ for
$n=2,\ldots,N-1$. Its construction principle by interpolation exploits
the optimal order approximation properties of the pressure solution
$\pth$ of Problem~\ref{Prob:FdS} in the sense of
Lemma~\ref{Lem:space_time_probl} in the midpoints $\tnb$ of the time
intervals $I_n$ (cf.\ Theorem~\ref{Thm:EstP}). In contrast to this numerical
type of construction, the post-processed quantities $(\utht,\ptht)$ of
Section~\ref{Sec:PPCol} are derived by collocation techniques for the model
equations based on physical principles.

Here, we  construct $\pthb$  such that 
it belongs to the set of the pressure solutions described in 
Lemma~\ref{Lem:space_time_probl}, i.e., that $(\uth,\pthb)$ is also a 
solution of Problem~\ref{Eq:FdS1} on the time
intervals $I_n$. To this end, we have to guarantee that 
$\pthb(\tnb) = \pth(\tnb)$ for all $n=1, \dots, N$.
The second degree of freedom to define the linear polynomial
$\pthb{}_{|I_n}\in\P_1(I_n; Q_h)$ is the 
coefficient $p_h^{n-1} \coloneqq  \pthb^+(t_{n-1})\in Q_h$.
For $n\ge 2$, we compute the value $p_h^{n-1}$ by means of linear 
interpolation between the values $\overline p_h^{n-1} \coloneqq  \pth(\overline t_{n-1})$ 
from $I_{n-1}$ and
$\overline p_h^n \coloneqq  \pth(\overline t_n)$ from $I_n$. 
For $n=1$, we compute  $p_h^{0}$ by means of linear extrapolation with
the values $\pth(\overline t_1)$ and $\pth(\overline t_2)$.  This is
sketched in Figure~\ref{Fig:PPI}. We note that for all different original pressure solutions $\pth$ in the
set of pressure solutions (cf.\ Lemma~\ref{Lem:space_time_probl}), the values $\overline p_h^n \coloneqq  \pth(\tnb)$ are
unique. This implies that also our post-processed pressure
trajectory $\pthb$, constructed from these values, 
will be unique for all different original pressure
solutions $\pth$. 

\begin{figure}[t]	   
	\centering	
		\begin{tikzpicture}[x=1cm, y=1cm, domain=0:12,smooth]
			\draw[->,thick] (0,0.5) -- (12.5,0.5) node[below] {$t$};
			\draw[->,thick] (0.5,0) -- (0.5,4.5) node[left] {};
			
			\foreach \c in {1.5,2.5,...,12}{\draw (\c,.4) -- (\c,.6);
				
				\node at (0.2,0.2){$t_0$};
				\node at (1.5,0.2){\color{blue} $\overline t_1$};
				\node at (2.5,0.2){$t_1$};
				\node at (3.5,0.2){\color{blue} $\overline t_2$};
				\node at (4.5,0.2){\color{blue} $t_2$};
				
				\node at (6.5,0.2){$t_{n-1}$};
				\node at (7.5,0.2){\color{blue} $\overline t_n$};
				\node at (8.5,0.2){$t_{n}$};
				\node at (9.5,0.2){\color{blue} $\overline t_{n+1}$};
				\node at (10.5,0.2){$t_{n+1}$};
				\node at (11.5,0.2){$\overline t_{n+2}$};
				
				\draw[name path=p1] plot (\x, {-0.09*pow(\x-5.5,2)+4.12});  
				
				\draw[red, thick] (0.5,2.02) -- (4.5,4.18);
				\draw[red, thick] (6.5,3.81) -- (8.5,3.47);
				\draw[red, thick] (8.5,3.10) -- (10.5,2.02);
				\draw[red, thick] (10.5,1.66) -- (12.0,0.31);

				\fill[blue] (1.5,2.56) circle (0.1);
				\fill[blue] (3.5,3.64) circle (0.1);
				\fill[blue] (7.5,3.64) circle (0.1);
				\fill[blue] (9.5,2.56) circle (0.1);
				\fill[blue] (11.5,0.76) circle (0.1);
				\fill[black] (4.5,4.18) circle (0.1);
				\fill[black] (8.5,3.47) circle (0.1);
				\fill[black] (10.5,2.02) circle (0.1);
				
			    \draw[blue, dashed] (9.5,2.56) -- (10.5,1.66);
				\draw[blue, dashed] (7.5,3.64)  -- (8.5,3.10) ;

				\draw[red, dashed] (2.5,0.5) -- (2.5,3.10);
				\draw[red, dashed] (6.5,0.5) -- (6.5,3.81);
				\draw[red, dashed] (8.5,0.5) -- (8.5,3.10);
				\draw[red, dashed] (10.5,0.5) -- (10.5,1.66);				
				
				\draw[red] (0.5,2.02) circle (0.1);
				\draw[red] (2.5,3.1) circle (0.1);
				\draw[red] (6.5,3.81) circle (0.1);
				\draw[red] (8.5,3.10) circle (0.1);
				\draw[red] (10.5,1.66) circle (0.1);
				
				\node at    (0.1,1.7)   {$p$};
				\node[blue] at  (1.5,2.13)  {$\overline p_h^1$};
				\node[blue] at  (3.7,3.30) {$\overline p_h^2$};
				
				\node[blue] at  (7.5,3.2) {$\overline p_h^{n}$};
				\node[blue] at  (9.2,2.1) {$\overline p_h^{n+1}$};
				\node[blue] at  (10.9,0.75)   {$\overline p_h^{n+2}$};
				
				\node[red] at      (0.8,1.7)   {$p_h^{0}$};
				\node[red] at      (2.8,2.8)   {$p_h^{1}$};
				
				\node[red] at     (6.1,3.45)  {$p_h^{n-1}$};
				\node[red] at     (8.2,2.8)  {$p_h^{n}$};
				\node[red] at    (10.1,1.3)   {$p_h^{n+1}$};

                 \node[black] at     (5.0,4.5)  {$\overline p_{\tau,h}(t_{2})$};
                 \node[black] at     (9.0,3.80)  {$\overline p_{\tau,h}(t_{n})$};
                 \node[black] at    (11.2,2.40)   {$\overline p_{\tau,h}(t_{n+1})$};
			}
		\end{tikzpicture}
      \caption{Construction \eqref{Def:IntPres1} of 
      ${\color{red} \overline p_{\tau,h}}$ 
      on ${\color{blue} I_{n+1}}$ by
      interpolation of ${\color{blue} \overline p_h^{n}}$ and
      ${\color{blue} \overline p_h^{n+1}}$.}
	\label{Fig:PPI}
\end{figure} 

\begin{defi}[Local interpolation operator $J_n$]
	\label{Def:phnm1_0}
	For a Banach space $B$ and the time mesh  $\mathcal M_\tau$, we
	define
	\begin{equation}
		\label{Eq:DefCmtau}
		C^{-1}_\tau (I; B) \coloneqq  
		\left\{{q : I \to B} : q_{| I_n} \in C^{0}(I_n;B), \; 
		q^+(t_{n-1}) \;\text{exists} \;\; \forall I_n\in \mathcal M_\tau \right\}\,,
	\end{equation}
      where  $q^+(t_{n-1})$ denotes the right-sided limit at the mesh
      point $t_{n-1}$ (cf.\ \eqref{Eq:Defw_In_bdr}). Now, for  $n\ge 2$,
      we define the
	local interpolation operator 
      $J_n : {C^{-1}_\tau (I; B)} \to \P_1( [t_{n-2},t_n]; B)$
	by means of linear interpolation of the values $q(\tnbm)$ and
	$q(\tnb)$, such that 
	\begin{equation}
		\label{Def:Jnop}
		J_n q (t) \coloneqq 
		q(\tnbm) + \frac{t-\tnbm}{\tnb -\tnbm} \big( q(\tnb) - q(\tnbm) \big),
		\qquad t\in [t_{n-2},t_n] \,.
	\end{equation}
\end{defi}
We note that $J_n$ is a linear operator. For $B=Q$, we can easily show 
\begin{equation}
	\label{Def:Jnest}
	\| J_n q (t) \| \lesssim \| q(\tnbm) \|  + \| q(\tnb) \|,
	\qquad t\in [t_{n-2},t_n] \,,\; n\ge 2 \,.
\end{equation}
From standard arguments of interpolation theory we get for a smooth function $q\in C^2(\overline I; Q)$ that
\begin{equation}
	\label{Def:Jnapprox}
	\| q(t) - J_n q (t) \| \lesssim \tau^2 \|q\|_{C^2(\overline I; Q)},
	\qquad t\in [t_{n-2},t_n] \,,\; n\ge 2 \,.
\end{equation}
Using the operators $J_n$, for $n=2,\dots,N$, with $B=Q_h$,  we can now
define the  post-processed and computationally cheap pressure trajectory $\pthb: \overline I \to Q_h $. 

\begin{defi}[Post-processed pressure trajectory $\overline p_{\tau,h}$]
	\label{Def:IntPres0}
      For a given solution
	$\pth \in Y_\tau^1(Q_h) \subset C^{-1}_\tau (I; Q_h)$ 
	from the set of the pressure solutions   
    according to Lemma~\ref{Lem:space_time_probl}, we define the
	post-processed pressure trajectory $\pthb: \overline I \to Q_h $ by
	\begin{equation}
		\label{Def:IntPres1}
		\pthb(t) \coloneqq
		\begin{cases}
			 J_2 \pth(t),  & t\in \overline I_1\,,  \\[1ex]
			 J_n \pth(t),  & t\in I_n\,, n\ge 2 \,. 
		\end{cases}
	\end{equation}
\end{defi}
The construction of $\pthb$ is sketched in Figure~\ref{Fig:PPI}. By the definition of $J_n \pth$ in \eqref{Def:Jnop} we have   $\pthb(\tnb) = \pth(\tnb) = \overline p_h^n  $ for all $n=1, \dots, N$,
which shows that the post-processed pressure $\pthb$ also belongs to the
set of the pressure solutions (cf.\ Lemma~\ref{Lem:space_time_probl}).

In the following, we will estimate the error of the post-processed pressure
\begin{equation}
	\label{Eq:ovlep}
	\epb(t) \coloneqq  p(t)- \pthb(t)	\,.
\end{equation}
In order to split $\epb(t)$, we define the time-dependent index function $n: \overline I\to \N$ by
\begin{equation}
	\label{Def:nt}
	n(t) \coloneqq
	\begin{cases}
		2,   & t\in \overline I_1\,,  \\[1mm]
   		m,   & t\in I_m\,, \; m\ge 2 \,. 
	\end{cases}
\end{equation}
For fixed $t\in\overline I$ and with the abbreviation $n=n(t)$, 
the value $J_{n}q(t)$ is well defined for any function $q\in
C^{-1}_\tau (\overline I; B)$ by~\eqref{Def:Jnop} since $t\in [t_{n-2}, t_n]$
for $n=n(t)$.  By the index function \eqref{Def:nt}, the
definition~\eqref{Def:IntPres1} of $\pthb$ can be written shortly as 
  $\pthb(t) \coloneqq  J_{n(t)} \pth(t) $ for each $t\in\overline I$. For
fixed $t\in \overline I$, we split the pressure error $\epb(t)$
according to
\begin{equation}
	\label{Eq:ovleps}
	\epb(t) = \overline \eta^p(t) + \eptb(t)   \,, \qquad
	\overline \eta^p(t) \coloneqq p(t) - J_{n(t)} p(t)\,,\quad
	\eptb(t)\coloneqq J_{n(t)} p(t) - J_{n(t)} \pth(t) \,.
\end{equation}
For the interpolation error $\overline \eta^p$,
estimate~\eqref{Def:Jnapprox} implies 
\begin{equation}
	\label{Eq:oetap}
	\| \overline \eta^p (t)\| = \| p(t) - J_{n(t)} p(t) \|
	\lesssim \tau^2,   \qquad t\in\overline I \,.
\end{equation}

\begin{lem}[Estimate of $\eptb$]
	\label{Eq:oethp}
      Suppose that the time step sizes $\tau_m$, $m=1, \dots, N$,
      satisfy the conditions
	\begin{equation}
		\label{Eq:CondTau}
		\tau_{1} \leq  c_1 \tau_2  \qquad\text{and}\qquad
		\tau_m \leq c_2 \tau_{m-1}\,, \quad m=2,\ldots,N\,,
	\end{equation}
	with some fixed positive constants $c_1$ and $c_2$. 
	Then, for the error $\eptb$ defined in~\eqref{Eq:ovleps}, it holds 
	\begin{equation}
		\label{Eq:oethp1}
		\sum_{m=1}^N \tau_m \| \eptb(t_m)\|^2  \lesssim \| \nabla \vec e_{\tau,h}^{\vec u}(0)\|^2  + \tau^4 + h^{2r}\,.
	\end{equation}
\end{lem}

\begin{mproof}
  For any fixed time $t=t_m \in[t_{n-2},t_n]$ with $m\in \{1,\dots,N\}$,
we use the abbreviation $n=n(t_m)$ and get for the error $\eptb$
defined in~\eqref{Eq:ovleps} that
	\[
	\eptb(t_m) = J_n p(t_m) - J_n \pth(t_m) = J_n(p-\pth)(t_m) =
	J_n q(t_m) \,,
	\]
where $q \coloneqq  p - \pth \in C^{-1}_\tau (I; Q)$.
  Then, the application of \eqref{Def:Jnest} for $t=t_m$ implies 
	\begin{equation}
		\label{Def:IntPres3}
		\| \eptb(t_m)\|  \lesssim  \| (p-\pth)(\overline
            t_{n-1}) \| +  \| (p-\pth)(\overline t_n)\| \,. 
	\end{equation}
We recall the notation $e_p(t) \coloneqq  (p-\pth)(t)$ from Theorem~\ref{Thm:EstP}. 
In the case $m\ge 2$, we have  $n=n(t_m)=m$ such that~\eqref{Def:IntPres3}
yields the estimate
	\begin{equation}
		\label{Def:IntPres4}
		\| \eptb(t_m)\|  \lesssim  \| e_p(\overline t_{m-1}) \| +  \| e_p(\overline t_m) \|  \,,
		\qquad m=2, \dots, N  \,.
	\end{equation}
In the case $m=1$, we have  $n=n(t_1)=2$ and get from~\eqref{Def:IntPres3}
  that
	\begin{equation}
		\label{Def:IntPres5}
		\| \eptb(t_1)\|  \lesssim  \| e_p(\overline t_1) \| +  \| e_p(\overline t_2) \|  \,.
	\end{equation}
	Combining the assumptions made in \eqref{Eq:CondTau} with \eqref{Def:IntPres4} and \eqref{Def:IntPres5},
	we  deduce the estimate 
	\begin{equation}	
		\label{Eq:oethp3}
		\sum_{m=1}^N \tau_m \| \eptb(t_m)\|^2  \lesssim 
		\sum_{m=1}^N \tau_m \| e_p(\overline t_m)\|^2  \,.
	\end{equation}
	From \eqref{Eq:oethp3} along with the results of Theorem~\ref{Thm:EstP} we directly conclude 
	the assertion~\eqref{Eq:oethp1}. 
\end{mproof}

{We note that condition~\eqref{Eq:CondTau} excludes an extremely rapid
change of the time step size of neighboring time intervals. }

\begin{thm}
	\label{Thm:EPPI}
	Suppose that condition \eqref{Eq:CondTau} is satisfied. Then, the
	post-processed pressure $\overline p_{\tau,h}$ of
	Def.~\ref{Def:IntPres0} and the right-sided limits $\pthb^+(t_{m-1})$ satisfy 
	\begin{subequations}	 
		\label{Eq:pbest0}
		\begin{alignat}{2}
		\label{Eq:pbest1}
		\sum_{m=1}^N \tau_m \| p(t_{m-1}) - \overline p_{\tau,h}^+(t_{m-1})\|^2  
	    & \lesssim   \| \nabla \vec e_{\tau,h}^{\vec u}(0)\|^2  + \tau^4 + h^{2r} \,, \\[1ex]
		\label{Eq:pbest2}
		\left( \int_0^T  \| p(t) - \overline p_{\tau,h}(t)\|^2 \ud t \right)^{1/2}
		& \lesssim  \| \nabla \vec e_{\tau,h}^{\vec u}(0)\|  + \tau^2 + h^{r}\,.
	\end{alignat}
	\end{subequations}
\end{thm}

\begin{mproof}
  For all $t\in I_m = (t_{m-1},t_m]$, we have from~\eqref{Def:nt} that
	$n(t)=2$ if $m=1$ and $n(t)=m$ if $m\ge 2$. 
      Therefore, the index function $n(t)$ has the constant value
      $n(t_m) \ge 2$ for all $t\in I_m$,
	which implies for $\eptb$ in the splitting~\eqref{Eq:ovleps}
	that 
	\begin{equation}	 
		\label{Eq:EPPI_1}
		\eptb(t) = J_n e_p(t)   \,, 
            \qquad n=n(t_m) \,,  \quad t\in I_m  \,,
	\end{equation}
	where $e_p \coloneqq  p - \pth$.
      Since $I_m\subset [t_{n-2}, t_n]$, where $n=n(t_m)$, we get from the 
      definition of $J_n$ in~\eqref{Def:Jnop} that  $\eptb{}_{| I_m} \in\P_1(I_m; Q)$.
	Thus, by interpolation arguments we can easily show  
	\begin{equation}	 
		\label{Eq:EPPI_2}
		\| \eptb(t) \|  \lesssim  \| \eptb(\overline t_m) \| + \| \eptb(t_m) \|\,,
		\qquad t\in I_m  \,.
	\end{equation}
	From~\eqref{Eq:EPPI_1} and the definition of $J_n$ we obtain  
	$\eptb(\overline t_m) = e_p(\overline t_m)$.
	Thus, applying~\eqref{Eq:oetap} and the splitting
	$\epb = \overline \eta^p + \eptb $  in~\eqref{Eq:ovleps}, we conclude
	\begin{equation}	 
		\label{Eq:EPPI_3}
		\| \epb(t) \|^2 \lesssim  \| \overline\eta^p(t) \|^2 + \| \eptb(t) \|^2 
		\lesssim  \tau^4 +  \| e_p(\overline t_m) \|^2 + \| \eptb(t_m) \|^2\,,
		\qquad t\in I_m  \,.
	\end{equation}
	If we take in~\eqref{Eq:EPPI_3} the right-sided limit
	for $t\to t_{m-1}+0$, we get 
	\begin{equation}	 
		\label{Eq:EPPI_4}
		\| {\epb}^+(t_{m-1}) \|^2 =  \| p(t_{m-1}) - \overline p_{\tau,h}^+(t_{m-1}) \|^2 
		\lesssim  \tau^4 +  \| e_p(\overline t_m) \|^2 + \| \eptb(t_m) \|^2 \,,
	\end{equation}
	where we have assumed global continuity of the pressure trajectory $p(t)$.
	Now, in order to get the assertion~\eqref{Eq:pbest1},
	we multiply~\eqref{Eq:EPPI_4} by $\tau_m$, take the sum over all
	$m$ from $1$ to $N$, and apply Theorem~\ref{Thm:EstP} together with
	estimate~\eqref{Eq:oethp1} of Lemma~\ref{Eq:oethp}.
	If we integrate \eqref{Eq:EPPI_3} over $t\in I_m$, we get 
	\begin{equation}	 
		\label{Eq:EPPI_5}
		\int_{I_m} \| \epb(t) \|^2 \ud t 
            = \int_{I_m} \| p(t) - \overline p_{\tau,h}(t) \|^2 \ud t  
		\lesssim  
		\tau_m \tau^4 +  \tau_m \| e_p(\overline t_m) \|^2 + \tau_m \| \eptb(t_m) \|^2 \,.
	\end{equation}
	Finally, assertion~\eqref{Eq:pbest2} follows by summation of~\eqref{Eq:EPPI_5}
	over all $m$ from $1$ to $N$ and applying Theorem~\ref{Thm:EstP} and 
	Lemma~\ref{Eq:oethp}.
\end{mproof}

\begin{rem}
We will comment on the results above.
	\label{Rem:PP}
	\begin{itemize}
        \item[a)]  In contrast to the approach of Section\ \ref{Sec:PPCol}
	based on post-processing by collocation, we do not need here an
	approximation $\vec a_{0,h} \in\Vhdiv$  of $\vec a_0 = \dt \vec u(0)$
	in order to
	determine the initial values $\vec u_{0,h}$  and   $p_{0,h}$  (cf.\
	Section\ \ref{Sec:DisIniVal}). It is sufficient to determine  $\vec u_{0,h}
	\in\Vhdiv$  such that $\| \nabla (\vec u_0 - \vec u_{0,h}) \| \lesssim h^r$.
	In particular, this is satisfied for the choice  $\vec u_{0,h} =
	R^{\vec u}_h(\vec u_0, 0)$ which corresponds to the solution of a Stokes
	problem. Then, for the resulting velocity error at $t=0$ in 
	estimates \eqref{Eq:pbest0}, we get  $  \| \nabla \vec e^{\vec u}_{\tau,h}(0) \|  \lesssim h^r $ .
        
        \item[b)]  Theorem~\ref{Thm:EPPI} shows the optimal order approximation of the
		post-processed pressure $\overline p_{\tau,h}$ in the norm of $L^2(I;Q)$
		and in a discrete $L^2$-norm in time. 
		These error norms will be also used in the next section to measure 
		the pressure error in our numerical experiments.
		The post-processed pressure trajectory 
               $\overline p_{\tau,h}$ 
		is in general discontinuous at the discrete time nodes $t_n$ 
            (cf. Figure~\ref{Fig:PPI}). 
		From estimate~\eqref{Eq:oethp1} together with the assumption
		$\tau_{m} \le c_2 \tau_{m-1}$,  $m=2,\dots,N$, and  estimate~\eqref{Eq:oetap} for
		$\|\overline\eta^p(t_{m-1})\|$, we get  
		\[
		\sum_{m=2}^N \tau_m \| p(t_{m-1}) - \overline p_{\tau,h}(t_{m-1})\|^2  
		\lesssim  \| \nabla \vec e_{\tau,h}^{\vec u}(0)\|^2   + \tau^4 + h^{2r} \,.
		\]
		Together with~\eqref{Eq:pbest1}
		this shows that the jumps in the discontinuous pressure $\overline p_{\tau,h}$ 
		at the discrete time nodes are small in the sense that
		\[
		\sum_{m=2}^N \tau_m 
		\| \overline p_{\tau,h}^+(t_{m-1}) - \overline p_{\tau,h}(t_{m-1})\|^2  
		\lesssim  \| \nabla \vec e_{\tau,h}^{\vec u}(0)\|^2  + \tau^4 + h^{2r} \,.
		\]
		We note that the jump at $t=0$ is zero since, by definition of $\pthb$, we have  $\pthb(0)=J_2 \pth(0)=\pthb^+(0)$.
		
        \item[c)] Alternatively, instead of computing $\overline p_{\tau,h}$ we can define a continuous pressure trajectory $\widehat p_{\tau,h} \in C(\overline I; Q_h)$ such that  
        $\widehat p_{\tau,h}{}_{|I_n} \in \P_1(I_n, Q_h)$  for $n=1, \dots, N$  and
        $\widehat p_{\tau,h}(t_n) = \overline p_{\tau,h}(t_n)$ for $n=0, \dots, N$. However, we note that the pair $(\uth, \widehat p_{\tau,h})$ will in general not satisfy the variational problem~\eqref{Eq:FdS1} on the subinterval $I_n$.
		On the other hand, we directly get from Theorem~\ref{Thm:EPPI}
		the corresponding error estimates also for the post-processed pressure 
		$\widehat p_{\tau,h}$ instead of $\pthb$.
	\end{itemize}
\end{rem}

\section{Numerical experiments}
\label{Sec:NumExp}

In this section we illustrate by numerical experiments the error
estimates proved in Theorem~\ref{Thm:EstP} for the pressure approximation
in the Gauss quadrature nodes $\overline t_n$, in Theorem\ \ref{Thm:ErrPPP}
for post-processing of the pressure variable by collocation and,
finally,  in Theorem\ \ref{Thm:EPPI} for post-processing of the pressure
variable by interpolation.
The required discrete initial quantities  at time $t=0$ were chosen in such
a way (see below) that the resulting initial error terms in the
corresponding error estimates are bounded as in 
\eqref{Eq:InitError} of Lemma~\ref{Lem:InitError}.
The implementation of the numerical scheme was
done in an in-house high-performance frontend solver \cite{AB21} for the
\texttt{deal.II} library \cite{Aetal21}. 

In the numerical experiments, the pressure error $p-p_{\tau,h}$ is
measured in the natural norm of $L^2(I; B)$, shortly denoted by $\|\cdot
\|_{L^2(B)}$, where $B= L^2(\Omega)$. Furthermore, the
velocity errors $\vec e_{\vec u} \coloneqq \vec u - \vec u_{\tau,h}$ and
$\dt\vec e_{\vec u}$ are also computed with respect to the norm
of $L^2(I;B)$, with $B=\vec H^1(\Omega)$ and $B=\vec L^2(\Omega)$,
respectively (see Theorem~\ref{Thm:Erruprime}). 
For the computation of the $L^2(I;B)$ 
norm of the error, a quadrature formula of high accuracy is used for
time integration.
Moreover,  errors are measured in the time-mesh
dependent norms (cf.~\eqref{Eq:EstP} and \eqref{Eq:pbest1})
\begin{alignat}{2}
\label{Eq:l2ErrNrms}
\| w \|_{\overline{l^2} (B)} & \coloneqq \bigg( \sum_{m=1}^N \tau_m \| w(\overline t_m) \|_B^2 \bigg)^{1/2}\,,
\qquad
\| w \|_{l^2(B)} & \coloneqq \bigg( \sum_{m=1}^N \tau_m \| w^+(t_{m-1}) \|_B^2 \bigg)^{1/2}\,,
\end{alignat}
defined for any function $w:\overline I \to B$ that is piecewise
sufficiently smooth with respect to the time mesh $\mathcal M_\tau$
such that the one-sided limits $w^+(t_{m-1})$ exist
(cf.~\eqref{Eq:Defw_In_bdr}).

In our numerical experiments, we study the Stokes system \eqref{Eq:SE}
for $\Omega=(0,1)^2$, $I=(0,T]$ with $T=2$, and the prescribed
analytical solution
\begin{align*}
	\vec{u}(\vec{x},t)
	&\coloneqq 
	\begin{pmatrix}
		\cos(x_2 \pi) \cdot \sin(x_1 \pi)^2  \cdot \sin(x_2 \pi) \cdot \sin(t) \\
		- \cos(x_1 \pi) \cdot \sin(x_2 \pi)^2  \cdot \sin(x_1 \pi) \cdot \sin(t)
	\end{pmatrix}\,,
	\\
	p(\vec{x},t)
	&\coloneqq 
	\cos(x_2 \pi) \cdot \sin(x_1 \pi) \cdot \cos(x_1 \pi) \cdot \sin(x_2 \pi) \cdot \sin(t)\,.
\end{align*}
In the convergence study, the domain $\Omega$ is decomposed into a
sequence of successively refined meshes of quadrilaterals. The spatial and temporal mesh sizes are halved in each of the
refinement steps. The step sizes of the coarsest space and time mesh are
$h_0=0.25\sqrt{2}$ and $\tau_0=1$. For the spatial discretization
of velocity and pressure, we use the
finite element pair $(\mathbb Q_2, \mathbb Q_1)$ of the Taylor--Hood family
which is inf-sup stable such that
Assumption~\ref{Ass:FES} is satisfied with $r=2$.

The calculated errors for the velocity and pressure approximation are
summarized in Tables \ref{Tab:Convergence_Lifting} and
\ref{Tab:Convergence_Interpolation}. For the post-processing by collocation,
the error is computed for the post-processed quantities
$\{\vec{\widetilde u}_{\tau,h},{\widetilde p}_{\tau,h}\}$ with their
restrictions to the interval $I_n$ being defined by \eqref{Eq:PP_up}.
For post-processing by interpolation, the error is computed for the
post-processed quantities $\{\vec{u}_{\tau,h},{\overline p}_{\tau,h}\}$
with ${\overline p}_{\tau,h}$ being defined by \eqref{Def:IntPres1}.

We note that in all numerical computations
the discrete initial velocity was chosen as 
$\vec u_{0,h}=\vec u_0 = \vec 0$
and not as the velocity solution of the Stokes problem~\eqref{Eq:u0hp0h1}
described in Section~\ref{Sec:DisIniVal}.
For the case of post-processing by collocation, we  determine the
post-processed solution  $(\utht,\ptht)$ on the first interval 
$\overline I_1$, as described in the first case of Remark~\ref{Rem:PPcoll_one}, by computing the
coefficients $(\vec{a}_h^{0} , \widetilde{p}_h^{0} )$
as solution
of problem~\eqref{Eq:probl_a_p} with $n=1$ and $\uth(t_0) \coloneqq \vec u_{0,h}$.
Then we can show  that the resulting initial value
$\partial_t\utht(0) = \vec{a}_h^{0}$ satisfies the same 
estimate~\eqref{Eq:InitError} for the error 
$\|\nabla \vec{E}_{\tau,h}^{\vec{u}}(0)\|$
as for the approach presented in Section~\ref{Sec:DisIniVal}. 

The tables nicely confirm the estimates of order $\tau^2 + h^2$
in Theorem~\ref{Thm:EstP} for the
pressure approximation in the Gauss quadrature nodes $\overline t_n$,
in Theorem\ \ref{Thm:ErrPPP} for post-processing by collocation and 
in Theorem\ \ref{Thm:EPPI} for post-processing by interpolation. 
For the approximation of the velocity field, 
the error bounds of order $\tau^2 + h^2$
in Theorem~\ref{Thm:Erruprime} with respect to the norms in \eqref{Eq:l2ErrNrms}
and the bound~\eqref{Eq:FErrPPu} in Remark~\ref{Rem:est_dt_utilde}
are confirmed as well. 
We note that the error of the time derivative of the velocity
$\dt\vec e_{\vec u} \coloneqq \dt \vec u - \dt\uth$
is limited by the order $\tau + h^2$, which is optimal with respect to
the power of $\tau$ due to the piecewise linear approximation in time of 
$\uth$. Table~\ref{Tab:Convergence_Interpolation} confirms this first order
of convergence for $\dt\vec e_{\vec u}$, measured in the $L^2(\vec L^2)$-norm and the
norm $\| \cdot \|_{l^2(B)}$ given in~\eqref{Eq:l2ErrNrms}.
However, we observe second order of convergence for $\dt\vec e_{\vec u}$,
measured in the  time-discrete norm $\| \cdot \|_{\overline{l^2} (B)}$ given in~\eqref{Eq:l2ErrNrms}. This confirms
superconvergence with respect to $\tau$ of $\dt\uth$ in the midpoints of the
time intervals,
which is proved in Theorem~\ref{Thm:Erruprime}.
 
\begin{table}[H]
	\centering
	\footnotesize
	\begin{tabular}{l}
		\begin{tabular}{ccp{2.6cm}cp{3.1cm}cp{2.6cm}c}
			\toprule
			$\tau$ & $h$ &
			 $\| \vec u - \vec{\widetilde u}_{\tau,h}  \|_{L^2(\vec H^1)}  \hspace{0.7em} $  & EOC &
			 $\| \partial_t \vec u - \partial_t \vec{\widetilde u}_{\tau,h}  \|_{L^2(\vec L^2)}  \hspace{0.9em} $ & EOC &
			 $\| p-\widetilde p_{\tau,h}  \|_{L^2(L^2)}   \hspace{0.9em}$ & EOC 
			\\
			\cmidrule(r){1-2}
			\cmidrule(lr){3-8}
$\tau_0/2^0$ & $h_0/2^0$ & 1.5106628370e+00 & -- & 6.4076648649e-01  & -- & 4.0701968487e-02 & -- \\ 
$\tau_0/2^1$ & $h_0/2^1$ & 2.3275917549e-01 & 2.70 & 1.9183003070e-01  & 1.74 & 7.4938877151e-03 & 2.44 \\
$\tau_0/2^2$ & $h_0/2^2$ & 3.2537111184e-02 & 2.84 & 5.1247103473e-02  & 1.90 & 1.8117135173e-03 & 2.05 \\
$\tau_0/2^3$ & $h_0/2^3$ & 5.4564310437e-03 & 2.58 & 1.3176859187e-02  & 1.96 & 4.5193083007e-04 & 2.00 \\
$\tau_0/2^4$ & $h_0/2^4$ & 1.1795859438e-03 & 2.21 & 3.3372082367e-03  & 1.98 & 1.1302787165e-04 & 2.00 \\
$\tau_0/2^5$ & $h_0/2^5$ & 2.8815509699e-04 & 2.03 & 8.3951890439e-04  & 1.99 & 2.8267540304e-05 & 2.00 \\
			\bottomrule
		\end{tabular}\\
		\mbox{}\\
		\begin{tabular}{ccp{2.6cm}cp{3.1cm}cp{2.6cm}c}
			\toprule
			$\tau$ & $h$ &
			 $\| \vec u - \vec{\widetilde u}_{\tau,h}  \|_{\overline{l^{2}}(\vec H^1)} \hspace{0.7em}$ & EOC &
			$\| \partial_t \vec u - \partial_t \vec{\widetilde u}_{\tau,h}  \|_{\overline{l^2}(\vec L^2)} \hspace{0.9em} $ & EOC &
			 $\| p-\widetilde p_{\tau,h}  \|_{\overline{l^2}(L^2)}  \hspace{0.9em}$ & EOC 
			\\
			\cmidrule(r){1-2}
			\cmidrule(lr){3-8}
			$\tau_0/2^0$ & $h_0/2^0$ & 2.0931146740e+00 & -- & 1.0795011327e-02 & -- & 3.3454424913e-02 & --  \\ 
$\tau_0/2^1$ & $h_0/2^1$ & 3.2285623070e-01 & 2.70 & 2.7723520259e-03 & 1.96 & 5.2469222478e-03 & 2.67  \\
$\tau_0/2^2$ & $h_0/2^2$ & 4.3994570438e-02 & 2.88 & 7.1122652548e-04 & 1.96 & 1.1921452162e-03 & 2.14  \\
$\tau_0/2^3$ & $h_0/2^3$ & 6.6982363507e-03 & 2.72 & 1.7864127424e-04 & 1.99 & 2.9261910060e-04 & 2.03  \\
$\tau_0/2^4$ & $h_0/2^4$ & 1.2904255347e-03 & 2.38 & 4.4694121036e-05 & 2.00 & 7.2863123672e-05 & 2.01  \\
$\tau_0/2^5$ & $h_0/2^5$ & 2.9740542344e-04 & 2.12 & 1.1175663260e-05 & 2.00 & 1.8198565535e-05 & 2.00  \\
			\bottomrule
		\end{tabular}\\
		\mbox{}\\
		\begin{tabular}{ccp{2.6cm}cp{3.1cm}cp{2.6cm}c}
			\toprule
			$\tau$ & $h$ &
			 $\| \vec u - \vec{\widetilde u}_{\tau,h}  \|_{l^{2}(\vec H^1)} $ & EOC &
			 $\| \partial_t \vec u - \partial_t \vec{\widetilde u}_{\tau,h}  \|_{l^2(\vec L^2)} $ & EOC &
			 $\| p- \widetilde p_{\tau,h}  \|_{l^2(L^2)}  $ & EOC 
			\\
			\cmidrule(r){1-2}
			\cmidrule(lr){3-8}
			$\tau_0/2^0$ & $h_0/2^0$ & 2.2659376009e-01 & -- & 1.1010666990e+00 & -- & 2.2556399926e-02 & --  \\ 
$\tau_0/2^1$ & $h_0/2^1$ & 6.7342263613e-02 & 1.75 & 3.2870366943e-01 & 1.74 & 8.6470802363e-03 & 1.38  \\
$\tau_0/2^2$ & $h_0/2^2$ & 1.7907425489e-02 & 1.91 & 8.7691721683e-02 & 1.91 & 2.4231261707e-03 & 1.84  \\
$\tau_0/2^3$ & $h_0/2^3$ & 4.5945110759e-03 & 1.96 & 2.2531306281e-02 & 1.96 & 6.2605791082e-04 & 1.95  \\
$\tau_0/2^4$ & $h_0/2^4$ & 1.1624232127e-03 & 1.98 & 5.7042965632e-03 & 1.98 & 1.5806711683e-04 & 1.99  \\
$\tau_0/2^5$ & $h_0/2^5$ & 2.9227716521e-04 & 1.99 & 1.4347418723e-03 & 1.99 & 3.9644243237e-05 & 2.00  \\
			\bottomrule
		\end{tabular}
	\end{tabular}
	\caption{%
		Calculated errors and experimental orders of convergence (EOC) for post-processing \eqref{Eq:PP_up} by collocation, with time-mesh dependent norms \eqref{Eq:l2ErrNrms}.
	}
	\label{Tab:Convergence_Lifting}
\end{table}

\begin{table}[H]
	\centering
	\footnotesize
	\begin{tabular}{l}
		\begin{tabular}{ccp{2.6cm}cp{3.1cm}cp{2.6cm}c}
			\toprule
			$\tau$ & $h$ &
			 $\| \vec u - \vec{u}_{\tau,h}  \|_{L^2(\vec H^1)}  \hspace{0.7em} $ & EOC &
			 $\| \partial_t \vec u - \partial_t \vec{u}_{\tau,h}  \|_{L^2(\vec L^2)}  \hspace{0.9em} $ & EOC &
			 $\| p-\overline p_{\tau,h}  \|_{L^2(L^2)}   \hspace{0.9em}$ & EOC 
			\\
			\cmidrule(r){1-2}
			\cmidrule(lr){3-8}
			$\tau_0/2^0$ & $h_0/2^0$ & 2.4103293998e-01 & -- & 9.8085008523e-02  & -- & 3.8869233248e-02 & -- \\ 
$\tau_0/2^1$ & $h_0/2^1$ & 6.0161939211e-02 & 2.00 & 4.8439406869e-02  & 1.02 & 6.9936716552e-03 & 2.47 \\
$\tau_0/2^2$ & $h_0/2^2$ & 1.5042373657e-02 & 2.00 & 2.4129431471e-02  & 1.01 & 1.6762471579e-03 & 2.06 \\
$\tau_0/2^3$ & $h_0/2^3$ & 3.7608570388e-03 & 2.00 & 1.2052596660e-02  & 1.00 & 4.1881864102e-04 & 2.00 \\
$\tau_0/2^4$ & $h_0/2^4$ & 9.4023407683e-04 & 2.00 & 6.0247593478e-03  & 1.00 & 1.0497742134e-04 & 2.00 \\
$\tau_0/2^5$ & $h_0/2^5$ & 2.3505983620e-04 & 2.00 & 3.0121866831e-03  & 1.00 & 2.6290992294e-05 & 2.00 \\
			\bottomrule
		\end{tabular}\\
		\mbox{}\\
		\begin{tabular}{ccp{2.6cm}cp{3.1cm}cp{2.6cm}c}
			\toprule
			$\tau$ & $h$ &
			 $\| \vec u - \vec{u}_{\tau,h}  \|_{\overline{l^{2}}(\vec H^1)} \hspace{0.7em}$ & EOC &
			 $\| \partial_t \vec u - \partial_t \vec{u}_{\tau,h}  \|_{\overline{l^2}(\vec L^2)} \hspace{0.8em} $ & EOC &
			 $\| p-\overline p_{\tau,h}  \|_{\overline{l^2}(L^2)}  \hspace{0.8em}$ & EOC 
			\\
			\cmidrule(r){1-2}
			\cmidrule(lr){3-8}
			$\tau_0/2^0$ & $h_0/2^0$ & 2.5481474586e-01 & -- & 1.0795011327e-02 & -- & 3.3454424913e-02 & --  \\ 
$\tau_0/2^1$ & $h_0/2^1$ & 6.1982493633e-02 & 2.04 & 2.7723520259e-03 & 1.96 & 5.2469222478e-03 & 2.67  \\
$\tau_0/2^2$ & $h_0/2^2$ & 1.5384702907e-02 & 2.01 & 7.1122652548e-04 & 1.96 & 1.1921452162e-03 & 2.14  \\
$\tau_0/2^3$ & $h_0/2^3$ & 3.8392598434e-03 & 2.00 & 1.7864127424e-04 & 1.99 & 2.9261910060e-04 & 2.03  \\
$\tau_0/2^4$ & $h_0/2^4$ & 9.5938540221e-04 & 2.00 & 4.4694121036e-05 & 2.00 & 7.2863123672e-05 & 2.01  \\
$\tau_0/2^5$ & $h_0/2^5$ & 2.3981953429e-04 & 2.00 & 1.1175663260e-05 & 2.00 & 1.8198565537e-05 & 2.00  \\
			\bottomrule
		\end{tabular}\\
		\mbox{}\\
		\begin{tabular}{ccp{2.6cm}cp{3.1cm}cp{2.6cm}c}
			\toprule
			$\tau$ & $h$ &
			 $\| \vec u - \vec{u}_{\tau,h}  \|_{l^{2}(\vec H^1)} $ & EOC &
			$\| \partial_t \vec u - \partial_t \vec{u}_{\tau,h}  \|_{l^2(\vec L^2)} $ & EOC &
			$\| p- \overline p_{\tau,h}  \|_{l^2(L^2)}  $ & EOC 
			\\
			\cmidrule(r){1-2}
			\cmidrule(lr){3-8}
			$\tau_0/2^0$ & $h_0/2^0$ & 2.2659376009e-01 & -- & 1.4476007458e-01 & -- & 3.9259453944e-02 & --  \\ 
$\tau_0/2^1$ & $h_0/2^1$ & 6.7342263613e-02 & 1.75 & 7.8393217938e-02 & 0.88 & 7.3029860451e-03 & 2.43  \\
$\tau_0/2^2$ & $h_0/2^2$ & 1.7907425489e-02 & 1.91 & 4.0493918949e-02 & 0.95 & 1.7500474215e-03 & 2.06  \\
$\tau_0/2^3$ & $h_0/2^3$ & 4.5945110759e-03 & 1.96 & 2.0558599482e-02 & 0.98 & 4.4166433190e-04 & 1.99  \\
$\tau_0/2^4$ & $h_0/2^4$ & 1.1624232127e-03 & 1.98 & 1.0356767161e-02 & 0.99 & 1.1142530510e-04 & 1.99  \\
$\tau_0/2^5$ & $h_0/2^5$ & 2.9227716521e-04 & 1.99 & 5.1977590232e-03 & 0.99 & 2.8000003707e-05 & 1.99  \\
			\bottomrule
		\end{tabular}
	\end{tabular}
	\caption{%
		Calculated errors and experimental orders of convergence (EOC) for post-processing \eqref{Def:IntPres1} by interpolation, with time-mesh dependent norms \eqref{Eq:l2ErrNrms}.
	}
	\label{Tab:Convergence_Interpolation}
\end{table}

\section*{Acknowledgements}

Computational resources (HPC-cluster HSUper) have been provided by the project hpc.bw, funded by dtec.bw --- Digitalization and Technology Research Center of the Bundeswehr. dtec.bw is funded by the European Union ---  NextGenerationEU.


\end{document}